\documentclass[11pt,a4paper]{article}
\usepackage[left=2.5cm, right=2.5cm, top=3cm, bottom=3cm]{geometry}

\usepackage[english]{babel} %
\usepackage{hyphenat} 
\usepackage{amsfonts,amsmath,amssymb,amsthm}
\usepackage{autobreak}
\usepackage{graphicx} 
\usepackage{subfig}
\usepackage{algorithm}
\usepackage{algorithmicx}%
\usepackage{algpseudocode}%
\usepackage{url}
\usepackage{xcolor}
\usepackage{multirow}
\usepackage[titletoc,title]{appendix}
\usepackage{threeparttable} 
\usepackage{authblk}
\usepackage{hyperref}
\providecommand{\keywords}[1]{\textbf{Keywords:} #1}
\usepackage{booktabs}

\usepackage{pifont}
\usepackage{tabularx}     
\usepackage{arydshln}
\usepackage{rotating}
\usepackage[numbers,sort&compress]{natbib}

\usepackage[right,mathlines]{lineno}
\theoremstyle{thmstyletwo}

\newtheorem{theorem}{Theorem}[section]
\newtheorem{assumption}{Assumption}[section]
\newtheorem{lemma}{Lemma}[section]
\newtheorem{corollary}{Corollary}[section]
\newtheorem{definition}{Definition}[section]
\newtheorem{remark}{Remark}[section]
\numberwithin{equation}{section}
\newcommand{\R}{\mathbb{R}}
\newcommand{\ve}{\varepsilon}%
\newcommand{\B}{\mathbb{B}}

\newcommand{\cL}{\mathcal{L}}

\newcommand{\veg}{\varepsilon_g}
\newcommand{\veh}{\varepsilon_h}

\title{On the complexity of proximal gradient and proximal gradient-Newton-CG methods with negative curvature for \(\ell_1\)-regularized optimization} 
\author{Hong~Zhu\thanks{zhuhongmath@126.com}}
\affil{School of Mathematical Sciences, Jiangsu University, Zhenjiang, 212013, Jiangsu, China.}
\date{}

\begin{document}

\maketitle

\begin{abstract}
In this paper, we propose two second-order methods for solving the \(\ell_1\)-regularized composite optimization problem, which are developed based on two distinct definitions of approximate second-order stationary points. We introduce a hybrid proximal gradient and negative curvature method, as well as an adaptive hybrid proximal gradient-Newton-conjugate gradient (Newton-CG) method with negative curvature directions, to find a strong* approximate second-order stationary point and a weak approximate second-order stationary point for \(\ell_1\)-regularized optimization problems, respectively. We provide comprehensive analyses of the iteration complexity and operation complexity, measured in terms of gradient evaluations and Hessian-vector products. We demonstrate that the proximal gradient-Newton-CG method achieves the best-known iteration complexity for attaining the proposed weak approximate second-order stationary point, consistent with results for finding an approximate second-order stationary point in unconstrained optimization. Through a toy example, we show that our proposed methods can effectively escape from a first-order approximate stationary point. Numerical experiments on the \(\ell_1\)-regularized Student's \(t\)-regression problem validate the effectiveness of both methods.

.
%
\end{abstract}

\section{Introduction}\label{intro}
In this paper, we study second-order methods for solving the optimization problem of the form 
\begin{equation}\label{eq:l1normcom}
\min_x \varphi(x) := f(x) + \lambda\|x\|_1,
\end{equation}
where \(f:~\R^n \to \R\) is twice continuously differentiable, with both its gradient and Hessian Lipschitz continuous, and may be nonconvex, and \(\lambda > 0\). Problem~\eqref{eq:l1normcom} is also known as \(\ell_1\)-regularized optimization, which has been widely applied in classification and regression models. 

Over the past two decades, numerous algorithms have been developed to address Problem~\eqref{eq:l1normcom} or a more general nonconvex composite formulation, \(\min_x~f(x) + h(x)\), where \(h(x)\) denotes a nonsmooth convex function. Examples include proximal gradient-type methods~\cite{BT09, BT10}, proximal Newton-type methods~\cite{BNO16, LW19, KL21, LPWY24, Z25}, and high-order regularized proximal Newton-type methods~\cite{DN21, DN22, LPQ24arxiv}. These methods essentially extend classical optimization frameworks---specifically, the gradient method, Newton's method, and the cubic-regularized Newton method~\cite{G81, NP06}---to the field of composite optimization. Similar to how the cubic-regularized Newton method achieves the best-known worst-case complexity results for first- and second-order optimality conditions in nonconvex unconstrained optimization when compared with its gradient-based and Newton counterparts, high-order regularized proximal Newton-type methods exhibit superior global iteration complexity for first-order optimality for nonconvex composite optimization relative to proximal gradient and proximal Newton methods~\cite{LPWY24}. However, the global iteration complexity for second-order optimality conditions remains underdeveloped for proximal-type methods in composite optimization settings.

The Newton-conjugate gradient (Newton-CG) method~\cite{RNW20} explicitly utilizes negative curvature directions derived from the Hessian of the objective function and achieves the best-known worst-case iteration complexity, with high-probability guarantees, for finding approximate second-order optimality conditions in nonconvex unconstrained optimization. Specifically, it finds an \((\varepsilon, \varepsilon^{1/2})\)-approximate second-order stationary point within \(\mathcal{O}(\varepsilon^{-3/2})\) iterations, requiring \(\widetilde{\mathcal{O}}(\varepsilon^{-7/4})\)\footnote{Logarithmic factors are omitted in \(\widetilde{\mathcal{O}}(\cdot)\)} gradient evaluations and Hessian-vector products. In contrast to the cubic-regularized Newton method, which obtains the next iterate by (approximately) solving a nonconvex cubic subproblem, the Newton-CG method utilizes the Capped-CG procedure to solve a regularized Newton equation. This procedure, which terminates in a finite number of iterations, yields either a descent direction or a negative curvature direction. This resulting direction is then employed within a line-search framework to ensure a decrease in the objective function of $\mathcal{O}(\varepsilon^{3/2})$, which is pivotal for achieving the $\mathcal{O}(\varepsilon^{-3/2})$ global iteration complexity. Recently, Xie and Wright~\cite{XW23} developed a projected Newton-CG method that establishes complexity guarantees for both first- and second-order approximate solutions in bounded-constrained optimization. This framework can be viewed as an adaptation of the Newton-CG method to composite optimization in which \(h(x)\) represents the indicator function of the box constraint set.

\textit{This paper aims to develop two second-order methods leveraging negative curvature information, with complexity guarantees for finding approximate second-order stationary points of Problem~\eqref{eq:l1normcom}.} There are three core issues to consider when applying second-order methods to \(\ell_1\)-regularized problems: (i) How can we  characterize the \((\veg, \veh)\)-approximate second-order stationary points for some \(\veg, \veh > 0\)? (ii) How can we exploit the second-order information of Problem~\eqref{eq:l1normcom}? (iii) How can we ensure sufficient descent in the objective function values during the iterations to achieve a global iteration complexity of \(\mathcal{O}(\veg^{-3/2})\) when \(\veh = \veg^{1/2}\)? We address these three issues, either fully or partially, for both methods. 
\begin{itemize}
\item[(i)] We show that the proposed hybrid proximal gradient and negative curvature method (HPGNCM) attains a strong* \((\varepsilon, \varepsilon^{1/2})\)-approximate second-order stationary point within \(\mathcal{O}(\varepsilon^{-2})\) iterations, requiring \(\mathcal{O}(\varepsilon^{-9/4})\) gradient evaluations and Hessian-vector products.  
\item[(ii)] We show that the proposed proximal gradient-Newton-CG method (PGN2CM) attains a weak \((\varepsilon, \varepsilon^{1/2})\)-approximate second-order stationary point within \(\mathcal{O}(\varepsilon^{-3/2})\) iterations, requiring \(\widetilde{\mathcal{O}}(\varepsilon^{-7/4})\) gradient evaluations and Hessian-vector products. 
\end{itemize}
To the best of our knowledge, these are the first practical methods that can directly solve the \(\ell_1\)-regularized optimization problem to find an approximate second-order stationary point while providing theoretical guarantees on both iteration and operation complexity. In particular, the complexity bounds of PGN2CM match those of the Newton-CG method for unconstrained optimization. 

\smallskip 

\noindent\textbf{Assumption and notation.} We assume that the solution set of Problem~\eqref{eq:l1normcom} is nonempty and let \(\varphi_{\rm low}\) denote the optimal objective value. We make the following key assumption on Problem~\eqref{eq:l1normcom}. 
\begin{assumption}\label{assum:problem}
For any given \(x^0\in\mathbb{R}^n\), suppose that 
\begin{itemize}
\item[(i)] the level set \(\cL_{\varphi}(x^0) = \{x\in\mathbb{R}^n\mid \varphi(x) \leq \varphi(x^0)\}\) is bounded;
\item[(ii)] \(f : \mathbb{R}^n \to (-\infty, +\infty)\) is twice uniformly Lipschitz continuously differentiable on an open neighborhood of \(\cL_{\varphi}(x^0)\).  
\end{itemize}
\end{assumption}
From Assumption~\ref{assum:problem}, there exist \(L_g \!\!>\! 0\) and \(L_H \!\!>\! 0\) such that,  for any \(x, y\!\in\!\cL_{\varphi}(x^0)\), 
\begin{subequations}
\begin{align}
&f(y) \leq f(x) + \nabla f(x)^\top (y - x) + \frac{L}{2}\|y - x\|^2, \quad \forall L \geq L_g; \label{eq:nablaf2}\\
& \|\nabla^2 f(x)(y - x) - (\nabla f(y) - \nabla f(x))\|\leq \frac{L_H}{2}\|y - x\|^2;\label{eq:p2}\\
& f(y) \leq f(x) + \nabla f(x)^\top(y - x) + \frac{1}{2}(y-x)^\top\nabla^2 f(x)(y - x) + \frac{L_H}{6}\|y - x\|^3. \label{eq:p3}
\end{align}
\end{subequations}
Throughout, \(\|\cdot\|\) denotes the Euclidean norm or its induced matrix norm. 
Moreover, there exists \(U_g > 0\), such that for any \(x\in \cL_{\varphi}(x^0)\), \(\|\nabla f(x)\| \leq U_g\), as well as \(\varphi(x) \geq \varphi_{\rm low}\) and \(\|\nabla^2 f(x) \| \leq L_g\). 

For any symmetric matrix \(A\), \(\lambda_{\min}(A)\) denotes its smallest eigenvalue, and \(A\succeq 0\) indicates that \(A\) is positive semidefinite. For the vector \(x\), \(x\geq 0\) means all components satisfying \(x_i \geq 0\). We use \(e\) and \(I\) to denote the vector with all entries equal to \(1\) and the identity matrix, respectively, without specifying their dimensions. \({\rm Diag}(s)\) denotes a diagonal matrix with the entries of the vector \(s\) on its diagonal. \({\rm sgn}(x)\) denotes the sign vector of \(x\), where \({\rm sgn}(0) = 1\). 
For \(x\in\mathbb{R}^n\), \(\max\{x, 0\} = (\max\{x_i, 0\})_{i=1}^n\). 
Given an index set \(J\subseteq [n] := \{1, \ldots, n\}\), \(x_J\in\mathbb{R}^{\vert J\vert}\) and \(A_J\in \mathbb{R}^{\vert J\vert \times \vert J\vert}\) denote the subvector of \(x\) and the principal submatrix of \(A\) restricted to the indices in \(J\), respectively, where \(\vert J\vert\) denotes the cardinality of \(J\). \(\B(x, r)\) represents the open Euclidean ball centered at \(x\) with radius \(r > 0\).

\noindent\textbf{Organization.}
The remainder of this paper is organized as follows. Section~\ref{sec:pre} introduces strong and weak approximate second-order stationary points for Problem~\eqref{eq:l1normcom} and reviews the Newton-CG method for unconstrained optimization. In Subsection~\ref{sec:hpnc}, we present the hybrid proximal gradient and negative curvature method, with its second-order complexity analysis detailed in Subsection~\ref{sec:so_hpgnc}. Subsection~\ref{sec:ncg} presents the proximal gradient-Newton-CG method, with its second-order complexity analysis detailed in Subsection~\ref{sec:scc}. Section~\ref{sec:num} includes numerical experiments conducted on a toy example as well as on the \(\ell_1\)-regularized Student's \(t\)-regression problem. Concluding remarks are provided in Section~\ref{sec:con}.

\section{Approximate second-order stationary points and related work}\label{sec:pre}
In this section, we introduce three distinct characterizations of \((\varepsilon_g, \varepsilon_h)\)-approximate second-order stationary points for Problem~\eqref{eq:l1normcom} with \(\veg, \veh > 0\), followed by a brief overview of the Newton-CG method for unconstrained optimization. We first summarize the main notation used in the stationarity conditions and subsequent analysis. Unless otherwise specified, \(x\in\mathbb{R}^n\). The index sets defined at a generic point $x$ are given in Table~\ref{tab:indexsets}. For notational convenience, when the reference point is clear from the context, we write \(I_{\star} := I_{\star}(x)\), \(I_{\star}^{\ve} := I_{\star}^{\ve}(x)\), where \(\star \in \{+, 0, -, \neq 0\}\). At iteration $k$, we further write \(I^k_{\star}:= I_{\star}(x^k)\) and \(I^{k\varepsilon}_{\star}:= I^{\varepsilon}_{\star}(x^k)\).

\begingroup  
\setlength{\tabcolsep}{9pt}  
\renewcommand{\arraystretch}{1.1}
\begin{table}[h!]
\centering
\caption{Main notation used in the stationarity conditions and subsequent analysis. }\label{tab:indexsets}
\begin{tabular*}{\textwidth}{c|c|l}\hline\hline
        Category  & Notation & \multicolumn{1}{c}{Description}      \\ \hline                                                                                             
 \multirow{4}{*}{Basic operators} & \(\partial h(x)\) & the subdifferential of a convex function $h$ at $x$ \\ 
& \({\rm prox}_h\) & \({\rm prox}_h(z) := \arg\min_y\{h(y) + \frac{1}{2}\|y - z\|^2\}\) \\ 
 & \({\rm dist}(x, \mathcal{X})\) & \({\rm dist}(x, \mathcal{X}) := \inf_{y\in\mathcal{X}}\|x - y\|\) \\ 
 & \({\rm supp}(x)\) &  \({\rm supp}(x) := \{i \mid x_i \neq 0\}\) \\ \hline
 \multirow{3}{*}{Residual mapping}  & \(\mathcal{G}_t(x)\) & \(\mathcal{G}_t(x) : = t(x - {\rm prox}_{\frac{\lambda}{t}\|x\|_1}(x - \frac{1}{t}\nabla f(x))), t > 0\) \\
 & \(g(x)\) &the exact first-order residual \\
 & \(g^{\varepsilon}(x)\) & the thresholded first-order residual \\ \hline
  \multirow{6}{*}{Index sets} &   \(I_+(x),~I_0(x)\) & \(I_+(x) := \{i \mid x_i > 0\}\);~~~~~~~~~~\(I_0(x) := \{i \mid x_i = 0\}\)\\
  &\(I_-(x),~I_{\neq 0}(x)\) & \(I_-(x) := \{i \mid x_i < 0\}\);~~~~~~~~~~\(I_{\neq 0}(x) := I_+(x) \cup I_-(x)\) \\
  &  \(I_+^{\varepsilon}(x),~I_0^{\varepsilon}(x)\)  &\(I_+^{\varepsilon}(x) := \{i \mid x_i > \veg^{1/2}\}\);~~~~~~\(I_0^{\varepsilon}(x) := \{i \mid \vert x_i\vert \leq \veg^{1/2}\}\)\\
  &\(I_-^{\varepsilon}(x),~I_{\neq 0}^{\varepsilon}(x)\) & \(I_-^{\varepsilon}(x) := \{i \mid x_i < -\veg^{1/2}\}\);~~~~\(I_{\neq 0}^{\varepsilon}(x) := I_+^{\varepsilon}(x) \cup I_-^{\varepsilon}(x)\)\\
  & \multirow{2}{*}{\((J_+^k, J_0^k, J_-^k)\)} &\(J^k_+ := \{i\in I^{k\varepsilon}_+ \mid d^k_i < 0\}\);~~~~~~\(J^k_0 := \{i\in I^{k\varepsilon}_0 \mid d^k_i < 0\}\)\\
  & & \(J^k_- := \{i\in I^{k\varepsilon}_- \mid d^k_i > 0\}\) \\ \hline\hline 
\end{tabular*}
\end{table}
\endgroup

\subsection{Approximate second-order stationary points}\label{subsec:asosp}

We first state the necessary conditions for a point \(\bar{x}\) to be a local minimizer of Problem~\eqref{eq:l1normcom}. Consider the following equivalent reformulation in nonnegative variables:
\begin{equation}\label{eq:x+x-}
\min_{x_+, x_-}F(x_+, x_-) + \lambda e^\top (x_+ + x_-),\quad {\rm s.t.}~x_+\geq 0,~x_-\geq 0, 
\end{equation}
where \(x_+ := \max\{x, 0\}\), \(x_- := \max\{-x, 0\}\), and \(F(x_+, x_-) = f(x_+ - x_-)\).
The second-order necessary condition for Problem~\eqref{eq:x+x-}, expressed in the original variable $\bar x$, requires 
\[
0 \in \nabla f(\bar{x})  + \lambda\partial \|\bar{x}\|_1, \quad
z^\top \nabla^2 f(\bar{x})z \geq 0, \quad \forall z\in \mathcal{C}_l(\bar{x}), 
\]
where 
\[
\mathcal{C}_l(\bar{x}) = \left\{
z \left\vert \begin{array}{l}
z_i \leq 0,~{\rm if}~\bar{x}_i = 0~{\rm and}~(\nabla f(\bar{x}))_i = \lambda \\
z_i = 0,~{\rm if}~\bar{x}_i = 0~{\rm and}~(\nabla f(\bar{x}))_i\in (-\lambda, \lambda)\\
z_i \geq 0,~{\rm if}~\bar{x}_i = 0~{\rm and}~(\nabla f(\bar{x}))_i = -\lambda
\end{array}\right.
\right\}. 
\]
However, verifying the latter condition is NP-hard in general~\cite{MK87}, since it amounts to checking the copositivity of \(\nabla^2 f(\bar{x})\) with respect to \(\mathcal{C}_l(\bar{x})\). A standard relaxation~\cite{ML16} replaces \(\mathcal{C}_l(\bar{x})\) with \(\mathcal{C}(\bar{x}) = \{z \mid z_i = 0~{\rm if}~\bar{x}_i = 0\}\), leading to the following weak second-order necessary optimality condition:
\begin{equation}\label{eq:localmin}
0 \in \nabla f(\bar{x})  + \lambda\partial \|\bar{x}\|_1, \quad
z^\top \nabla^2 f(\bar{x})z \geq 0, \quad \forall z\in \mathcal{C}(\bar{x}). 
\end{equation}
For any \(t > 0\), \(\mathcal{G}_t(x) = 0\) is equivalent to the first condition in~\eqref{eq:localmin}; see~\cite[Theorem 10.7]{B17}. The second condition in~\eqref{eq:localmin} can be stated as \((\nabla^2 f(\bar{x}))_{I_{\neq 0}} \succeq 0\). Notice that, when \(I_{\neq 0} = \emptyset\), we have \(\mathcal{C}(\bar{x}) = \{0\}\), and the second condition in~\eqref{eq:localmin} holds vacuously. Throughout this paper, for any symmetric matrix \(A\), we set \(\lambda_{\min}(A_{\emptyset}) := +\infty\), 
where \(A_{\emptyset}\) denotes the corresponding \(0\times 0\) principal submatrix. The same convention applies to scaled submatrices of the form \(S_{I_{\neq 0}}A_{I_{\neq 0}}S_{I_{\neq 0}}\).
Accordingly, a natural definition of approximate second-order stationary point of Problem~\eqref{eq:l1normcom} can be given as follows. 

\begin{definition}[\textbf{strong \((\varepsilon_g, \varepsilon_h)\)-2o point}]\label{def:eps-2o}
Given \(\varepsilon_g, \veh > 0\), we say that \(\bar{x}\in\mathbb{R}^n\) is a \textbf{strong \(\varepsilon_g\)-1o point} of Problem~\eqref{eq:l1normcom} if 
\[
\|\mathcal{G}_t(\bar{x})\| \leq \veg, 
\]
for some \(t > 0\). 
We say \(\bar{x}\) is a \textbf{strong \((\varepsilon_g, \varepsilon_h)\)-2o point} of Problem~\eqref{eq:l1normcom} if it is a strong \(\varepsilon_g\)-1o point and 
\begin{equation}\label{eq:so}
\lambda_{\min}((\nabla^2 f(\bar{x}))_{I_{\neq 0}}) \geq -\varepsilon_h.
\end{equation}  
\end{definition}

The strong $(\varepsilon_g, \varepsilon_h)$-2o point is motivated by the second-order necessary condition in~\ref{eq:localmin}. It is not the direct target of either proposed algorithm, but rather  serves as a natural theoretical reference. The strong* and weak notions introduced below are algorithmically motivated relaxations of this condition. In the limiting case \(\veg = \veh = 0\), a strong \((\veg, \veh)\)-2o point satisfies~\eqref{eq:localmin}. The following statement holds.   

\begin{lemma}\label{lem:epstos}
Suppose there exist sequences \(\{\veg^k\} \downarrow 0\) and \(\{\veh^k\} \downarrow 0\), and a sequence \(\{w^k\} \to w^*\) such that for each \(k\), \(w^k\) is a strong \((\veg^k, \veh^k)\)-2o point of Problem~\eqref{eq:l1normcom} with respect to a sequence \(\{t_{k+1}\}\subseteq [t_{\min}, t_{\max}]\), where \(0 < t_{\min} \leq t_{\max} < +\infty\). Then \(w^*\) satisfies~\eqref{eq:localmin}. 
\end{lemma}
\begin{proof}
Define \(I^k_+\), \(I^k_0\), \(I^k_-\), and \(I^k_{\neq 0}\) as in Table~\ref{tab:indexsets} with \(x := w^k\). Notice that 
\[
\|\mathcal{G}_{t_{k+1}}(w^k)\| \leq \veg^k \quad {\rm and}\quad (u^k)^\top (\nabla^2 f(w^k))_{I^k_{\neq 0}}u^k \geq -\veh^k\|u^k\|^2, \quad \forall u^k \in\mathbb{R}^{\vert I^k_{\neq 0}\vert}, \quad \forall k\in\mathbb{N}
\]
for some \(t_{k+1} > 0\). 
Define \(I^*_+\), \(I^*_-\), \(I^*_{0}\), and \(I^*_{\neq 0}\) as in Table~\ref{tab:indexsets} with \(x: = w^*\). 
Since \(w^k \to w^*\), \(\veg^k \downarrow 0\), and \(\veh^k\downarrow 0\), there exists \(\bar{k}\) such that for any \(k \geq \bar{k}\), we have \(I^k_{0} \subseteq I^*_0\), \(I^k_{+} \supseteq I^*_+\), and \(I^k_{-} \supseteq I^*_-\). The following statements hold. 

(i) Taking the limit in \(\|\mathcal{G}_{t_{k+1}}(w^k)\| \leq \veg^k\), and selecting a subsequence if necessary, we may assume that \(\{t_{k+1} \to t_*\}\) for some \(t_* \in [t_{\min}, t_{\max}]\). Since \(\mathcal{G}_t(x)\) is continuous in \((t, x)\) for \(t > 0\), it follows that \(\|\mathcal{G}_{t_*}(w^*)\| = 0\), which implies 
\[
{\rm dist}(0, \nabla f(x) + \lambda\partial \|x\|_1)\big\vert_{x = w^*} = 0.
\]  
Hence, the first condition in~\eqref{eq:localmin} holds with \(\bar{x} = w^*\). 

(ii) Let \(\mathcal{C}(w^*) = \{z \mid z_i = 0~{\rm if}~i\in I^*_0\}\). Then, for any \(z\in\mathcal{C}(w^*)\), we have \(z_i = 0\) for every \(i \in I_0^k\) with \(k \geq \bar{k}\), since \(I_0^k \subseteq I^*_0\). Moreover, since \(I^*_{\neq 0} \subseteq I^k_{\neq 0}\) and \(z_i = 0\) for \(i \in I^k_{\neq 0}\setminus I^*_{\neq 0}\), we have 
\[
z^\top \nabla^2 f(w^*)z = z_{I^*_{\neq 0}}^\top (\nabla^2 f(w^*))_{I^*_{\neq 0}}z_{I^*_{\neq 0}} = z_{I^k_{\neq 0}}^\top (\nabla^2 f(w^*))_{I^k_{\neq 0}}z_{I^k_{\neq 0}}, \quad \forall k\geq \bar{k}. 
\]
Recall that \(z_{I^k_{\neq 0}}^\top (\nabla^2 f(w^k))_{I^k_{\neq 0}}z_{I^k_{\neq 0}} \geq -\veh^k\|z_{I^k_{\neq 0}}\|^2\), \(\forall k\in\mathbb{N}\). By taking limits and using the continuity of \(\nabla^2 f\), the second condition in~\eqref{eq:localmin} holds with \(\bar{x} = w^*\).  
\end{proof}

\begin{remark}\label{remark:2o}
As can be seen, the optimality condition~\eqref{eq:localmin} for Problem~\eqref{eq:l1normcom} was derived by reformulating it as a nonnegative-constrained optimization problem, for which the scaled Hessian has conventionally been employed to characterize the second-order optimality condition, account for the geometry of the feasible region, and restrict curvature analysis to feasible directions~\cite{CXY10,BCY15,HLY18,OW20,XW23}.
For \(x\in\mathbb{R}^n\), let \(S = {\rm Diag}(s)\), where \(s\in\mathbb{R}^n\) is defined by \(s_i = 1\) if \(\vert x_i\vert > \veg^{1/2}\) and \(s_i = x_i\)  otherwise. Let \(I_{\neq 0}\neq \emptyset\) be defined as in Table~\ref{tab:indexsets}. Then the scaled Hessian \(S\nabla^2 f(x)S\) satisfies 
\[
\lambda_{\min}(\nabla^2 f(x)_{I_{\neq 0}}) \geq \frac{1}{\min\{1, \min_{i\in I_{\neq 0}\setminus I^\ve_{\neq 0}}\{\vert x_i\vert\}\}^2}\min\{\lambda_{\min}(S_{I_{\neq 0}}(\nabla^2 f(x))_{I_{\neq 0}}S_{I_{\neq 0}}), 0\}.
\]

The above inequality can be derived as follows. Let \(H_{\neq 0} := (\nabla^2 f(x))_{I_{\neq 0}} = U\Lambda U^\top\), where \(\Lambda = {\rm Diag}(\lambda_1, \ldots, \lambda_r, 0, \ldots, 0)\), \(\lambda_i \neq 0\), \(i = 1, \ldots, r\), is the eigen-decomposition of \(H_{\neq 0}\). Denote \(S_{\neq 0} := S_{I_{\neq 0}}\) and \(l := \vert I_{\neq 0}\vert\). Then  
\begin{align*}
\lambda_{\min}(H_{\neq 0}) = \min_{y\neq 0}\{\frac{y^\top H_{\neq 0}y }{y^\top y}\} \geq \min\{\min_{z\neq 0}\{\frac{z^\top \widehat{\Lambda} z}{z^\top z}\}, 0\},
\end{align*}
where \(\widehat{\Lambda} = {\rm Diag}(\lambda_1, \ldots, \lambda_r)\). 
Let \(U_1\) be the first \(r\) column of \(U\) and \(\widehat{S} := S_{\neq 0}U_1\). Then \(\widehat{S}\) has full columns rank and \(S_{\neq 0}H_{\neq 0}S_{\neq 0} = \widehat{S}\widehat{\Lambda}\widehat{S}^\top\). For any \(z\in\mathbb{R}^r\), there exists \(w \in\mathbb{R}^l\) such that \(z = \widehat{S}^\top w\). 

If \(\Omega:=\{w\mid w^\top\widehat{S}\widehat{\Lambda}\widehat{S}^\top w<0\}\) is empty, then we have \(S_{\neq 0}H_{\neq 0}S_{\neq 0} = \widehat{S}\widehat{\Lambda}\widehat{S}^\top \succeq 0\). Since \(S_{\neq 0}\) is nonsingular, this further implies \(H_{\neq 0} \succeq 0\). Consequently, both sides of the desired inequality are nonnegative, and the inequality holds trivially. Otherwise, 
\begin{align*}
&\lambda_{\min}(H_{\neq 0})  \geq \min\{\min_{w\in\Omega}\{\frac{w^\top \widehat{S}\widehat{\Lambda} \widehat{S}^\top w}{w^\top \widehat{S} \widehat{S}^\top w}\}, 0\} = \min\{\min_{w\in\Omega}\{\frac{w^\top S_{\neq 0}H_{\neq 0}S_{\neq 0}w}{w^\top \widehat{S} \widehat{S}^\top w}\}, 0\}\\
=&\min\{\min_{w\in\Omega}\{\frac{w^\top S_{\neq 0}H_{\neq 0}S_{\neq 0}w}{w^\top w}\frac{w^\top w}{w^\top \widehat{S} \widehat{S}^\top w}\}, 0\}\!\geq\! \min\{\min_{w\in\Omega}\{\frac{w^\top S_{\neq 0}H_{\neq 0}S_{\neq 0}w}{w^\top w}\frac{1}{\sigma^2_{\min}(\widehat{S})}\}, 0\}\\
\geq&\frac{1}{\sigma^2_{\min}(\widehat{S})}\min\{\lambda_{\min}(S_{\neq 0}H_{\neq 0}S_{\neq 0}), 0\}\geq \frac{1}{\min_{i\in I_{\neq 0}}\{\vert s_i\vert^2\}}\min\{\lambda_{\min}(S_{\neq 0}H_{\neq 0}S_{\neq 0}), 0\},
\end{align*}
where the last inequality follows from 
\[
\sigma_{\min}(\widehat{S}) = \sigma_{\min}(S_{\neq 0}U_1) \geq \min_{i\in I_{\neq 0}}\{\vert s_i\vert\} = \min\{1, \min_{i\in I_{\neq 0}\setminus I^\ve_{\neq 0}}\{\vert x_i\vert\}\}.
\] 

Therefore, inequality~\eqref{eq:so} can be verified by checking whether 
\[
\lambda_{\min}(S_{I_{\neq 0}}(\nabla^2 f(x))_{I_{\neq 0}}S_{I_{\neq 0}})\geq -(\min\{1, \min_{i\in I_{\neq 0}\setminus I^\ve_{\neq 0}}\{\vert x_i\vert\}\})^2\veh
\] 
holds.  
\end{remark}

We consider the following definition as an alternative to Definition~\ref{def:eps-2o}.
\begin{definition}[\textbf{strong* \((\varepsilon_g, \varepsilon_h)\)-2o point}]\label{def:eps-2ostar}
For \(\bar{x}\in\mathbb{R}^n\) and \(\varepsilon_g, \varepsilon_h > 0\), define \(I_{\neq 0} = \{i \mid \bar{x}_i\neq 0\}\) and \(S = {\rm Diag}(s)\), where \(s\in\mathbb{R}^n\) is defined by \(s_i = 1\) if \(\vert \bar{x}_i\vert > \veg^{1/2}\) and \(s_i = \bar{x}_i\) otherwise. 
We say that \(\bar{x}\) is a \textbf{strong* \((\varepsilon_g, \varepsilon_h)\)-2o point} of Problem~\eqref{eq:l1normcom} if it is a strong \(\varepsilon_g\)-1o point and 
\begin{equation}\label{eq:sostar}
\lambda_{\min}(S_{I_{\neq 0}}(\nabla^2 f(\bar{x}))_{I_{\neq 0}}S_{I_{\neq 0}}) \geq -\varepsilon_h.
\end{equation} 
\end{definition}

Let \(g: = g(x)\in\mathbb{R}^n\) be the vector defined as 
\begin{equation}\label{eq:gx}
g_i=\left\{
\begin{array}{ll}
(\nabla f(x))_i + \lambda, & {\rm if}~i\in I_+;\\
(\nabla f(x))_i - \lambda, & {\rm if}~i\in I_-;\\
(\nabla f(x))_i - \min\{\max\{-\lambda, (\nabla f(x))_i\}, \lambda\}, & {\rm if}~i\in I_0. 
\end{array}
\right.
\end{equation}
Then we have \({\rm dist}(0, \nabla f(x) + \lambda\partial\|x\|_1) = \|g(x)\|\). Both \(\|\mathcal{G}_t(x)\| \leq \veg\) and \(\|g(x)\| \leq \veg\) can be viewed as approximations of \(0 \in \nabla f(x) + \lambda \partial \|x\|_1\). The following lemma establishes the relationship between \(\mathcal{G}_t(x)\) and \(g(x)\). The proof is provided in Appendix~\ref{app:proofggt}. 

\begin{lemma}\label{lem:proofggt}
Suppose \(f(x)\) satisfies Assumption~\ref{assum:problem} (ii). Let \(g(x)\) be the vector defined as in~\eqref{eq:gx} for some \(x\in\mathbb{R}^n\).  Then the following statements hold. 
\begin{itemize}
\item[(a)] \(\|\mathcal{G}_t(x)\| \leq \|g(x)\|\), for all \(t > 0\).

\item[(b)] Let \(\hat{J}_+ = \{i \in I_+ \mid g_i > 0\}\) and \(\hat{J}_- = \{i\in I_- \mid g_i < 0\}\). Then 
\[ 
\mathcal{G}_t(x) = g(x), \quad \forall t \geq \hat{t} := \max\{0, \max_{i\in \hat{J}_+}\{\frac{g_i}{x_i}\}, \max_{i\in \hat{J}_-}\{\frac{g_i}{x_i}\}\}. 
\] 
Here and below, the maximum over an empty index set is understood as zero. 

\item[(c)] Let $\bar{x}$ be a stationary point of Problem \eqref{eq:l1normcom} satisfying the strict complementarity (SC) property, i.e., $0 \in \nabla f(\bar{x}) + \lambda \partial \|\bar{x}\|_1$ and $|(\nabla f(\bar{x}))_i| < \lambda$ for all $i \in \{i \mid \bar{x}_i = 0\}$. Then, for any fixed $t > 0$, there exists \(\delta > 0\) such that 
\[
    g(x) = \mathcal{G}_t(x), \quad \forall x \in \mathbb{B}(\bar{x}, \delta) \cap \mathcal{M},
\]
where $\mathcal{M} := \{x \in \mathbb{R}^n \mid \text{supp}(x) \subseteq \text{supp}(\bar{x})\}$.
\end{itemize}
\end{lemma}

To address the final question raised in the Introduction regarding the attainment of an $\mathcal{O}(\veg^{-3/2})$ global iteration complexity, it is necessary to relax the definition of the strong $\veg$-1o stationary point given in Definition~\ref{def:eps-2o}. Specifically, inspired by the Newton-CG framework for unconstrained optimization, a key requirement for achieving such a complexity bound is that the objective function decreases by at least $\mathcal{O}(\veg^{3/2})$ between successive iterates. For Problem~\eqref{eq:l1normcom}, a fundamental obstacle is that the index partition in the standard measure $g(x)$ is based on the exact condition $x_i = 0$, which prevents the line search from establishing a uniform positive lower bound on the step size. 
To overcome this, we consider a relaxed gradient vector $g^\varepsilon$ by replacing the exact support with a threshold-based partition: $\vert x_i\vert > \veg^{1/2}$ and $\vert x_i\vert \leq \veg^{1/2}$. For indices satisfying $\vert x_i\vert > \veg^{1/2}$, we set $g^{\varepsilon}_i = g_i$, keeping $g^\varepsilon$ as close to the original subgradient as possible. For coordinates near the origin, namely, those satisfying $\vert x_i\vert \leq \veg^{1/2}$, we introduce a gradient tolerance of $\veg^{3/4}$. This scaling between the numerical threshold $\veg^{1/2}$ and the tolerance $\veg^{3/4}$ is mathematically balanced to ensure that a step of sufficient length can be taken to yield the required $\mathcal{O}(\veg^{3/2})$ decrease in the objective function. This leads to the following definition.

\begin{definition}[\textbf{weak \(\varepsilon_g\)-1o point}]\label{def:eps-1o-w}
For \(x\in\mathbb{R}^n\) and \(\varepsilon_g > 0\),  let \(g^{\varepsilon} := g^{\varepsilon}(x)\in\mathbb{R}^n\) be the vector defined by 
\begin{equation}\label{eq:gxeps}
g^{\varepsilon}_i=\left\{
\begin{array}{ll}
(\nabla f(x))_i + \lambda, & {\rm if}~i\in I_+^{\varepsilon};\\
(\nabla f(x))_i - \lambda, & {\rm if}~i\in I_-^{\varepsilon};\\
(\nabla f(x))_i - \min\{\max\{-\lambda - \veg^{3/4}, (\nabla f(x))_i\}, \lambda + \veg^{3/4}\}, & {\rm if}~i\in I_0^{\varepsilon}. 
\end{array}
\right.
\end{equation}
We say that \(x\) is a \textbf{weak \(\varepsilon_g\)-1o point} of Problem~\eqref{eq:l1normcom} if \(\|g^{\varepsilon}\| \leq \veg\). 
\end{definition}
It can be verified that \(\|g^{\varepsilon}(x)\| \leq \|g(x)\|\) for any \(x\in\R^n\). As demonstrated in Section~\ref{sec:ncgethods}, this modification enables an \(\mathcal{O}(\veg^{-3/2})\) iteration complexity for attaining a weak \(\varepsilon_g\)-1o point. 

\begin{remark}
Several remarks regarding the design of $g^\varepsilon$ and the choice of exponents are in order:
\begin{itemize}
\item Consistency on Smooth Components: For $i \in I_{\neq 0}^\varepsilon$, the objective $\varphi$ is locally differentiable with respect to $x_i$. We set $g_i^\varepsilon = g_i$ to maintain consistency with the standard optimality measure. 
\item Choice of Exponent $3/4$: The exponent $3/4$ is specifically chosen to achieve the $\mathcal{O}(\varepsilon_g^{-3/2})$ complexity. As shown in Lemma~\ref{lem:dphiik0}, a proximal gradient step yields an objective decrease proportional to the square of the gradient violation, i.e., $(\varepsilon_g^{3/4})^2 = \varepsilon_g^{3/2}$. 
Replacing $3/4$ with another exponent would destroy this balance. For example, using the exponent \(1\) would yield the suboptimal $\mathcal{O}(\varepsilon_g^{-2})$ complexity typical of standard first-order methods.
\end{itemize}
\end{remark}

\begin{definition}[\textbf{weak \((\varepsilon_g, \varepsilon_h)\)-2o point}]\label{def:eps-2o-w}
For \(x\in\mathbb{R}^n\) and \(\varepsilon_g, \varepsilon_h > 0\), define \(S := {\rm Diag}(s)\), where \(s\in\mathbb{R}^n\) satisfies \(s_i = 1\) if \(\vert x_i\vert > \veg^{1/2}\) and \(s_i = x_i\) otherwise. 
We say that \(x\) is a \textbf{weak \((\varepsilon_g, \varepsilon_h)\)-2o point} of Problem~\eqref{eq:l1normcom} if it is a {weak \(\varepsilon_g\)-1o point} and 
\[
\lambda_{\min}(S_{I_{\neq 0}}(\nabla^2 f(x))_{I_{\neq 0}}S_{I_{\neq 0}}) \geq -\varepsilon_h.
\]
\end{definition}

We summarize the roles and relationships of the three first-order residuals below. The proximal-gradient mapping $\mathcal{G}_t(x)$ is used in the definition of strong first-order stationarity and in the analysis of HPGNCM. The residual \(g(x)\) is used to characterize exact first-order stationarity in the sense that \(\|g(x)\| = {\rm dist}(0, \nabla f(x) + \lambda\partial\|x\|_1)\) and to motivate the definition of the thresholded residual $g^\varepsilon(x)$. The thresholded residual $g^\varepsilon(x)$ is used in the definition of weak first-order stationarity and in the analysis of PGN2CM. Their relationships are 
\[
\|\mathcal{G}_t(x)\|\leq \|g(x)\|\quad {\rm and}\quad 
\|g^\varepsilon(x)\|\leq \|g(x)\|.
\]
Moreover, Lemma~\ref{lem:proofggt}(c) shows that
$\mathcal{G}_t(x)=g(x)$ locally under the strict complementarity condition and the support condition specified therein. No general ordering between $\|\mathcal{G}_t(x)\|$ and $\|g^\varepsilon(x)\|$ is assumed. In the limiting case \(\veg = \veh = 0\), a weak \((\veg, \veh)\)-2o point satisfies~\eqref{eq:localmin}. The following statement holds.   

\begin{lemma}\label{lem:epstoseps}
Suppose that \(\{\veg^k\} \downarrow 0\), \(\{\veh^k\} \downarrow 0 \), and a sequence \(\{w^k\}\to w^*\).
\begin{itemize}
\item[(a)] If, for every \(k\),  \(w^k\) is a strong* \((\veg^k, \veh^k)\)-2o point of Problem~\eqref{eq:l1normcom} with respect to the sequence \(\{t_{k+1}\}\subseteq [t_{\min}, t_{\max}]\), where \(0 < t_{\min} \leq t_{\max} < +\infty\), then \(w^*\) satisfies~\eqref{eq:localmin}. 
\item[(b)] If, for every \(k\), \(w^k\) is a weak \((\veg^k, \veh^k)\)-2o point of Problem~\eqref{eq:l1normcom}, then \(w^*\) satisfies~\eqref{eq:localmin}.
\end{itemize}  
\end{lemma}
\begin{proof}
(a) The first-order condition for the strong* $(\veg^k, \veh^k)$-2o point follows from Lemma~\ref{lem:epstos}. We next verify the second-order condition. Let $S^k = {\rm Diag}(s^k))$, where \(s^k_i = 1\) if \(\vert w^k_i\vert > (\veg^k)^{1/2}\) and \(s^k_i = w^k_i\) otherwise. For any $z \in \mathcal{C}(w^*)$, we have $z_i = 0$ for all $i \notin I_{\neq 0}^* := \{i \mid w^*_i \neq 0\}$. Since $w^k \to w^*$ and $\veg^k \to 0$, for every $i \in I_{\neq 0}^*$, we have $\vert w_i^k\vert > (\veg^k)^{1/2}$ for all sufficiently large $k$. Hence,  $s_i^k = 1$ for all $i \in I_{\neq 0}^*$. Moreover, for all sufficiently large \(k\), \(I^*_{\neq 0} \subseteq I^k_{\neq 0}\). Thus, by the second-order condition in the definition of a strong* $(\veg^k, \veh^k)$-2o point, 
\[
z^\top \nabla^2 f(w^k) z = z^\top (S^k)^\top \nabla^2 f(w^k) S^k z \geq -\veh^k \|z\|^2
\]
for all sufficiently large \(k\). By the continuity of \(\nabla^2 f\) and the fact that \(\veh^k \to 0\), letting \(k \to \infty\) gives $z^\top \nabla^2 f(w^*) z \geq 0$. Hence, \(w^*\) satisfies~\eqref{eq:localmin}.

(b) Since the second-order conditions in the definitions of strong* and weak $(\varepsilon_g,\varepsilon_h)$-2o points are identical, the above argument also applies to the weak $(\veg^k, \veh^k)$-2o point. It remains only to verify the first-order condition. For simplicity, let \(\varepsilon_k := \varepsilon_g^k\) and \(g^{k\varepsilon} := g^{\varepsilon_g^k}(w^k)\). Since  \(w^k\) is a weak \((\veg^k, \veh^k)\)-2o point, we have \(\|g^{k\varepsilon}\| \leq \varepsilon_k\). Consequently, 
\[
\vert g^{k\varepsilon}_i \vert \leq \varepsilon_k, \quad \forall i\in[n].
\]
Fix \(i \in [n]\). If \(w_i^* > 0\), then, since \(w_i^k \to w_i^* > 0\) and \(\varepsilon_k^{1/2} \to 0\), we have \(w_i^k > \varepsilon_k^{1/2}\) for all sufficiently large \(k\). Hence, \(i \in I^{k\epsilon}_+\), and therefore \(g^{k\varepsilon}_i = (\nabla f(w^k))_i + \lambda\). It follows that
\[
\vert (\nabla f(w^k))_i + \lambda\vert \leq \varepsilon_k.
\]
By the continuity of \(\nabla f\), we have \((\nabla f(w^*))_i = -\lambda\). 

Similarly, if \(w_i^* < 0\), then \((\nabla f(w^*))_i = \lambda\). 

It remains to consider the case \(w_i^* = 0\). We show that  
\[
\vert (\nabla f(w^k))_i\vert \leq \lambda + \varepsilon_k^{3/4} + \varepsilon_k, \quad \forall k\in\mathbb{N}. 
\]
Indeed, if \(i \in I^{k\varepsilon}_+\), then \(\vert (\nabla f(w^k))_i + \lambda \vert = \vert g^{k\varepsilon}_i\vert \leq \varepsilon_k\), which implies \(\vert  (\nabla f(w^k))_i \vert \leq \lambda + \varepsilon_k\). If \(i \in I^{k\varepsilon}_-\), then \(\vert (\nabla f(w^k))_i - \lambda \vert = \vert g^{k\varepsilon}_i\vert \leq \varepsilon_k\), which implies \(\vert  (\nabla f(w^k))_i \vert \leq \lambda + \varepsilon_k\). Finally, if \(i \in I^{k\varepsilon}_0\), then \(\vert (\nabla f(w^k))_i - \min\{\max\{-\lambda - \varepsilon_k^{3/4}, (\nabla f(w^k))_i\}, \lambda + \varepsilon_k^{3/4}\}\vert = \vert g^{k\varepsilon}_i\vert \leq \varepsilon_k\), which implies that \(\vert (\nabla f(w^k))_i\vert \leq \lambda + \varepsilon_k^{3/4} + \varepsilon_k\). 
Letting \(k \to \infty\), we obtain \(\vert (\nabla f(w^*))_i\vert \leq \lambda\). 

Consequently, we have $0 \in \nabla f(w^*) + \lambda \partial \|w^*\|_1$. Combining this first-order condition with the second-order condition established above, we conclude that \(w^*\) satisfies~\eqref{eq:localmin}.
\end{proof}

\begin{remark}
As stated at the beginning of this section, Problem~\eqref{eq:l1normcom} can be reformulated as the nonnegativity-constrained optimization problem~\eqref{eq:x+x-}, for which second-order algorithms are discussed in~\cite{XW23}. We note that when solving Problem~\eqref{eq:l1normcom} via~\eqref{eq:x+x-}, Definition~\ref{def:eps-2o} is strictly stronger than~\cite[Def. 1]{XW23}  in both the first-order and second-order conditions, while Definition~\ref{def:eps-2o-w} is weaker and has a much closer formal structure to it. We show that Definition~\ref{def:eps-2o-w} is weaker than the counterpart presented in~\cite[Def. 1]{XW23}. 

For a point \((\hat{x}_+, \hat{x}_-)\), let \(J^+ = \{i_+ \mid 0 \leq \hat{x}_+^i \leq \epsilon^{1/2}\}\cup \{i_- \mid 0 \leq \hat{x}_-^i \leq \epsilon^{1/2}\}\) and \(J^- = [n] \setminus J^+\). In~\cite{XW23}, a point \((\hat{x}_+, \hat{x}_-)\) is called an \(\epsilon\)-1o point of Problem~\eqref{eq:x+x-} if it satisfies \(\|\widehat{S}\nabla F(\hat{x}_+, \hat{x}_-)\| \leq 2\epsilon\) and \(\nabla_i F(\hat{x}_+, \hat{x}_-) \geq \epsilon^{3/4}\) for all \(i \in J^+\), where \(\widehat{S} = {\rm Diag}(\hat{s})\), with entries defined as follows: \(\hat{s}_{i_+} = \hat{x}_+^i\) if \(i_+ \in J^+\); \(\hat{s}_{i_-} = \hat{x}_-^i\) if \(i_- \in J^+\); and \(\hat{s}_i = 1\) otherwise. It follows that for any \(i_+\in J^+\) and \(i_-\in J^+\), we have \(\hat{x}_i = \hat{x}_+^i - \hat{x}_-^i \in [-\epsilon^{1/2}, \epsilon^{1/2}]\) and \(\nabla_if(\hat{x}) \in [-\lambda - \epsilon^{3/4}, \lambda + \epsilon^{3/4}]\). This implies that \(i \in I_0^{\epsilon}\) and \(g_i^{\epsilon}(\hat{x}) = 0\). When \(i_+ \in J^-\), \(i_- \notin J^-\), and \(\hat{x}_i = \hat{x}_+^i - \hat{x}_-^i < -\epsilon^{1/2}\), an additional requirement is imposed on \(\hat{x}_i(\nabla_if(\hat{x}) + \lambda)\) compared to Definition~\ref{def:eps-1o-w}; similarly, a restriction applies to \(\hat{x}_i(\nabla_if(\hat{x}) - \lambda)\) when \(i_+\notin J^-\), \(i_- \in J^-\), and \(\hat{x}_i = \hat{x}_+^i - \hat{x}_-^i >\epsilon^{1/2}\).
According to~\cite[Def. 1]{XW23}, \((\hat{x}_+, \hat{x}_-)\) is an \((\epsilon, \epsilon^{1/2})\)-2o point of Problem~\eqref{eq:x+x-}  if it additionally satisfies \(\lambda_{\min}(\widehat{S}\nabla^2 F(\hat{x}_+, \hat{x}_-)\widehat{S}) \geq -\epsilon^{1/2}\). This condition differs from Inequality~\eqref{eq:sostar} since for any  index set \(I\), \(\lambda_{\min}(\widehat{S}\nabla^2 F(\hat{x}_+, \hat{x}_-)\widehat{S}) \leq \lambda_{\min}(\widehat{S}_{I}\nabla^2 F(\hat{x}_+, \hat{x}_-)_{I}\widehat{S}_{I})\).  
Hence, if \((\hat{x}_+, \hat{x}_-)\) is an \((\epsilon, \epsilon^{1/2})\)-2o point of Problem~\eqref{eq:x+x-} as defined in~\cite[Def. 1]{XW23}, then it is also an \((\varepsilon_g, \veg^{1/2})\)-2o point by Definition~\ref{def:eps-2o-w} for some \(\veg\in(0, 1)\). 

Solving Problem~\eqref{eq:l1normcom} via the reformulation in Problem~\eqref{eq:x+x-} not only increases the problem dimensionality but also introduces a redundant representation. This redundancy results in linearly dependent gradient components and the rank-deficiency Hessian matrix.
\end{remark}

\begin{figure}[h!]
\begin{minipage}[t]{1\linewidth}
\centering
\includegraphics[width = 0.80\textwidth]{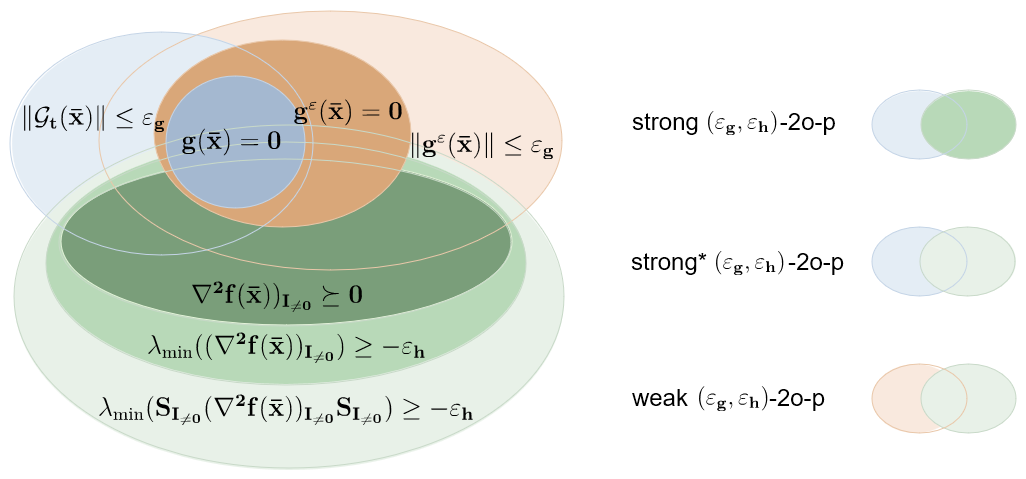}
\end{minipage}
\caption{Venn diagram illustrating the relationships among the various stationary points.}\label{fig:sp}
\end{figure}

Figure~\ref{fig:sp} illustrates the relationships among these stationary points. Table~\ref{tab:summary} summarizes the definitions of the various stationary points, along with the corresponding algorithms and theorems. 

\begingroup  
\setlength{\tabcolsep}{9pt}  
\renewcommand{\arraystretch}{1.1}
\begin{threeparttable}[h!]
\centering
\caption{Summary of the approximate stationary points and main results. ``SP" stands for stationary point, ``1o-p" and ``2o-p" refer to first-order point and second-order point, respectively. }
\small 
\label{tab:summary}
\begin{tabular*}{\textwidth}{c|c|c|cc}\hline\hline
{SP}            &       {conditions}                          & {Defs.} &Algs. & Ths.     \\ \hline \hline
strong   & \multirow{2}{*}{\(\|\mathcal{G}_t(\bar{x})\| \leq \veg\)~for some \(t > 0\)}& \multirow{2}{*}{Def.~\ref{def:eps-2o}}   & \multirow{2}{*}{Alg.~\ref{alg:hpgnc}} &   \multirow{2}{*}{Th.~\ref{th:foicpgnc}} \\             
\(\veg\)-1o-p & &&&\\ \hline                                                                                                                        
weak   &\multirow{2}{*}{\(\|g^{\varepsilon}(\bar{x})\|\leq \veg\)~~~~\eqref{eq:gxeps}}  & \multirow{2}{*}{Def.~\ref{def:eps-1o-w}}   & \multirow{2}{*}{Alg.~\ref{alg:pncg}} & \multirow{2}{*}{Th.~\ref{th:foic}}  \\ 
\(\veg\)-1o-p &   & &&\\ \hline  
strong   & ``strong \(\veg\) 1o-p" \&  & \multirow{2}{*}{Def.~\ref{def:eps-2o}}  & \multirow{2}{*}{\tnote{a}} &   \\ 
\((\veg, \veh)\)-2o-p&  ``\(\lambda_{\min}((\nabla^2 f(\bar{x}))_{I_{\neq 0}}) \geq -\varepsilon_h\)" & & 
\\ \hline        
strong*   & ``strong \(\veg\) 1o-p" \& & \multirow{2}{*}{Def.~\ref{def:eps-2ostar}}   &\multirow{2}{*}{Alg.~\ref{alg:hpgnc}} &  \multirow{2}{*}{Th.~\ref{th:soicpgnc}} \\   
\((\veg, \veh)\)-2o-p & ``\(\lambda_{\min}(S_{I_{\neq 0}}(\nabla^2 f(\bar{x}))_{I_{\neq 0}}S_{I_{\neq 0}}) \geq -\varepsilon_h\)" & &\\ \hline  
weak   & ``weak \(\veg\) 1o-p" \&  & \multirow{2}{*}{Def.~\ref{def:eps-2o-w}}  &\multirow{2}{*}{Alg.~\ref{alg:pncg}}  &  \multirow{2}{*}{Th.~\ref{th:soic}}   \\
\((\veg, \veh)\)-2o-p & ``\(\lambda_{\min}(S_{I_{\neq 0}}(\nabla^2 f(\bar{x}))_{I_{\neq 0}}S_{I_{\neq 0}}) \geq -\varepsilon_h\)" &  & 
\\ \hline  \hline
\end{tabular*}
\begin{tablenotes}
    \small
    \item[a] The strong $(\varepsilon_g, \varepsilon_h)$-2o point is not directly targeted by either algorithm.
\end{tablenotes}
\end{threeparttable}
\endgroup

\subsection{Newton-CG method}\label{subsec:newtoncgm}

The Newton-CG method~\cite{RNW20} was proposed for finding an \((\veg, \veh)\)-approximate second-order stationary point of the problem \(\min_xf(x)\) under Assumption~\ref{assum:problem} with \(\varphi = f\), where an \((\veg, \veh)\)-approximate second-order stationary point satisfies 
\[
\|\nabla f(\bar{x})\|\leq \veg \quad {\rm and}\quad \lambda_{\min}(\nabla^2 f(\bar{x}))\geq -\veh.  
\]
Given the current iterate \(x^k\) for \(k\in\mathbb{N}\), the iteration of the Newton-CG method can be divided into two cases. 

\noindent\textbf{Case I. \(\|\nabla f(x^k)\| > \veg\).} Address the regularized Newton equation
\begin{equation}\label{eq:ncg_rne}
(\nabla^2 f(x^k) + 2\veh I)d = -\nabla f(x^k)
\end{equation}
using the Capped CG method~\cite{RNW20}. Let \(d\) and \(d_{\rm type}\) be the outputs of the Capped CG method. Then \(d_{\rm type} \in \{\text{`SOL'},  \text{`NC'}\}\).  
\begin{itemize}
\item If \(d_{\rm type} =\) `SOL', then \(d\) is an approximate solution of Equation~\eqref{eq:ncg_rne}. Set \(d_k = d\). 
\item  If \(d_{\rm type} =\) `NC', then \(d\) is a negative curvature direction of \(\nabla^2 f(x^k)\). Set 
\[
d_k = -{\rm sgn}(\nabla f(x^k)^\top d)\frac{\vert d^\top \nabla^2 f(x^k)d\vert}{\|d\|^2}\frac{d}{\|d\|}.
\] 
\end{itemize}
Update \(x^{k+1}\) by \(x^{k+1} = x^k + \theta^{j_k} d_k\) for some \(\theta \in (0, 1)\), where \(j_k\) is the smallest nonnegative integer \(j\) such that 
\begin{equation}\label{eq:ncg_ls}
f(x^k + \theta^jd_k) < f(x^k) - \frac{\eta}{6}\theta^{3j}\|d_k\|^3 \quad \text{for~some}~\eta > 0. 
\end{equation}

\smallskip
\noindent\textbf{Case II. \(\|\nabla f(x^k)\| \leq \veg\).} Check whether \(\lambda_{\min}(\nabla^2 f(x^k)) > -\veh\) by performing the Minimum Eigenvalue Oracle (MEO) on \(\nabla^2 f(x^k)\). MEO returns a unit vector \(d\) that satisfies \(\nabla^2 f(x^k)d = \hat{\lambda} d\) for some \(\hat{\lambda} < -\varepsilon_h/2\), or else returns a certificate that all eigenvalues of \(\nabla^2 f(x^k)\) are greater than \(-\varepsilon_h\). In the former case, set \(d_k\) as a rescaled negative curvature direction. Update \(x^{k+1}\) by \(x^{k+1} = x^k + \theta^{j_k} d_k\), where \(j_k\) is the smallest nonnegative integer such that~\eqref{eq:ncg_ls} holds. In the latter case, the certificate may be wrong with probability \(\sigma\) for some \(\sigma \in (0, 1)\); that is, the algorithm may  terminated incorrectly if \(\hat{\lambda} > -\veh\) but \(\lambda_{\min}(\nabla^2 f(x^k)) < -\veh\). 

Notice that the regularization term \(2\veh I\) in \eqref{eq:ncg_rne} depends only on the approximate accuracy \(\veh\). 
Some other regularization terms have been discussed in the literature~\cite{HL23,ZX24,ZXBDZ25}. 
In~\cite{ZX24}, the authors introduced the regularized Newton equation
\begin{equation}\label{eq:ncg_rng}
(\nabla^2 f(x^k) + \tau_k\|\nabla f(x^k)\|^{\delta})d = -\nabla f(x^k), \quad \tau_k \in [\frac{2\veg^{1/2}}{\|\nabla f(x^k)\|^{\delta}}, \frac{2\hat{\tau}\veg^{1/2}}{\|\nabla f(x^k)\|^{\delta}}],
\end{equation}
where \(\hat{\tau} \geq 1\) and \(\delta \in [0, 1]\).  
The directions generated by solving Equation~\eqref{eq:ncg_rng} using the Capped CG method and combining it with the line search condition
\begin{equation}\label{eq:hirnc_ls}
f(x^k + \theta^jd_k) < f(x^k) - \frac{\eta}{6}\theta^{2j}\veg^{1/2}\|d_k\|^2
\end{equation}
lead to a function value decrease of \(c\veg^{3/2}\) for some \(c > 0\), matching the decrease of the Newton-CG method with \(\veh = \veg^{1/2}\). Line search condition~\eqref{eq:hirnc_ls} has been adopted in the projected Newton-CG method~\cite{XW23}, which is mainly used to reduce the order of \(L_H\) that appears in the iteration complexity result \(\mathcal{O}(\max\{\veg^{-3}\veh^3, \veh^{-3}\})\). 
Adaptive update rules for the regularization term, which depend on the estimate of the smoothness parameter \(L_H\) and can  further reduce the order of \(L_H\) in the iteration complexity result, have been presented in~\cite{HL23, ZXBDZ25}. 

\section{The hybrid proximal gradient and negative curvature method}\label{sec:hpgncmethod}

In this section, we present the hybrid proximal gradient and negative curvature method, which is designed to find a strong* $(\varepsilon_g,\varepsilon_h)$-2o point of Problem~\eqref{eq:l1normcom}.

\subsection{Algorithm}\label{sec:hpnc}

Given the current iterate \(x^k\) for \(k\in\mathbb{N}\), we first compute a proximal-gradient trial stepsize $t_{k+1}=\beta^{j_k}$, where $j_k$ is the smallest nonnegative integer \(j\) such that 
\begin{equation}\label{eq:lsprox}
\varphi({\rm prox}_{\frac{\lambda}{\beta^j}\|x\|_1}(x^k - \frac{1}{\beta^j}\nabla f(x^k))) < \varphi(x^k) - \frac{\bar{\eta}}{\beta^j}\|\mathcal{G}_{\beta^j}(x^k)\|^2 
\end{equation}
for some \(\beta > 1\) and \(\bar{\eta}\in (0, 1)\).
Then we evaluate the corresponding residual $\|\mathcal{G}_{t_{k+1}}(x^k)\|$. This computation is used for phase selection in HPGNCM. 

\smallskip
\noindent\textbf{Case I. the tested residual does not satisfy the first-order threshold (that is, \(\|\mathcal{G}_{t_{k+1}}(x^k)\| > \veg\)).} In this case, the \textbf{proximal gradient} trial \textbf{step} (ProxG step) is accepted, and we set
\[
x^{k+1} = {\rm prox}_{\frac{\lambda}{t_{k+1}}\|x\|_1}(x^k - \frac{1}{t_{k+1}}\nabla f(x^k)). 
\]

\smallskip
{\noindent\textbf{Case II. the tested residual does  satisfy the first-order threshold (that is, \(\|\mathcal{G}_{t_{k+1}}(x^k)\| \leq \veg\)).} We check whether the second-order approximate condition~\eqref{eq:sostar} holds at \(x^k\). We perform the Minimum Eigenvalue Oracle (MEO;~see Algorithm~\ref{alg:meo} in Appendix~\ref{appendix:meo}) on \(S^k_{\neq 0}(\nabla^2 f(x^k))_{\neq 0}S^k_{\neq 0}\), where \(S^k\) is defined as in Definition~\ref{def:eps-2ostar} with \(x := x^k\), \(S^k_{\neq 0} := S^k_{I^k_{\neq 0}}\), and \((\nabla^2 f(x^k))_{\neq 0} := (\nabla^2 f(x^k))_{I^k_{\neq 0}}\). MEO returns a certificate that all eigenvalues of \(S^k_{\neq 0}(\nabla^2 f(x^k))_{\neq 0}S^k_{\neq 0}\) are greater than \(-\varepsilon_h\), or else returns a unit vector \(u\) satisfying \(S^k_{\neq 0}(\nabla^2 f(x^k))_{\neq 0}S^k_{\neq 0}u = \lambda_{\min}u\) for some \(\lambda_{\min}  < -\varepsilon_h/2\). 
In the former case, the certificate may be wrong with probability \(\sigma\) for some \(\sigma \in (0, 1)\). In the latter case, we define 
\[
d^k_{I^k_{\neq 0}} = -{\rm sgn}((g^k_{I^k_{\neq 0}})^\top S^k_{\neq 0}u)\vert u^\top S^k_{\neq 0}(\nabla^2 f(x^k))_{\neq 0}S^k_{\neq 0}u\vert u
\]
and \(d^k_{I^k_0} = 0\). Then we update  \(x^k\) to \(x^{k + 1}\) by \(x^{k+1} = x^k + \theta^{j_k}S^kd^k\), where \(j_k\) is the smallest nonnegative integer \(j\) such that
\begin{equation}\label{eq:lsreq2}
\varphi(x^k + \theta^{j}S^kd^k) < \varphi(x^k) - \eta\theta^{2j}\|d^k\|^3
\end{equation}
for some \(\eta \in (0, \frac{1}{2})\). }

We summarize the HPGNCM in Algorithm~\ref{alg:hpgnc}.

\begin{algorithm}[h!]
\caption{\underline{H}ybrid \underline{P}roximal \underline{G}radient and \underline{N}egative \underline{C}urvature \underline{M}ethod (HPGNCM) for Problem~\eqref{eq:l1normcom}.}\label{alg:hpgnc}
\begin{algorithmic}[1]
\Require{ \(0 < \varepsilon_g, \varepsilon_h < 1\), \(\beta > 1\), \(\bar{\eta}\in (0, 1)\), \(\eta \in (0, \frac{1}{2})\), \(\theta \in (0, 1)\).}
\Ensure{\(\{x^k\}\)}
\For{\(k = 0, 1, 2, \ldots\)}
\State{Compute \(I^k_+\), \(I^k_0\), \(I^k_-\), and \(I^k_{\neq 0}\);}
\State{Choose \(t_{k+1} = \beta^{j_k}\) and \(j_k\) is the smallest nonnegative integer \(j\) such that~\eqref{eq:lsprox} holds and compute \(\|\mathcal{G}_{t_{k+1}}(x^k)\|\).}
\If{\(\|\mathcal{G}_{t_{k+1}}(x^k)\| > \veg\)} \Comment{\textbf{first phase}}\label{line3}
\State{Let \(x^{k+1} = {\rm prox}_{\frac{\lambda}{t_{k+1}}\|x\|_1}(x^k - \frac{1}{t_{k+1}}\nabla f(x^k))\).}\Comment{ProxG step}
\Else \Comment{\textbf{second phase}}\label{line6}
\If{\(I_{\neq 0}^k = \emptyset\)}
\State{Stop and output $x^k$.}
\Else
\State{Compute \(S^k\) as in Definition~\ref{def:eps-2ostar} with \(x = x^k\). Call Algorithm~\ref{alg:meo} (MEO) with \(H \gets S^k_{\neq 0}(\nabla^2 f(x^k))_{\neq 0}S^k_{\neq 0}\) and \(\epsilon \gets \varepsilon_h\), together with an upper bound \(M\) on \(\|H\|\) if such a bound is available. Let \((\lambda_{\min}, u)\) be the output of Algorithm~\ref{alg:meo}.}\label{line8}
\If{\(\lambda_{\min} \geq -\varepsilon_h\)}
\State{Stop and output \(x^k\).}
\Else\Comment{Negative curvature step}
\State{Let \(d^k_{I^k_{\neq 0}} \!\!\!\!=\! -{\rm sgn}((g^k_{I^k_{\neq 0}})^\top S^k_{\neq 0}u)\vert \!u^\top\! S^k_{\neq 0}(\nabla^2 f(x^k))_{\neq 0}S^k_{\neq 0}u\vert u\) and \(d^k_{I^k_{0}} \!=\! 0\).} \label{linedk}
\State{Let \(x^{k+1} = x^k + \alpha_kS^kd^k\), where \(\alpha_k = \theta^{j_k}\) and \(j_k\) is the smallest nonnegative integer \(j\) such that}\label{line11}
\[
\varphi(x^k + \theta^{j}S^kd^k) < \varphi(x^k) - \eta\theta^{2j}\|d^k\|^3.
\]
\EndIf\label{line12}
\EndIf
\EndIf
\EndFor
\end{algorithmic}
\end{algorithm}

\subsection{Global complexity results of HPGNCM}\label{sec:so_hpgnc}

In this subsection, we present the iteration and computational complexity results of Algorithm~\ref{alg:hpgnc} for finding a strong* \((\varepsilon_g, \varepsilon_h)\)-2o point by estimating lower bounds on the decrease in \(\varphi\) per iteration in each of the following two cases:
\begin{itemize}
\item[(i)] a proximal gradient step is taken (Lemma~\ref{lem:dphiik0_hpgnc});
\item[(ii)] Algorithm~\ref{alg:meo} is invoked and returns a direction of negative curvature for \(S^k_{\neq 0}(\nabla^2 f(x^k))_{\neq 0}S^k_{\neq 0}\) and a negative curvature step is taken (Lemma~\ref{lem:dmeo}).
\end{itemize}

\begin{lemma}\label{lem:dphiik0_hpgnc}
Suppose that Assumption~\ref{assum:problem} holds and that the proximal gradient step is taken. 
Then \(j_k < +\infty\) and  
\[
\varphi(x^k) - \varphi(x^{k+1}) \geq c^1_{pg}\varepsilon_g^2,
\]
where \(c^1_{pg} = \frac{\bar{\eta}}{\max\{\beta L_g/(2(1 - \bar{\eta})), 1\}}\). 
\end{lemma}
\begin{proof}
The statement follows directly from the standard theory of the classical proximal gradient method; see~\cite[Section 10.3]{B17}. 
\end{proof}

\begin{theorem}\label{th:foicpgnc}
Suppose that Assumption~\ref{assum:problem} holds for Problem~\eqref{eq:l1normcom}. Consider Algorithm~\ref{alg:hpgnc} with \(0< \veg < 1\). Define 
\begin{equation*}
K^{fo}_{hpgncm} := \left\lceil \frac{\varphi(x^0) - \varphi_{\rm low}}{c^1_{pg}}\veg^{-2}\right\rceil + 1,
\end{equation*}
where \(c^1_{pg}\) is the same as in Lemma~\ref{lem:dphiik0_hpgnc}. 
 Then Algorithm~\ref{alg:hpgnc} reaches a strong \(\varepsilon_g\)-1o point in at most \(K^{fo}_{hpgncm}\) iterations. 
\end{theorem}
\begin{proof}
Suppose, for contradiction, that \(\|\mathcal{G}_{t_{k+1}}(x^k)\| > \veg\) for \(k = 0, 1, \ldots, K^{fo}_{hpgncm}\). Then from Lemma~\ref{lem:dphiik0_hpgnc}, we have 
\[
\varphi(x^0) - \varphi(x^{K^{fo}_{hpgncm} + 1}) =\!\!\!\!\sum_{l = 0}^{K^{fo}_{hpgncm}}\!\!\varphi(x^l) - \varphi(x^{l + 1})\geq (K^{fo}_{hpgncm}+1)c_{pg}^1\veg^{2} > \varphi(x^0) - \varphi_{\rm low},
\]
where the last inequality follows from the definition of \(K^{fo}_{hpgncm}\). This contradicts the definition of \(\varphi_{\rm low}\). 
\end{proof}

Recall that \(\varphi(x)\) is lower bounded, Algorithm~\ref{alg:hpgnc} cannot get stuck in proximal gradient steps, according to Lemma~\ref{lem:dphiik0_hpgnc}.
When Algorithm~\ref{alg:meo} is invoked and identifies a direction with curvature less than or equal to \(-\frac{1}{2}\varepsilon_h\), the following property holds.

\begin{lemma}\label{lem:dphi}
Suppose Algorithm~\ref{alg:meo} is invoked and identifies a direction with curvature less than or equal to \(-\frac{1}{2}\varepsilon_h\). Let \(S^k\) and \(d^k\) be generated in steps~\ref{line8} and~\ref{linedk}, respectively, of  Algorithm~\ref{alg:hpgnc}, \(H^k_{\neq 0}:= (\nabla^2 f(x^k))_{I^k_{\neq 0}}\) and \(g^k_{\neq 0} := g^k_{I^k_{\neq 0}}\), and define  the following auxiliary bounds for the trial stepsize:
\begin{equation}\label{eq:tk+-0}
\begin{split}
& t^k_+ = \left\{
\begin{array}{ll}
\min_{i\in J^k_+}\{-\frac{x^k_i}{d^k_i}\}, & {\rm if}~J^k_+\neq \emptyset,\\
+\infty, & {\rm otherwise},
\end{array}
\right.~t^k_- = \left\{
\begin{array}{ll}
\min_{i\in J^k_-}\{-\frac{x^k_i}{d^k_i}\}, & {\rm if}~J^k_-\neq \emptyset,\\
+\infty, & {\rm otherwise}, 
\end{array}
\right. \\
&t^k_0 = \left\{
\begin{array}{ll}
\min_{i\in J^k_0}\{-\frac{1}{d^k_i}\}, & {\rm if}~J^k_0\neq \emptyset,\\
+\infty, & {\rm otherwise}.
\end{array}
\right.
\end{split}
\end{equation}
Then, for any \(t \in (0, \min\{t^k_+, t^k_-, t^k_0\}]\), we have 
\[
\varphi(x^{k} + tS^kd^k) \!\leq\! \varphi(x^k) + t(g_{\neq 0}^k)^\top S^k_{\neq 0}d^k_{\neq 0} + \frac{t^2}{2}(d^k_{\neq 0})^\top S^k_{\neq 0}H^k_{\neq 0}S^k_{\neq 0}d^k_{\neq 0}+ \frac{L_H}{6}t^3\|d^k_{\neq 0}\|^3.
\] 
\end{lemma}
\begin{proof}
From~\eqref{eq:p3}, for any \(t > 0\), we have 
\begin{align*}
\varphi(x^k + tS^kd^k) \leq & \varphi(x^k) + t\nabla f(x^k)^\top S^kd^k+ \frac{t^2}{2}(d^k)^\top S^k\nabla^2 f(x^k)S^kd^k + \frac{L_H}{6}t^3\|S^kd^k\|^3  \\
&+ \lambda(\|x^k + tS^kd^k\|_1 - \|x^k\|_1) \\
\leq &\varphi(x^k) + t(g^k)^\top S^kd^k + \frac{t^2}{2}(d^k)^\top S^k\nabla^2 f(x^k)S^kd^k + \frac{L_H}{6}t^3\|d^k\|^3  \\
&+ \lambda(\|x^k + tS^kd^k\|_1 - \|x^k\|_1) + t(\nabla f(x^k) - g^k)^\top S^kd^k, 
\end{align*}
where the second inequality follows from \(\|S^kd^k\|\leq \|d^k\|\). 

On the one hand, from the definitions of \(t^k_+\), \(t^k_-\), and \(t^k_0\), we have
\[
x_i^k + ts^k_id_i^k \geq 0,~{\rm if}~x^k_i > 0,~{\rm and}~x_i^k + ts^k_id_i^k \leq 0,~{\rm if}~x^k_i < 0, ~\forall t \in (0, \min\{t^k_+, t^k_-, t^k_0\}],
\]
which yields 
\[
\|x^k \!+\! tS^kd^k\|_1 \!-\! \|x^k\|_1 \!=\! t\!\!\sum_{i\in I^{k\varepsilon}_{\neq 0}}\!{\rm sgn}(x^k_i)d^k_i \!+\! t\!\!\!\sum_{i\in I^{k\varepsilon}_{0}}{\rm sgn}(x^k_i)x^k_id^k_i,~\forall t \in (0, \min\{t^k_+, t^k_-, t^k_0\}].
\]

On the other hand, we have 
\begin{align*}
(\nabla f(x^k) - g^k)^\top S^kd^k = -\lambda \sum_{i\in I^{k\varepsilon}_{\neq 0}}{\rm sgn}(x^k_i)d^k_i - \lambda\sum_{i\in I^{k\varepsilon}_{0}}{\rm sgn}(x^k_i)x^k_id^k_i.
\end{align*}
Therefore, for any \(t \in (0, \min\{t^k_+, t^k_-, t^k_0\}]\), we have 
\begin{align*}
\varphi(x^k \!+\! tS^kd^k) \leq  &\varphi(x^k) \!+\! t(g^k)^\top S^kd^k \!+\! \frac{t^2}{2}(d^k)^\top S^k\nabla^2 f(x^k)S^kd^k \!+\! \frac{L_H}{6}t^3\|d^k\|^3\\
=&\varphi(x^k)\!+\! t(g_{\neq 0}^k)^\top S^k_{\neq 0}d^k_{\neq 0} \!+\! \frac{t^2}{2}(d^k_{\neq 0})^\top S^k_{\neq 0}H^k_{\neq 0}S^k_{\neq 0}d^k_{\neq 0} \!+\! \frac{L_H}{6}t^3\|d^k_{\neq 0}\|^3. 
\end{align*} 
This completes the proof. 
\end{proof}

\begin{lemma}\label{lem:dmeo}
Let Assumption~\ref{assum:problem} hold and suppose that, at iteration \(k\), Algorithm~\ref{alg:meo} is invoked by Algorithm~\ref{alg:hpgnc} and identifies a direction with curvature less than or equal to \(-\frac{1}{2}\varepsilon_h\). Then, we have 
\[
\varphi(x^k) - \varphi(x^{k+1}) > \frac{1}{8}c_{nc}\min\{\veg\veh, \veh^3\},
\]
where \(c_{nc} = \eta\theta^2\min\{1, \frac{9(1 - 2\eta)^2}{L_H^2}\}\). 
\end{lemma}
\begin{proof}
Let \((\lambda_{\min}, u)\) be the outputs of Algorithm~\ref{alg:meo}. Then, we have 
\[
u^\top S^k_{\neq 0}H^k_{\neq 0}S^k_{\neq 0}u = \lambda_{\min} \leq -\frac{1}{2}\varepsilon_h \quad {\rm and}\quad \|u\|\!=\! 1.
\]
Moreover,  
\begin{equation}
\begin{split}\label{eq:meo}
&(g^k_{\neq 0})^\top S^k_{\neq 0}d^k_{\neq 0} = -\vert (g^k_{\neq 0})^\top S^k_{\neq 0}u\vert\vert u^\top S^k_{\neq 0}H^k_{\neq 0}S^k_{\neq 0}u\vert \leq 0; \\
&(d_{\neq 0}^k)^\top S^k_{\neq 0}H^k_{\neq 0}S^k_{\neq 0}d_{\neq 0}^k = \vert u^\top S^k_{\neq 0}H^k_{\neq 0}S^k_{\neq 0}u\vert^3 = \lambda_{\min}^3 = -\|d^k\|^3.
\end{split}
\end{equation}
Moreover, \(\|d^k\| = -\lambda_{\min} \geq \frac{1}{2}\varepsilon_h\).  
From Lemma~\ref{lem:dphi}, for any \(0 \!<\! t \!<\! \min\{t^k_+, t^k_-, t^k_0\}\), we have 
\begin{align*}
\varphi(x^k \!+\! tS^kd^k) \leq& \varphi(x^k) \!+\! t(g^k_{\neq 0})^\top \!S^k_{\neq 0}d^k_{\neq 0} \!+\! \frac{1}{2}t^2(d^k_{\neq 0})^\top \!S^k_{\neq 0}H^k_{\neq 0}S^k_{\neq 0}d^k_{\neq 0} \!+\! \frac{L_H}{6}t^3\|d^k_{\neq 0}\|^3\\
\overset{\eqref{eq:meo}}{\leq}& \varphi(x^k) + (\frac{L_H}{6}t - \frac{1}{2})t^2\|d_{\neq 0}^k\|^3,
\end{align*} 
which yields 
\[
\varphi(x^k + tS^kd^k)  \leq \varphi(x^k) + (\frac{L_H}{6}t - \frac{1}{2})t^2\|d^k\|^3 < \varphi(x^k)  - \eta t^2\|d^k\|^3
\]
for any \(t\in(0, \min\{t^k_+, t^k_-, t^k_0, \frac{3(1 - 2\eta)}{L_H}\})\). Noting that \(\vert d^k_i\vert \leq \|d^k\| = -\lambda_{\min}\), we have \(\min\{t^k_+, t^k_-, t^k_0\} \geq \frac{\veg^{1/2}}{\vert\lambda_{\min}\vert}\). Hence, \(j_k < +\infty\) and \(\theta^{j_k} \geq \min\{\theta\min\{t^k_+, t^k_-, t^0_k,  \frac{3(1 - 2\eta)}{L_H}\}, 1\}\). Define $\alpha_k := \theta^{j_k}$. Therefore,  
\[
\alpha_k\|d^k\| = \theta^{j_k}\|d^k\| \geq \min\{\theta\min\{t^k_+\|d^k\|, t^k_-\|d^k\|, t^0_k\|d^k\|,  \frac{3(1 - 2\eta)\|d^k\|}{L_H}\}, \|d^k\|\}.
\]
From the definition of \(t^k_+\), \(t^k_-\), and \(t^k_0\), together with \(\vert d_i^k\vert \leq \|d^k\|\), we have \(\min\{t^k_+\|d^k\|, t^k_-\|d^k\|\} \geq \veg^{\frac{1}{2}}\) and \(t^k_0\|d^k\| \geq 1\). Recalling that \(\|d^k\| = -\lambda_{\min} \geq \frac{1}{2}\veh\), we have 
\[
\alpha_k\|d^k\| \geq \min\{\theta\min\{\veg^{\frac{1}{2}}, 1,  \frac{3(1 - 2\eta)\veh}{2L_H}\}, \frac{1}{2}\veh\} \geq \theta\min\{\varepsilon_g^{\frac{1}{2}}, \frac{1}{2}\min\{1, \frac{3(1 - 2\eta)}{L_H}\}\varepsilon_h\}.
\]
Therefore, 
\begin{align*}
\varphi(x^k) - \varphi(x^k + \alpha_k S^kd^k) >& \eta \alpha_k^2\|d^k\|^3 \geq \frac{\eta\theta^2}{2}\min\{\varepsilon_g, \frac{1}{4}\min\{1, \frac{9(1 - 2\eta)^2}{L_H^2}\}\varepsilon_h^2\}\varepsilon_h\\
\geq&\frac{\eta\theta^2}{8}\min\{1, \frac{9(1 - 2\eta)^2}{L_H^2}\}\min\{\veg\veh, \veh^3\}.
\end{align*}
The desired result holds. 
\end{proof}

\begin{theorem}\label{th:soicpgnc}
Suppose that Assumption~\ref{assum:problem} holds for Problem~\eqref{eq:l1normcom}. Consider Algorithm~\ref{alg:hpgnc} with \(0< \veg, \veh < 1\). Define 
\begin{equation*}
K^{so}_{hpgncm} := \left\lceil \frac{\varphi(x^0) - \varphi_{\rm low}}{\min\{c_{nc}/8, c^1_{pg}\}}\max\{\veg^{-1}\veh^{-1}, \veh^{-3},  \veg^{-2}\}\right\rceil + 1,
\end{equation*}
where \(c^1_{pg}\) and \(c_{nc}\) are the same as those in Lemmas~\ref{lem:dphiik0_hpgnc} and~\ref{lem:dmeo}, respectively. 
 Then, with probability at least \((1 - \sigma)^{K^{so}_{hpgncm}}\), Algorithm~\ref{alg:hpgnc} terminates at a strong* \((\varepsilon_g, \varepsilon_h)\)-2o point within at most \(K^{so}_{hpgncm}\) iterations, where \(\sigma \in [0, 1)\) is the probability of failure of Algorithm~\ref{alg:meo}. In particular, if \(\veg = \varepsilon\) and \(\veh = \varepsilon^{\nu}\) with \(\nu \in (0, \frac{2}{3}]\), then Algorithm~\ref{alg:hpgnc} outputs a strong* \((\varepsilon, \varepsilon^{\nu})\)-2o point with probability at least \((1 - \sigma)^{\overline{K}^{so}_{hpgncm}}\) within 
 \[
\overline{K}^{so}_{hpgncm} := \left\lceil \frac{\varphi(x^0) - \varphi_{\rm low}}{\min\{c_{nc}/8, c^1_{pg}\}}\varepsilon^{-2}\right\rceil + 1
\]
iterations. 
\end{theorem}
\begin{proof}
Suppose, for contradiction that Algorithm~\ref{alg:hpgnc} runs for \(K^{so}_{hpgncm}\) iterations without terminating. We consider the following partition of iteration indices:
\[
\mathcal{K}_1 \cup \mathcal{K}_2 = \{0, 1, \cdots, K^{so}_{hpgncm}- 1\},
\]
where \(\mathcal{K}_1\) and \(\mathcal{K}_2\) are defined as follows. 

\smallskip 

\noindent\textbf{Case 1.} \(\mathcal{K}_1 := \{l = 0, 1, \cdots, K^{so}_{hpgncm} - 1: {\rm the}~\textbf{second~phase}~{\rm is~ invoked~at~the}~l\text{-}th~{\rm step} \}\). At each iteration \(l \in \mathcal{K}_1\), Algorithm~\ref{alg:hpgnc} calls Algorithm~\ref{alg:meo} and identifies a direction of negative curvature for \(S^l_{\neq 0}(\nabla^2 f(x^l))_{\neq 0}S^l_{\neq 0}\) since \(l < K^{so}_{hpgncm}\) and the algorithm continues to iterate. By Lemma~\ref{lem:dmeo}, we have 
\begin{equation}\label{eq:dvarphi1pgnc}
\varphi(x^l) - \varphi(x^{l+1}) > \frac{1}{8}c_{nc}\min\{\veg\veh, \veh^3\}, \quad \forall l\in\mathcal{K}_1.
\end{equation}

\smallskip 

\noindent\textbf{Case 2.} \(\mathcal{K}_2 := \{l \!=\! 0, 1, \cdots, K^{so}_{hpgncm} \!-\! 1: {\rm the}~\textbf{first~phase}~{\rm is~invoked}\}\). According to Lemma~\ref{lem:dphiik0_hpgnc}, we have 
\begin{align}\label{eq:dvarphi2pgnc}
\varphi(x^l) - \varphi(x^{l+1}) >&c^1_{pg}\veg^{2}, \quad \forall l\in\mathcal{K}_2. 
\end{align}

Using the monotonicity of \(\{\varphi(x^l)\}_{l\in\mathbb{N}}\), \eqref{eq:dvarphi1pgnc}, and \eqref{eq:dvarphi2pgnc}, we have
\begin{align*}
&\varphi(x^0) - \varphi_{\rm low} \geq \sum_{l \in \mathcal{K}_1}(\varphi(x^l) - \varphi(x^{l+1})) \geq \frac{1}{8}c_{nc}\min\{\veg\veh, \veh^3\}\!\vert\mathcal{K}_1\vert; ~{\rm and}\\ 
&\varphi(x^0) - \varphi_{\rm low} \geq \sum_{l \in \mathcal{K}_2}(\varphi(x^l) - \varphi(x^{l+1})) \geq c^1_{pg}\varepsilon_g^{2}\vert\mathcal{K}_2\vert,
\end{align*}
which yields
\[
\vert\mathcal{K}_1\vert  \leq \frac{\varphi(x^0) - \varphi_{\rm low}}{c_{nc}/8}\max\{\veg^{-1}\veh^{-1}, \veh^{-3}\}\quad {\rm and}\quad\vert\mathcal{K}_2\vert \leq \frac{\varphi(x^0) - \varphi_{\rm low}}{c^1_{pg}}\veg^{-2}.
\]
Hence, we have 
\begin{align*}
K^{so}_{hpgncm} =& \vert \mathcal{K}_1\vert + \vert \mathcal{K}_2\vert \leq \frac{\varphi(x^0) - \varphi_{\rm low}}{\min\{c_{nc}/8, c^1_{pg}\}}\max\{\veg^{-1}\veh^{-1}, \veh^{-3},  \veg^{-2}\} \leq K^{so}_{hpgncm} - 1,
\end{align*}
where the last inequality follows from the definition of \(K^{so}_{hpgncm}\). This leads to a contradiction. 
 
 The probability that the output \(x^k\) fails to satisfy \(\lambda_{\min}(S^k_{\neq 0}(\nabla^2 f(x^k))_{\neq 0}S^k_{\neq 0}) \geq -\veh\) within \(K^{so}_{hpgncm}\) iteration is bounded above by \(1 - (1 - \sigma)^{K^{so}_{hpgncm}}\), follows the same augment used in~\cite[Theorem 3]{XW23} and is therefore omitted.

If \(\veg = \varepsilon\) and \(\veh = \varepsilon^{\nu}\) with \(\nu\in(0, \frac{2}{3}]\), then we have \(\max\{\veg^{-1}\veh^{-1}, \veh^{-3},  \veg^{-2}\} = \varepsilon^{-2}\).  
\end{proof}

We now give the operation complexity of Algorithm~\ref{alg:hpgnc} in terms of gradient evaluations (\(\nabla f(x)\), or equivalently, \(\mathcal{G}_t(x)\)) and Hessian-vector products (\(\nabla^2 f(x)\)). The proof depends on the complexity analysis of MEO in~\cite{RNW20}. We make the following assumption on MEO.

\begin{assumption}\label{assume:meo}
At every iteration \(k\) at which MEO (Algorithm~\ref{alg:meo}) is invoked, for a specified failure probability \(\sigma\) with \(0\leq \sigma \ll 1\), MEO either certifies that \(H \succeq -\epsilon I\) or finds a vector of curvature smaller than \(-\frac{1}{2}\epsilon\) in at most 
\[
N^{\rm meo} := \min\{n, 1 + \lceil \mathcal{C}_{meo}\epsilon^{-\frac{1}{2}}\rceil\}
\]
Hessian-vector products, with probability \(1 - \sigma\), where \(\mathcal{C}_{meo}\) depends at most logarithmically on \(\sigma\) and \(\epsilon\). 
\end{assumption}

\begin{corollary}\label{cor:hpgnc}
Suppose that the assumptions of Theorem~\ref{th:soicpgnc} and Assumption~\ref{assume:meo} hold. Let \(K^{so}_{hpgncm}\) be defined as in Theorem~\ref{th:soicpgnc}. 
Then, with probability at least \((1 - \sigma)^{K^{so}_{hpgncm}}\), Algorithm~\ref{alg:hpgnc} terminates at a strong* \((\veg, \veh)\)-2o point after at most
\[
\mathcal{O}(K^{so}_{hpgncm}\min\{n, \veh^{-1/2}\ln(n\sigma^{-1})\})
\]
Hessian-vector products and/or gradient evaluations. The probability that the algorithm terminates incorrectly, outputting a strong \(\veg\)-1o point that is not a strong* \((\veg, \veh)\)-2o point, is at most \(1 - (1 - \sigma)^{K^{so}_{hpgncm}}\).

For \(n\) sufficiently large, the following statements hold:
\begin{itemize}
\item[(i)] With probability at least \((1 - \sigma)^{K^{so}_{hpgncm}}\), the bound is 
\[
\mathcal{O}(\max\{\veg^{-1}\veh^{-3/2}, \veh^{-7/2}, \veg^{-2}\veh^{-1/2}\}).
\] 
\item[(ii)] If \(\veg = \varepsilon\) and \(\veh = \varepsilon^{\nu}\) with \(\nu \in (0, \frac{2}{3}]\), then, with probability at least \((1 - \sigma)^{\overline{K}^{so}_{hpgncm}}\), the bound is \({\mathcal{O}}(\varepsilon^{-(2 + \frac{1}{2}\nu)})\). 
\end{itemize}
\end{corollary}
\begin{proof}
According to \cite[Lemma 2]{RNW20}, when the Lanczos method is used to estimate the smallest eigenvalue of \(H\) at iterate \(k\), starting with a random vector uniformly generated on the unit sphere, Assumption~\ref{assume:meo} holds with
\[
\mathcal{C}_k^{\rm meo} = \min\{n, 1 + \lceil \frac{1}{2}\ln(2.75n/\sigma^2)\sqrt{L_g}\rceil\} = \min\{n, \mathcal{O}(\ln(n/\sigma))\}. 
\]  
Let \(\epsilon = \veh\). Then the bound on the number of (scaled) Hessian-vector products required by MEO at iteration \(k\) is 
\[
N_k^{\rm meo} := \min\{n, \mathcal{O}(\veh^{-1/2}\ln(n/\sigma))\}.
\] 

Therefore, using Theorem~\ref{th:soicpgnc}, the total number of Hessian-vector products of Algorithm~\ref{alg:hpgnc} can be bounded by
\begin{align*}
\sum_{k=0}^{K^{so}_{hpgncm} - 1}\!\!\!\!N_k^{\rm meo} =&\sum_{k=0}^{K^{so}_{hpgncm} - 1}\!\!\!\!\min\{n, \mathcal{O}(\veh^{-1/2}\ln(n/\sigma))\} = \mathcal{O}(K^{so}_{hpgncm}\min\{n, \veh^{-1/2}\ln(n/\sigma)\})\\
=&\mathcal{O}(\max\{\veg^{-1}\veh^{-1}, \veh^{-3}, \veg^{-2}\}\min\{n, \veh^{-1/2}\ln(n\sigma^{-1})\}).
\end{align*}
(i) For \(n\) sufficiently large, the above bound becomes 
\[
\mathcal{O}(\max\{\veg^{-1}\veh^{-3/2}, \veh^{-7/2}, \veg^{-2}\veh^{-1/2}\}\ln(n\sigma^{-1})) = \mathcal{O}(\max\{\veg^{-1}\veh^{-3/2}, \veh^{-7/2}, \veg^{-2}\veh^{-1/2}\}).
\]
(ii) When \(\veg = \varepsilon\) and \(\veh = \varepsilon^{\nu}\) with \(\nu \in (0, \frac{2}{3}]\), the bound in statement (i) becomes  
\[
\mathcal{O}(\max\{\varepsilon^{-(1 + 3/2\nu)}, \varepsilon^{-7/2\nu}, \varepsilon^{-(2 + \frac{1}{2}\nu)}\}\ln(n\sigma^{-1})) = \mathcal{O}(\varepsilon^{-(2 + \frac{1}{2}\nu)}). 
\]
The desired statements hold.
\end{proof}

We note that the term $\mathcal{O}(\ln(n/\sigma))$ in $N_k^{\rm meo}$ does not lead to any computational blow-up in practice. The failure probability $\sigma$ for MEO is treated as a fixed positive parameter. Because the operation complexity depends on $\sigma$ only logarithmically, choosing a highly conservative value (e.g., $\sigma = 10^{-6}$) increases the bound only by a small constant factor. Therefore, the logarithmic dependence on $1/\sigma$ does not lead to any practical divergence in computation time.

\section{The proximal gradient-Newton-CG method}\label{sec:ncgethods}

In this section, we present the proximal gradient-Newton-CG method, which is designed to find a weak \((\veg, \veh)\)-2o point of Problem~\eqref{eq:l1normcom}.

\subsection{Algorithm}\label{sec:ncg}

Let \(x^k\) be the current iterate for \(k\in\mathbb{N}\). Let \(g^{k\varepsilon} := g^{\varepsilon}(x^k)\) be the vector defined in~\eqref{eq:gxeps}. Then, the iteration of PGN2CM is divided into two cases.

\smallskip
\noindent\textbf{Case I. \(x^k\) is not a weak \(\veg\)-1o point.}
\begin{itemize}
\item[i).] Suppose \(I^{k\varepsilon}_0\neq \emptyset\) and \(g^{k\varepsilon}_i \neq 0\) for some \(i\in I^{k\varepsilon}_0\). We perform the \textbf{proximal gradient step}, that is, 
\[
x^{k+1} = {\rm prox}_{\frac{\lambda}{\alpha_k}\|x\|_1}(x^k - \frac{1}{\alpha_k}\nabla f(x^k)), \qquad \alpha_k = \beta^{j_k},
\]
where \(\beta > 1\) and \(j_k\) is the smallest nonnegative integer \(j\) such that~\eqref{eq:lsprox} holds. 

\item[ii).] Suppose \(I^{k\varepsilon}_0 = \emptyset\) or \(g^{k\varepsilon}_i = 0\) for all \(i \in I^{k\varepsilon}_0\), \(I^{k\varepsilon}_{\neq 0} \neq \emptyset\) and \(\|g^{k\varepsilon}_{I^{\varepsilon}_{\neq 0}}\| > \veg\).  
In this case, we set the line-search direction \(d^k\) to satisfy \(d^k_{I^{k\varepsilon}_0} = 0\) and determine \(d^k_{\neq 0\varepsilon} (:=d^k_{I^{k\varepsilon}_{\neq 0}})\) as follows. 
We perform the \textbf{Newton-CG step} on the following regularized Newton equation:
\begin{equation}\label{eq:rnecopycopy}
(H^k_{\neq 0\varepsilon} + \tau_k\|g^{k}_{\neq 0\varepsilon}\|^{\delta}I)d = -g^{k}_{\neq 0\varepsilon},
\end{equation}
where \(H^k_{\neq 0\varepsilon} := (\nabla^2 f(x^k))_{I^{k\varepsilon}_{\neq 0}}\), \(g^{k}_{\neq 0\varepsilon} := g^{k\varepsilon}_{I^{k\varepsilon}_{\neq 0}}\), and \(\tau_k \in [ \frac{2\varepsilon_h}{\|g^{k}_{\neq 0\varepsilon}\|^{\delta}}, \frac{2\hat{\tau}\varepsilon_h}{\|g^{k}_{\neq 0\varepsilon}\|^{\delta}}]\) for some \(\hat{\tau} \geq 1\) and \(\delta \in [0, 1]\).
We use the Capped CG method (see Algorithm~\ref{alg:ccg} in Appendix~\ref{appendix:ccg}) on equation~\eqref{eq:rnecopycopy} as in the Newton-CG method (see also~\cite{ZX24}). Let \(d\) and \(d_{\rm type}\) be the outputs of the Capped CG method. 
Then the following properties hold.
\begin{itemize}
\item[a).] \(d_{\rm type} = \) `SOL', which implies that \(d\) is an approximate solution to equation~\eqref{eq:rnecopycopy}. According to~\cite[Lemma~1]{ZX24} (similarly to~\cite[Lemma~3]{RNW20}), we have  
\begin{subequations}
\begin{align}
d^\top(H^k_{\neq 0\varepsilon} + \tau_k\|g^k_{\neq 0\varepsilon}\|^{\delta}I)d \geq \veh\|d\|^2&; \label{eq:dsolacopy}\\
\|d\| \leq  1.1\veh^{-1}\|g^k_{\neq 0\varepsilon}\|&; \label{eq:dsolbcopy}\\
\|\hat{r}_k\| \leq \frac{1}{2}\veh\zeta\|d\|&, \label{eq:dsolccopy}
\end{align}
\end{subequations}
where \(\hat{r}_k = (H^k_{\neq 0\varepsilon} + \tau_k\|g^k_{\neq 0\varepsilon}\|^{\delta}I)d + g^k_{\neq 0\varepsilon}\) and \(\zeta > 0\) is an input parameter of the Capped CG method.
In this case, we set \(d^k_{\neq 0\varepsilon} \triangleq d^k_{I^{k\varepsilon}_{\neq 0}}= d\). We update \(x^{k+1}\) as \(x^{k+1} = x^k + \theta^{j_k} d^k\) for some \(\theta \in (0, 1)\), where \(j_k\) is the smallest nonnegative integer \(j\) such that 
\begin{equation}\label{eq:lsreq}
\varphi(x^k + \theta^{j}d^k) < \varphi(x^k) - \eta\theta^{2j}\veh\|d^k\|^2
\end{equation}
for some \(\eta \in (0, 1)\). 

\item[b).] \(d_{\rm type} = \) `NC', which implies that \(d\) is a direction of negative curvature for \(H^k_{\neq 0\varepsilon}\).  According to~\cite[Lemma~1]{ZX24} (similarly to~\cite[Lemma~3]{RNW20}), we have  
\[
d^\top(H^k_{\neq 0\varepsilon} + \tau_k\|g^k_{\neq 0\varepsilon}\|^{\delta}I)d \leq \veh\|d\|^2.
\]
In this case, we set 
\[
d^k_{\neq 0\varepsilon} = -{\rm sgn}(d^\top g^k_{\neq 0\varepsilon})\frac{\vert d^\top H^k_{\neq 0\varepsilon} d\vert}{\|d\|^2}\frac{d}{\|d\|}.
\]
From~\cite[Lemma~1]{ZX24}, \(d^k_{\neq 0\varepsilon}\) satisfies
\begin{equation}\label{eq:pdnc}
(d^k_{\neq 0\varepsilon})^\top g^k_{\neq 0\varepsilon} \leq 0 \quad {\rm and}\quad \frac{(d^k_{\neq 0\varepsilon})^\top H^k_{\neq 0\varepsilon}d^k_{\neq 0\varepsilon}}{\|d^k_{\neq 0\varepsilon}\|^2} = -\|d^k_{\neq 0\varepsilon}\| \leq -\veh. 
\end{equation}
We update \(x^{k+1}\) as \(x^{k+1} = x^k + \theta^{j_k} d^k\), where \(j_k\) is the smallest nonnegative integer \(j\) such that~\eqref{eq:lsreq} holds. 
\end{itemize}
\end{itemize}

When neither (i) nor (ii) holds, we have 
\[
\|g^{\varepsilon}(x^k)\|^2 = \|g^{k\varepsilon}\|^2 = \|g^{k\ve}_{I^{k\varepsilon}_0}\|^2 + \|g^{k\ve}_{I^{k\varepsilon}_{\neq 0}}\|^2 \leq \veg^2,
\]
which implies that \(x^k\) is a weak \(\varepsilon_g\)-1o point. 

\smallskip
\noindent\textbf{Case II. \(x^k\) is a weak \(\veg\)-1o point.} We check whether the second-order approximate condition in Definition~\ref{def:eps-2o-w} holds for \(x^k\) similarly to \textbf{Case II} of HPGNCM.  

We summarize PGN2CM in Algorithm~\ref{alg:pncg}.

\begin{algorithm}[h!]
\caption{\underline{P}roximal \underline{G}radient-\underline{N}ewton-{CG} with \underline{N}egative \underline{C}urvature \underline{M}ethod (PGN2CM) for Problem~\eqref{eq:l1normcom}.}\label{alg:pncg}
\begin{algorithmic}[1]
\setlength{\itemsep}{.01pt}
\Require{\(0 < \veg < 1\), \(\veg^{1/2} \leq \veh < 1\), \(\beta >1\), \(\delta \in [0, 1]\), \(\hat{\tau} \geq 1\), \(\zeta \in (0, 1)\), \(\bar{\eta} \in (0, 1)\), \(\eta \in (0, \frac{1 - \zeta}{2})\), \(\theta \in (0, 1)\).}
\Ensure{\(\{x^k\}\)}
\For{\(k = 0, 1, 2, \ldots\)}
\State{Compute \(g^{k\varepsilon}  \gets g^{\varepsilon}(x^k)\), \(I^k_0\), \(I^k_{\neq 0}\), \(I^{k\varepsilon}_0\), and \(I^{k\varepsilon}_{\neq 0}\);}
\If{\(I^{k\varepsilon}_0 \neq \emptyset\) and \(\|g^{k\ve}_{I^{k\varepsilon}_0}\| \neq 0\)}\Comment{\textbf{first phase}}\Comment{{ProxG step}}
\State{Let  \(x^{k+1} = {\rm prox}_{\frac{\lambda}{\beta^{j_k}}\|x\|_1}(x^k - \frac{1}{\beta^{j_k}}\nabla f(x^k))\), where \(j_k\) is the smallest nonnegative integer \(j\) such that~\eqref{eq:lsprox} holds.}
\ElsIf{(\(I^{k\varepsilon}_0 = \emptyset\) or \(\|g^{k\ve}_{I^{k\varepsilon}_0}\| = 0\)) and \(\|g^{k\ve}_{I^{k\varepsilon}_{\neq 0}}\| > \varepsilon_g\)}\label{line5}
\State{Compute \(H^k_{\neq 0\varepsilon} = (\nabla^2 f(x^k))_{I^{k\varepsilon}_{\neq 0}}\) and \(g^k_{\neq 0\varepsilon} = g^{k\varepsilon}_{I^{k\varepsilon}_{\neq 0}}\). Let \(d^k_{I_0^{k\varepsilon}} = 0\); Call Algorithm~\ref{alg:ccg} with \(H \gets H^k_{\neq 0\varepsilon}\), \(g \gets g^k_{\neq 0\varepsilon} \), \(\veh\), \(\zeta\), \(\delta\), and \(\tau_k \in [ \frac{2\veh}{\|g^k_{\neq 0\varepsilon}\|^{\delta}}, \frac{2\hat{\tau}\veh}{\|g^k_{\neq 0\varepsilon}\|^{\delta}}]\) to obtain outputs \(d\), \(d_{\rm type}\);}\Comment{Newton-CG step}
\If{\(d_{\rm type} = NC\)}
\State{Let \(d^k_{\neq 0\varepsilon} := -{\rm sgn}(d^\top g^k_{\neq 0\varepsilon})\frac{\vert d^\top H^k_{\neq 0\varepsilon} d\vert}{\|d\|^2}\frac{d}{\|d\|}\);}
\Comment{Negative curvature step}
\Else
\State{Let \(d^k_{\neq 0\varepsilon} := d\);}
\Comment{Approx. solution}
\EndIf 
\State{Let \(x^{k+1} \!=\! x^k \!+\! \alpha_kd^k\), where \(\alpha_k \!=\! \theta^{j_k}\) and \(j_k\) is the smallest nonnegative integer \(j\) such that~\eqref{eq:lsreq} holds.}
\Else\Comment{\textbf{second phase}}\label{line15}
\If{\(I_{\neq 0}^k = \emptyset\)}
\State{Stop and output $x^k$.}
\Else
\State{Compute \(S^k\) as in Definition~\ref{def:eps-2ostar} with \(x = x^k\). Call Algorithm~\ref{alg:meo} (MEO) with \(H \gets S^k_{\neq 0}(\nabla^2 f(x^k))_{\neq 0}S^k_{\neq 0}\) and \(\epsilon \gets \varepsilon_h\), together with an upper bound \(M\) on \(\|H\|\) if such a bound is available. Let \((\lambda_{\min}, u)\) be the output of Algorithm~\ref{alg:meo}.}
\If{\(\lambda_{\min} \geq -\varepsilon_h\)}
\State{Stop and output \(x^k\).}
\Else\Comment{Negative curvature step}
\State{Let \(d^k_{I^k_{\neq 0}} \!\!\!=\! -{\rm sgn}((g^k_{I^k_{\neq 0}})^\top\! S^k_{\neq 0}u)\vert u^\top\! S^k_{\neq 0}(\nabla^2 f(x^k))_{\neq 0}S^k_{\neq 0}u\vert u\) and \(d^k_{I^k_{0}} \!\!=\! 0\).}
\State{Let \(x^{k+1} = x^k + \alpha_kS^kd^k\), where \(\alpha_k = \theta^{j_k}\) and \(j_k\) is the smallest nonnegative integer \(j\) such that~\eqref{eq:lsreq2} holds.}
\EndIf\label{line21}
\EndIf
\EndIf
\EndFor
\end{algorithmic}
\end{algorithm}

\subsection{Global complexity results of PGN2CM}\label{sec:scc}

In this subsection, we present the iteration and operation complexity results of Algorithm~\ref{alg:pncg} for finding a weak \((\varepsilon_g, \varepsilon_h)\)-2o point by estimating lower bounds on the decrease in \(\varphi\) per iteration for each of the following four cases: 
\begin{itemize}
\item[(i)] A proximal gradient step is taken (Lemma~\ref{lem:dphiik0});
\item[(ii)] A Newton-CG step with negative curvature direction is taken (Lemma~\ref{lem:dphiik+-nc});
\item[(iii)] A Newton-CG step with an approximate solution to the regularized Newton equation is taken (Lemma~\ref{lem:dphiik+-sol});
\item[(iv)] Algorithm~\ref{alg:meo} is invoked and returns a direction of negative curvature for \(S^k\nabla^2 f(x^k)S^k\),  followed by a negative curvature step (Lemma~\ref{lem:dmeo2}).
\end{itemize}

\begin{lemma}\label{lem:dphiik0}
Let Assumption~\ref{assum:problem} hold and suppose that \(I^{k\varepsilon}_0\neq \emptyset\) at iteration \(k\), and \(\|g^{k\ve}_{I^{k\varepsilon}_0}\| \neq 0\), so that the proximal gradient step is taken. Then we have \(j_k < +\infty\) and  
\[
\varphi(x^k) - \varphi(x^{k+1}) \geq c^1_{pg}\varepsilon_g^{3/2},
\]
where \(c^1_{pg}\) is defined as in Lemma~\ref{lem:dphiik0_hpgnc}. 
\end{lemma}
\begin{proof}
By \cite[Remark 10.13]{B17}, the line search condition~\eqref{eq:lsprox} holds for any \(j_k\) such that \(\beta^{j_k} \geq \frac{L_g}{2(1 - \bar{\eta})}\), which implies \(\beta^{j_k - 1} < \frac{L_g}{2(1 - \bar{\eta})}\), as well as \(j_k < +\infty\) and \(t_{k+1} = \beta^{j_k} < \max\{\frac{\beta L_g}{2(1 - \bar{\eta})}, 1\}\). Next, we show that \(\|\mathcal{G}_{t_{k+1}}(x^k)\|^2 \geq \veg^{3/2}\). 

From the definition of \(g_i^{\varepsilon}\), we have \(\vert (\nabla f(x^k))_i\vert > \lambda + \veg^{3/4}\) for any \(i \!\in\! I^{k\varepsilon}_0\) such that \(g^{k\varepsilon}_i \!\neq\! 0\).  

\noindent \textbf{Case a).} Suppose that there exists  \(i \in I^{k\varepsilon}_0\) such that \((\nabla f(x^k))_i > \lambda + \veg^{3/4}\). \begin{itemize}
\item If \(x_i^k > \frac{1}{t_{k+1}}((\nabla f(x^k))_i + \lambda)\), then by~\eqref{eq:dxx}, we have \(\mathcal{G}_{t_{k+1}}(x^k)_i^2 = ((\nabla f(x^k))_i + \lambda)^2 >  (2\lambda + \veg^{3/4})^2 > \veg^{3/2}\).

\item If \(\vert x^k_i - \frac{1}{t_{k+1}}(\nabla f(x^k))_i\vert \leq \frac{\lambda}{t_{k+1}}\), then \(x^k_i \geq \frac{1}{t_{k+1}}((\nabla f(x^k))_i - \lambda) > \frac{1}{t_{k+1}}\veg^{3/4}\). Hence, by~\eqref{eq:dxx}, we have \(\mathcal{G}_{t_{k+1}}(x^k)_i^2= t_{k+1}^2(x^k_i)^2 \geq \veg^{3/2}\).

\item If \(x^k_i < \frac{1}{t_{k+1}}((\nabla f(x^k))_i - \lambda)\), then by~\eqref{eq:dxx}, we have \(\mathcal{G}_{t_{k+1}}(x^k)_i^2 = ((\nabla f(x^k))_i - \lambda)^2 >\veg^{3/2}\).
\end{itemize}

\noindent \textbf{Case b).} Suppose that there exists  \(i \in I^{k\varepsilon}_0\) such that \((\nabla f(x^k))_i < -(\lambda + \veg^{3/4})\).
\begin{itemize}
\item If \(x_i^k > \frac{1}{t_{k+1}}((\nabla f(x^k))_i + \lambda)\), then by~\eqref{eq:dxx}, we have \(\mathcal{G}_{t_{k+1}}(x^k)_i^2 = ((\nabla f(x^k))_i + \lambda)^2 > \veg^{3/2}\).

\item If \(\vert x^k_i - \frac{1}{t_{k+1}}(\nabla f(x^k))_i\vert \leq \frac{\lambda}{t_{k+1}}\), then \(x^k_i \leq \frac{1}{t_{k+1}}((\nabla f(x^k))_i +\lambda) < -\frac{1}{t_{k+1}}\veg^{3/4}\). Hence, by~\eqref{eq:dxx}, we have \(\mathcal{G}_{t_{k+1}}(x^k)_i^2 = t_{k+1}^2(x^k_i)^2 \geq \veg^{3/2}\).

\item If \(x^k_i < \frac{1}{t_{k+1}}((\nabla f(x^k))_i - \lambda)\), then by~\eqref{eq:dxx}, we have \(\mathcal{G}_{t_{k+1}}(x^k)_i^2 = ((\nabla f(x^k))_i - \lambda)^2 > (2\lambda + \veg^{3/4})^2 \geq \veg^{3/2}\).
\end{itemize}
This completes the proof.  
\end{proof}

The following property holds when the Newton-CG step is invoked. 
\begin{lemma}\label{lem:dphi2}
Suppose that a Newton-CG step is invoked. Let \((t^k_+, t^k_-)\) be defined as in~\eqref{eq:tk+-0}. Then, for any \(t \in (0, \min\{t^k_+, t^k_-\}]\), we have 
\[
\varphi(x^{k} + td^k) \leq \varphi(x^k) + t(g^k_{\neq 0\varepsilon})^\top d^k_{\neq 0\varepsilon} + \frac{1}{2}t^2(d^k_{\neq 0\varepsilon})^\top H^k_{\neq 0\varepsilon}d^k_{\neq 0\varepsilon} + \frac{L_H}{6}t^3\|d^k_{\neq 0\varepsilon}\|^3. 
\] 
\end{lemma}
\begin{proof}
From~\eqref{eq:p3}, for any \(t > 0\), we have 
\begin{align*}
\varphi(x^k \!+\! td^k) \leq & \varphi(x^k) \!+\! t\nabla f(x^k)^\top d^k \!+\! \frac{t^2}{2}(d^k)^\top \nabla^2 f(x^k)d^k \!+\! \frac{L_H}{6}t^3\|d^k\|^3  \!+\! \lambda(\|x^k \!+\! td^k\|_1 \!-\! \|x^k\|_1) \\
=&\varphi(x^k) + t(g^k_{\neq 0\varepsilon})^\top d^k_{\neq 0\varepsilon} + \frac{t^2}{2}(d^k_{\neq 0\varepsilon})^\top H^k_{\neq 0\varepsilon}d^k_{\neq 0\varepsilon} + \frac{L_H}{6}t^3\|d^k_{\neq 0\varepsilon}\|^3 \\
&+ t((\nabla f(x^k))_{\neq 0\varepsilon} - g^k_{\neq 0\varepsilon})^\top d^k_{\neq 0\varepsilon} + \lambda\sum_{i\in I^{k\varepsilon}_{\neq 0}}(\vert x^k_i + td^k_i\vert - \vert x^k_i\vert).  
\end{align*}

On the one hand, from the definitions of \(t^k_+\) and \(t^k_-\), we have
\[
x_i^k + td_i^k \geq 0, \text{if}~x^k_i > \veg^{1/2} \quad\text{and}\quad x_i^k + td_i^k \leq 0, {\rm if}~x^k_i < -\veg^{1/2}, \quad \forall t \in (0, \min\{t^k_+, t^k_-\}].
\]
Therefore,  
\[
\sum_{i\in I^{k\varepsilon}_{\neq 0}}(\vert x^k_i + td^k_i\vert - \vert x^k_i\vert) = t\sum_{i\in I^{k\varepsilon}_{\neq 0}}{\rm sgn}(x^k_i)d^k_i, \quad \forall t \in (0, \min\{t^k_+, t^k_-\}]. 
\]

On the other hand, from the definition of \(g^k_{\neq 0\varepsilon}\), we have 
\begin{align*}
((\nabla f(x^k))_{\neq 0\varepsilon} - g^k_{\neq 0\varepsilon})^\top d^k_{\neq 0\varepsilon} = -\lambda \sum_{i\in I^{k\varepsilon}_{\neq 0}}{\rm sgn}(x^k_i)d^k_i.
\end{align*}
Therefore, for any \(t \in (0, \min\{t^k_+, t^k_-\}]\), the desired result follows. 
\end{proof}

\begin{lemma}\label{lem:dphiik+-nc}
Let Assumption~\ref{assum:problem} hold and suppose that at iteration \(k\), a Newton-CG step is invoked and Capped CG returns \(d_{\rm type} =\) `NC'.
Then we have \(j_k < +\infty\) and 
\[
\varphi(x^k) - \varphi(x^{k + 1}) \geq  c_{nc}\min\{\varepsilon_g\veh, \veh^3\},
\]
where \(c_{nc} = \eta\theta^2\min\{1, \frac{9(1 - 2\eta)^2}{L_H^2}\}\).
\end{lemma}
\begin{proof}
From Lemma~\ref{lem:dphi2} and~\eqref{eq:pdnc}, for any \(t \in (0, \min\{t^k_+, t^k_-\}]\), we have 
\[
\varphi(x^k + td^k) \leq \varphi(x^k) - \frac{t^2}{2}(1 - \frac{L_H}{3}t)\|d^k_{\neq 0\varepsilon}\|^3.
\]
Hence, for \(0 < t < \min\{\min\{t^k_+, t^k_-\}, \frac{3(1 - 2\eta)}{L_H}\}\), we have 
\[
\varphi(x^k + td^k) < \varphi(x^k) - \eta t^2\|d^k_{\neq 0\varepsilon}\|^3 \leq \varphi(x^k) - \eta t^2\veh\|d^k\|^2,
\]
where the second inequality follows from \(\|d^k\| = \|d^k_{\neq 0\varepsilon}\| \geq \veh\). 
Therefore, the line search condition~\eqref{eq:lsreq} is satisfied whenever \(\theta^{j} < \min\{\min\{t^k_+, t^k_-\}, \frac{3(1 - 2\eta)}{L_H}\}\). Thus, \(j_k < +\infty\). By the line search rule, 
\[
\alpha_k = \theta^{j_k} \geq \min\{\theta\min\{\min\{t^k_+, t^k_-\}, \frac{3(1 - 2\eta)}{L_H}\}, 1\}. 
\]
From the definitions of \(t^k_+\) and \(t^k_-\), we have 
\begin{equation}\label{eq:mint+t-}
\min\{t^k_+, t^k_-\}\|d^k\| \geq \veg^{1/2}. 
\end{equation}
Since \(\|d^k\| = \|d^k_{\neq 0\varepsilon}\| = \frac{\vert (d^k_{\neq 0\varepsilon})^\top H^k_{\neq 0\varepsilon} d^k_{\neq 0\varepsilon}\vert}{\|d^k_{\neq 0\varepsilon}\|^2} \leq \|H^k_{\neq 0\varepsilon}\| \leq \|\nabla^2 f(x^k)\| \leq L_g\), it follows that   
\[
\alpha_k  \geq \min\{\theta\min\{\frac{\varepsilon_g^{1/2}}{L_g}, \frac{3(1 - 2\eta)}{L_H}\}, 1\}.
\]
Hence, \(j_k\) is finite. 

In addition, recalling that \(\|d^k\| \geq \veh\), we have  
\[
\alpha_k\|d^k\| \geq \min\{\theta\min\{\varepsilon_g^{1/2}, \frac{3(1 - 2\eta)\|d^k\|}{L_H}\}, \|d^k\|\} \geq \min\{\theta\min\{\varepsilon_g^{1/2}\veh^{-1}, \frac{3(1 - 2\eta)}{L_H}\}, 1\}\veh.
\]
Consequently,  
\[
\eta\!\alpha_k^2\veh\|d^k\|^2 \!\geq\! \eta\!\theta^2\!\min\{\veg\!\veh, \min\{\frac{9(1 \!-\! 2\eta)^2}{L_H^2}, \frac{1}{\theta^2}\}\veh^3\}\!\geq\!\eta\!\theta^2\!\min\{1, \frac{9(1 \!-\! 2\eta)^2}{L_H^2}\}\!\min\{\veg\!\veh, \veh^3\}.
\]
The desired result follows. 
\end{proof}

\begin{lemma}\label{lem:dphiik+-sol}
Let Assumption~\ref{assum:problem} hold for Problem~\eqref{eq:l1normcom}, and suppose that \(\veh \geq \veg^{1/2}\) in Algorithm~\ref{alg:pncg}. Suppose that, at iteration \(k\), a Newton-CG step is invoked and the Capped CG method returns \(d_{\rm type} =\) `SOL'. Then \(j_k < +\infty\) and 
\begin{equation*}
\varphi(x^k) - \varphi(x^{k+1}) \geq c_{sol}\min\{\|g^{k+1}_{\neq 0\varepsilon}\|^2\veh^{-1}, \veh^3, \veg\veh\},
\end{equation*}
where \(c_{sol} = \eta\min\{(\frac{4}{\sqrt{(\zeta+ 4\hat{\tau})^2 + 8L_H} + (\zeta+ 4\hat{\tau})})^2, \frac{9(1 - \zeta - 2\eta)^2\theta^2}{L_H^2}, \theta^2, \frac{(1 - \zeta)^2\theta^2 }{4\max\{L_H/3, 2\eta\}^2}\}\).
\end{lemma}
\begin{proof}
Define 
\begin{align*}
l_k \!\triangleq& \min\{l\in\mathbb{N}\mid x_i^k + \theta^ld_i^k > 0,~{\rm if}~i\in I^{k\varepsilon}_+ ~{\rm and}~x_i^k + \theta^ld_i^k < 0,~{\rm if}~i\in I^{k\varepsilon}_-\};\\
s_k\! \triangleq& \min\{s \!\geq\! l_k, s\!\in\!\mathbb{N}\!\mid \theta^s(g^k_{\neq 0\varepsilon})^\top \!d^k_{\neq 0\varepsilon} \!+\! \frac{\theta^{2s}}{2}(d^k_{\neq 0\varepsilon})^\top\! H^k_{\neq 0\varepsilon}d^k_{\neq 0\varepsilon} \!+\! \frac{L_H}{6}\theta^{3s}\|d^k_{\neq 0\varepsilon}\|^3 \!<\! -\eta\theta^{2s}\veh\|d^k\|^2\}.
\end{align*}
According to the definitions of \(t^{k}_+\), \(t^{k}_-\), and \(d^k_{I^{k\varepsilon}_0} = 0\), we have \(\theta^{l_k} < \min\{t^{k}_+, t^{k}_-\}\leq \theta^{l_k - 1}\). Hence, \(l_k\) is well-defined. For any \(s \geq l_k\), Lemma~\ref{lem:dphi2} gives 
\begin{equation} \label{eq:dvarphicopy}
\varphi(x^{k} + \theta^sd^k) - \varphi(x^k) \leq  \theta^s(g^k_{\neq 0\varepsilon})^\top d^k_{\neq 0\varepsilon} + \frac{1}{2}\theta^{2s}(d^k_{\neq 0\varepsilon})^\top H^k_{\neq 0\varepsilon}d^k_{\neq 0\varepsilon} + \frac{L_H}{6}\theta^{3s}\|d^k_{\neq 0\varepsilon}\|^3. 
\end{equation} 
Moreover, 
\begin{align}
&\theta^s(g^k_{\neq 0\varepsilon})^\top \!d^k_{\neq 0\varepsilon} \!+\! \frac{\theta^{2s}}{2}(d^k_{\neq 0\varepsilon})^\top\! H^k_{\neq 0\varepsilon}d^k_{\neq 0\varepsilon} \!+\! \frac{L_H}{6}\theta^{3s}\|d^k_{\neq 0\varepsilon}\|^3 \nonumber \\ 
=\!& -\theta^s((H^k_{\neq 0\varepsilon} + \tau_k\|g^k_{\neq 0\varepsilon}\|^{\delta}I)d^k_{\neq 0\varepsilon} - \hat{r}_k)^\top d^k_{\neq 0\varepsilon} + \frac{\theta^{2s}}{2}(d^k_{\neq 0\varepsilon})^\top H^k_{\neq 0\varepsilon}d^k_{\neq 0\varepsilon}  + \frac{L_H}{6}\theta^{3s}\|d^k_{\neq 0\varepsilon}\|^3 \nonumber \\
=\!& -\theta^s(1 - \frac{\theta^s}{2})(d^k_{\neq 0\varepsilon})^\top(H^k_{\neq 0\varepsilon} + \tau_k\|g^k_{\neq 0\varepsilon}\|^{\delta}I)d^k_{\neq 0\varepsilon} - \frac{\theta^{2s}}{2}\tau_k\|g^k_{\neq 0\varepsilon}\|^{\delta}\|d^k_{\neq 0\varepsilon}\|^2 \nonumber\\
\!&+ \theta^s\hat{r}_k^\top d^k_{\neq 0\varepsilon} + \frac{L_H}{6}\theta^{3s}\|d^k_{\neq 0\varepsilon}\|^3\nonumber\\
\overset{\eqref{eq:dsolacopy}}{\leq}\!& -\theta^s(1 - \frac{\theta^s}{2})\veh\|d^k_{\neq 0\varepsilon}\|^2 + \theta^s\|\hat{r}_k\| \|d^k_{\neq 0\varepsilon}\| + \frac{L_H}{6}\theta^{3s}\|d^k_{\neq 0\varepsilon}\|^3\nonumber\\
\overset{\eqref{eq:dsolccopy}}{\leq}\!&   - \frac{\theta^s}{2}\veh\|d^k_{\neq 0\varepsilon}\|^2 + \frac{\theta^s}{2}\zeta\veh\|d^k_{\neq 0\varepsilon}\|^2 + \frac{L_H}{6}\theta^{3s}\|d^k_{\neq 0\varepsilon}\|^3\nonumber\\
=\!&- \frac{\theta^s}{2}(1 - \zeta)\veh\|d^k_{\neq 0\varepsilon}\|^2 + \frac{L_H}{6}\theta^{3s}\|d^k_{\neq 0\varepsilon}\|^3 \label{eq:ddvarphi}\\
\overset{\eqref{eq:dsolbcopy}}{\leq}\!&  - \frac{\theta^s}{2}(1 - \zeta)\veh\|d^k_{\neq 0\varepsilon}\|^2 + \frac{1.1L_H}{6}\theta^{3s}\veh^{-1}\|g^k_{\neq 0\varepsilon}\|\|d^k_{\neq 0\varepsilon}\|^2\nonumber\\
\leq\!&  - \frac{\theta^{2s}}{2}(1 - \zeta)\veh\|d^k_{\neq 0\varepsilon}\|^2 + \frac{1.1L_H(U_g + \lambda\sqrt{n})}{6}\theta^{3s}\veh^{-1}\|d^k_{\neq 0\varepsilon}\|^2, \nonumber
\end{align}
where the last inequality follows from \(\theta \in(0, 1)\) and \(\|g^k_{\neq 0\varepsilon}\| \leq \|(\nabla f(x^k))_{I^{k\varepsilon}_{\neq 0}}\| + \lambda\cdot\|{\rm sgn}((x^k)_{I^{k\varepsilon}_{\neq 0}})\| \leq U_g + \lambda\sqrt{n}\). 

For any \(\bar{s}\in\mathbb{N}\) satisfying \(\bar{s} \geq l_k\) and \(\theta^{\bar{s}} < \frac{3(1-\zeta - 2\eta)\veh^2}{1.1L_H(U_g + \lambda\sqrt{n})}\), the preceding inequality implies  
\[
- \frac{\theta^{2\bar{s}}}{2}(1 - \zeta)\veh\|d^k_{\neq 0\varepsilon}\|^2 + \frac{1.1L_H(U_g + \lambda\sqrt{n})}{6}\theta^{3\bar{s}}\veh^{-1}\|d^k_{\neq 0\varepsilon}\|^2 < -\eta\theta^{2\bar{s}}\veh\|d^k\|^2. 
\]
Hence, \(s_k\) is well-defined and \(j_k \leq s_k < +\infty\). 

\noindent\textbf{Case a).} Suppose that \(s_k = l_k = 0\), then \(j_k = 0\) and \(\alpha_k = \theta^{j_k} = 1\).  Moreover, \({\rm sgn}(x^{k+1}_{I^{k\varepsilon}_{\neq 0}}) = {\rm sgn}(x^k_{I^{k\varepsilon}_{\neq 0}})\). Let \(\hat{g}^{k+1}_{\neq 0\varepsilon} := (g(x^{k+1}))_{I^{k\varepsilon}_{\neq 0}}\). Using the definitions of \(g^{k}_{\neq 0\varepsilon}\) and \(\hat{r}^k\), we have  
\begin{align*}
\|\hat{g}^{k+1}_{\neq 0\varepsilon}\| =& \|(\nabla f(x^k + \alpha_kd^k))_{I^{k\varepsilon}_{\neq0}} - (\nabla f(x^k))_{I^{k\varepsilon}_{\neq0}}  + (\nabla f(x^k))_{I^{k\varepsilon}_{\neq0}} + \lambda\cdot{\rm sgn}((x^k)_{I^{k\varepsilon}_{\neq0}})\| \nonumber \\
=& \|(\nabla f(x^k + \alpha_kd^k))_{I^{k\varepsilon}_{\neq0}} - (\nabla f(x^k))_{I^{k\varepsilon}_{\neq0}} - H^k_{\neq 0\varepsilon}d^k_{\neq 0\varepsilon} - \tau_k\|g^k_{\neq 0\varepsilon}\|^{\delta}d^k_{\neq 0\varepsilon} + \hat{r}^k\|\nonumber \\
\leq& \|(\nabla f(x^k + \alpha_kd^k))_{I^{k\varepsilon}_{\neq0}} - (\nabla f(x^k))_{I^{k\varepsilon}_{\neq0}} - H^k_{\neq 0\varepsilon}d^k_{\neq 0\varepsilon}\|  + (\frac{1}{2}\zeta + 2\hat{\tau})\veh\|d^k_{\neq 0\varepsilon}\|,
\end{align*}
where the inequality follows from \eqref{eq:dsolccopy} and the fact \(\tau_k \|g^k_{\neq 0\varepsilon}\|^{\delta} \leq 2\hat{\tau}\veh\). Let \(E_k\in\mathbb{R}^{\vert I_{\neq 0\varepsilon}^k\vert \times n}\) be the projection matrix such that \(E_k z = z_{I_{\neq 0}^{k\varepsilon}}\) for any \(z\in\R^n\). Since \(\alpha_k = 1\), 
\begin{align*}
&\|(\nabla f(x^k + d^k))_{I^{k\varepsilon}_{\neq0}} - (\nabla f(x^k))_{I^{k\varepsilon}_{\neq0}} -H^k_{\neq 0\varepsilon}d^k_{\neq 0\varepsilon}\| \\
=&\|E_k \nabla f(x^k + \alpha_kd^k) - E_k\nabla f(x^k) - E_k\nabla^2 f(x^k)E_k^\top E_kd^k\|\\
\leq & \|E_k\|\|\nabla f(x^k + \alpha_kd^k) - \nabla f(x^k) - \nabla^2 f(x^k)d^k\|\\
\leq& \frac{L_H}{2}\|d^k\|^2 = \frac{L_H}{2}\|d_{\neq 0\varepsilon}^k\|^2,
\end{align*}
where the first inequality follows from \(E_k^\top E_kd^k = d^k\), since \(d^k_i = 0\) for all \(i \in I^{k\varepsilon}_0\), and the second inequality follows from~\eqref{eq:p2}. Since the Newton-CG step is invoked at iteration \(k\), \(d^k_{I^{k\ve}_0} = 0\). Hence, \(I^{k\ve}_0 \subseteq I^{(k+1)\ve}_0\) and \( I^{(k+1)\ve}_{\neq 0} \subseteq I^{k\ve}_{\neq 0}\). By the definitions of \(g\) and \(g^{\ve}\), \(g^{\ve}_i(x^{k+1}) = g_i(x^{k+1})\), \(\forall i \in I^{(k+1)\ve}_{\neq 0}\). Therefore,  
\[
\|g^{k+1}_{\neq 0\varepsilon}\|^2 = \sum_{i\in I^{(k+1)\ve}_{\neq 0}}\vert g_i^\ve(x^{k+1})\vert^2 = \sum_{i\in I^{(k+1)\ve}_{\neq 0}}\vert g_i(x^{k+1})\vert^2 \leq \sum_{i\in I^{k\ve}_{\neq 0}}\vert g_i(x^{k+1})\vert^2 = \|\hat{g}^{k+1}_{\neq 0\ve}\|^2. 
\]
Consequently,  
\begin{equation}\label{eq:ngkp1}
\|g^{k+1}_{\neq 0\varepsilon}\| \leq \|\hat{g}^{k+1}_{\neq 0\varepsilon}\| \leq \frac{L_H}{2}\|d_{\neq 0\varepsilon}^k\|^2 +  (\frac{1}{2}\zeta + 2\hat{\tau})\veh\|d^k_{\neq 0\varepsilon}\|. 
\end{equation}

Applying the quadratic formula to~\eqref{eq:ngkp1}, we obtain 
\begin{align*}
\|d^k_{\neq 0\varepsilon}\| \geq&  (-(\frac{1}{2}\zeta + 2\hat{\tau}) + \sqrt{(\frac{1}{2}\zeta + 2\hat{\tau})^2 + 2L_H\|g^{k+1}_{\neq 0\varepsilon}\|/\veh^2})\frac{\veh}{L_H}\\
\geq& (-(\frac{1}{2}\zeta + 2\hat{\tau}) + \sqrt{(\frac{1}{2}\zeta + 2\hat{\tau})^2 + 2L_H\min\{\|g^{k+1}_{\neq 0\varepsilon}\|/\veh^2, 1\}})\frac{\veh}{L_H}\\
=&\frac{2\min\{\|g^{k+1}_{\neq 0\varepsilon}\|\veh^{-1}, \veh\}}{\sqrt{(\frac{1}{2}\zeta + 2\hat{\tau})^2 + 2L_H\min\{\|g^{k+1}_{\neq 0\varepsilon}\|/\veh^2, 1\}} + (\frac{1}{2}\zeta + 2\hat{\tau})}\\
\geq& \frac{4}{\sqrt{(\zeta+ 4\hat{\tau})^2 + 8L_H} + (\zeta+ 4\hat{\tau})}\min\{\|g^{k+1}_{\neq 0\varepsilon}\|\veh^{-1}, \veh\}.
\end{align*}
Therefore, by~\eqref{eq:lsreq}, 
\begin{equation}\label{eq:dvarphisol1}
\varphi(x^k) - \varphi(x^{k+1}) \geq \eta \alpha_k^2\veh\|d^k\|^2  \geq \hat{c}_{sol}\min\{\|g^{k+1}_{\neq 0\varepsilon}\|^2\veh^{-1}, \veh^3\},
\end{equation}
where \(\hat{c}_{sol} \!=\! \eta(\frac{4}{\sqrt{(\zeta+ 4\hat{\tau})^2 + 8L_H} + (\zeta+ 4\hat{\tau})})^2\). 

\noindent\textbf{Case b).} Suppose that \(s_k > l_k = 0\). If $j_k=0$, then \(\alpha_k = 1\). Since $l_k=0$, the sign-preservation property used in \textbf{Case a)} remains valid. Hence, the argument in \textbf{Case a)} applies verbatim and yields~\eqref{eq:dvarphisol1}. 

If $j_k > 0$, then \(\alpha_k < 1\) and the line search condition~\eqref{eq:lsreq} does not hold for \(j = 0\) or \(j = j_k - 1\). Thus, 
\[
-\eta\theta^{2j}\veh\|d^k\|^2 \leq \varphi(x^k + \theta^{j}d^k) - \varphi(x^k)\quad {\rm with}~j = 0~{\rm and}~j = j_k - 1.
\]
Combined this inequality with~\eqref{eq:dvarphicopy} and~\eqref{eq:ddvarphi}, we obtain 
\begin{subequations} 
\begin{align}
& -\eta\veh\|d^k\|^2 \leq -\frac{1}{2}(1-\zeta)\veh\|d^k\|^2 + \frac{L_H}{6}\|d^k\|^3; \quad {\rm and} \label{eq:j0} \\
&-\eta\theta^{2(j_k-1)}\veh\|d^k\|^2 \leq -\frac{\theta^{j_k - 1}}{2}(1-\zeta)\veh\|d^k\|^2 + \frac{L_H}{6}\theta^{3(j_k - 1)}\|d^k\|^3. \label{eq:jk-1}
\end{align}
\end{subequations}
From~\eqref{eq:j0}, we obtain \(\|d^k\| \geq \frac{6}{L_H}(\frac{1-\zeta}{2} - \eta)\veh\). Combining this with~\eqref{eq:jk-1}, we have 
\[
\frac{1 - \zeta}{2}\veh \leq \frac{L_H}{6}\theta^{2(j_k - 1)}\|d^k\| + \eta\veh\theta^{j_k-1} \leq (\frac{L_H}{6} + \frac{\eta\veh}{\|d^k\|})\theta^{j_k - 1}\|d^k\|\leq \frac{L_H(1 - \zeta)}{6(1 - \zeta - 2\eta)}\theta^{j_k - 1}\|d^k\|.
\]
Hence, we further have \(\theta^{j_k - 1}\|d^k\| \geq \frac{3(1 - \zeta - 2\eta)\veh}{L_H}\). By the line search condition~\eqref{eq:lsreq}, 
\begin{equation}\label{eq:dvarphisol2}
\varphi(x^k) - \varphi(x^{k+1}) \geq \eta \frac{9(1 - \zeta - 2\eta)^2\theta^2}{L_H^2}\veh^3.
\end{equation}

\noindent\textbf{Case c).} Suppose that \(s_k = l_k \geq 1\), then \(j_k \leq s_k = l_k\) and  \(\theta^{2j_k}\|d^k\|^2 \geq \theta^{2l_k}\|d^k\|^2 \geq (\theta\min\{t^{k}_+, t^{k}_-\}\|d^k\|)^2 \overset{\eqref{eq:mint+t-}}{\geq} \theta^2\veg\). 
Hence,
\begin{equation}\label{eq:dvarphisol3}
\varphi(x^k) - \varphi(x^{k+1}) \geq \eta\theta^{2j_k}\veh\|d^k\|^2 \geq \eta\theta^2\veg\veh.
\end{equation}

\noindent\textbf{Case d).} Suppose that \(s_k > l_k \geq 1\), then \(\min\{t^{k}_+, t^{k}_-\} \leq \theta^{l_k - 1} \leq 1\), which, together with~\eqref{eq:mint+t-}, yields \(\veg^{1/2} \leq \min\{t^{k}_+, t^{k}_-\}\|d_k\| \leq \|d_k\|\).  

If \(j_k \leq l_k\), then \(\theta^{j_k} \geq \theta^{l_k}\) and \(\theta^{l_k}\|d_k\| \geq \theta\min\{t^k_+, t^k_-\}\|d_k\| \geq \theta\veg^{1/2}\). It follows from~\eqref{eq:lsreq} that 
\begin{equation}\label{eq:eq:dvarphisoladd}
\varphi(x^k) - \varphi(x^{k+1}) \geq \eta\theta^{2j_k}\veh\|d_k\|^2 \geq  \eta\theta^{2l_k}\veh\|d_k\|^2 \geq \eta\theta^2\veg\veh. 
\end{equation}

If \(j_k \geq l_k + 1\), then \(\alpha_k < 1\). By an argument similar to that used in \textbf{Case b)}, it follows from~\eqref{eq:jk-1} that  
\[
\theta^{j_k}\|d^k\| \geq \frac{(1 - \zeta)\theta \veh}{L_H/3 + 2\eta \veh/\|d^k\|} \geq \frac{(1 - \zeta)\theta \veh}{L_H/3 + 2\eta \veh\veg^{-1/2}} \geq \frac{(1 - \zeta)\theta \veg^{1/2}}{2\max\{L_H/3, 2\eta\}},
\]
where the last inequality follows from \(\veh \geq \veg^{1/2}\). 
Hence, 
\begin{equation}\label{eq:dvarphisol4}
\varphi(x^k) - \varphi(x^{k+1}) \geq \eta\theta^{2j_k}\veh\|d^k\|^2 \geq  \frac{\eta(1 - \zeta)^2\theta^2 }{4\max\{L_H/3, 2\eta\}^2}\veg\veh.
\end{equation}
Combining~\eqref{eq:dvarphisol1},~\eqref{eq:dvarphisol2}--\eqref{eq:dvarphisol4}, the desired result follows. 
\end{proof}

We first state our result for the first-order iteration complexity. 
\begin{theorem}\label{th:foic}
Suppose that Assumption~\ref{assum:problem} holds for Problem~\eqref{eq:l1normcom}. Consider Algorithm~\ref{alg:pncg} with \(0 < \veg < 1\) and \(\veg^{1/2} \leq \veh < 1\). Define 
\begin{equation*}
K^{fo}_{pgn2cm} := \left\lceil \frac{3(\varphi(x^0) - \varphi_{\rm low})}{\min\{c_{nc}, c^1_{pg}, c_{sol}\}}\max\{\veg^{-1}\veh^{-1}, \veh^{-3}, \veg^{-2}\veh, \veg^{-3/2}\}\right\rceil + 2,
\end{equation*}
where \(c^1_{pg}\), \(c_{nc}\), and \(c_{sol}\) are the same as in Lemmas~\ref{lem:dphiik0},~\ref{lem:dphiik+-nc}, and~\ref{lem:dphiik+-sol}, respectively. Then Algorithm~\ref{alg:pncg} reaches a weak \(\varepsilon_g\)-1o point in at most \(K^{fo}_{pgn2cm}\) iterations. In particular, if \(\veg = \varepsilon\) and \(\veh = \varepsilon^{1/2}\), then Algorithm~\ref{alg:pncg} outputs a weak \(\varepsilon\)-1o point within 
 \[
\overline{K}^{fo}_{pgn2cm} := \left\lceil \frac{3(\varphi(x^0) - \varphi_{\rm low})}{\min\{c_{nc}, c^1_{pg}, c_{sol}\}}\varepsilon^{-3/2}\right\rceil + 2
\]
iterations. 
\end{theorem}

\begin{proof}
Suppose, for contradiction, that \(\|g^{k\varepsilon}\| > \veg\) for \(k = 0, 1, \ldots, K^{fo}_{pgn2cm} + 1\). Consider the following sets of iteration indices:
\[
\mathcal{K}_1 \cup \mathcal{K}_2 \cup \mathcal{K}_3  = \{0, 1, \cdots, K^{fo}_{pgn2cm}\},
\]
where \(\mathcal{K}_1\), \(\mathcal{K}_2\), and \(\mathcal{K}_3\) are defined as follows. 

\noindent\textbf{Case 1.} \(\mathcal{K}_1 := \left\{l = 0, 1, \cdots, K^{fo}_{pgn2cm}: {\rm the}~\textbf{Proximal~gradient~step}~{\rm is~invoked}\right\}\). It follows from Lemma~\ref{lem:dphiik0} that
\begin{align}\label{eq:dvarphif2}
\varphi(x^l) - \varphi(x^{l+1}) >&c^1_{pg}\veg^{3/2}, \quad \forall l\in\mathcal{K}_1. 
\end{align}

\noindent\textbf{Case 2.} \(\mathcal{K}_2 := \left\{l \!=\! 0, \ldots, K^{fo}_{pgn2cm}: {\rm the}~\textbf{Newton-CG step}~{\rm is~invoked~and}~\|g^{l+1}_{\neq 0\varepsilon}\| >\veg\right\}\). By  Lemmas~\ref{lem:dphiik+-nc} and~\ref{lem:dphiik+-sol}, for any \(l \in \mathcal{K}_2\), we have 
\begin{align}\label{eq:dvarphif3}
\varphi(x^l) - \varphi(x^{l+1}) >&\min\{c_{nc}, c_{sol}\}\min\{\veg\veh, \veh^3, \|g^{l+1}_{\neq0\varepsilon}\|^2\veh^{-1}\}\nonumber \\
\geq&\min\{c_{nc}, c_{sol}\}\min\{\veg\veh, \veh^3, \veg^2\veh^{-1}\}. 
\end{align}

\noindent\textbf{Case 3.} \(\mathcal{K}_3 := \left\{l = 0,\ldots, K^{fo}_{pgn2cm}: {\rm the}~\textbf{Newton-CG step}~{\rm is~invoked~but}~\|g^{l+1}_{\neq 0\varepsilon}\| \leq \veg\right\}\). By the contradiction hypothesis, \(\|g^{(l+1)\ve}\| > \veg\), whereas \(\|g^{l+1}_{\neq 0\varepsilon}\| \leq \veg\). Hence, \(I^{(l+1)\ve}_0 \neq \emptyset\) and \(\|g^{(l+1)\ve}_{I_0^{(l+1)\ve}}\| \neq 0\). Therefore, \(l + 1\in \mathcal{K}_1\) whenever \(l+1 \leq K^{fo}_{pgn2cm}\), and consequently \(\vert \mathcal{K}_3\vert \leq \vert\mathcal{K}_1\vert + 1\). 

Using the monotonicity of \(\{\varphi(x^l)\}_{l\in\mathbb{N}}\) and \eqref{eq:dvarphif2}--\eqref{eq:dvarphif3}, we obtain
\begin{equation}
\begin{split}\label{eq:dvarphikf}
&\varphi(x^0) - \varphi_{\rm low} \geq  \sum_{l \in \mathcal{K}_1}(\varphi(x^l) - \varphi(x^{l+1})) \geq c^1_{pg}\varepsilon_g^{3/2}\vert\mathcal{K}_1\vert; \\
&\varphi(x^0) - \varphi_{\rm low} \geq \sum_{l \in \mathcal{K}_2}(\varphi(x^l) - \varphi(x^{l+1}))  \geq \min\{c_{nc}, c_{sol}\}\min\{\veg\veh, \veh^3, \veg^2\veh^{-1}\}\vert\mathcal{K}_2\vert.  
\end{split}
\end{equation}
From~\eqref{eq:dvarphikf}, we have 
\begin{align*}
&\vert\mathcal{K}_1\vert \leq \frac{\varphi(x^0) - \varphi_{\rm low}}{c^1_{pg}}\veg^{-3/2},~~~~\vert\mathcal{K}_2\vert  \leq \frac{\varphi(x^0) - \varphi_{\rm low}}{\min\{c_{nc}, c_{sol}\}}\max\{\veg^{-1}\veh^{-1}, \veh^{-3}, \veg^{-2}\veh\}.
\end{align*}
Hence, we have 
\begin{align*}
K^{fo}_{pgn2cm} + 1 =& \sum_{i=1}^3\vert \mathcal{K}_i\vert \leq 2\vert \mathcal{K}_1\vert + \vert \mathcal{K}_2\vert + 1\leq 3\max\{\vert \mathcal{K}_1\vert, \vert \mathcal{K}_2\vert\} + 1\\
\leq& \frac{3(\varphi(x^0) - \varphi_{\rm low})}{\min\{c_{nc}, c^1_{pg}, c_{sol}\}}\max\{\veg^{-1}\veh^{-1}, \veh^{-3}, \veg^{-2}\veh, \veg^{-3/2}\} + 1 \leq K^{fo}_{pgn2cm} - 1,
\end{align*}
which yields the required contradiction. 
\end{proof}
If Problem~\eqref{eq:l1normcom} is an \(\ell_1\)-regularized convex composite optimization problem, then the Capped CG method returns an approximate solution to the regularized Newton equation whenever a Newton-CG step is invoked. It follows from Theorem~\ref{th:foic} that Algorithm~\ref{alg:pncg} finds a weak \(\veg\)-1o point within at most \(\lceil \frac{3(\varphi(x^0) - \varphi_{\rm low})}{\min\{c^1_{pg}, c_{sol}\}}\veg^{-3/2}\rceil\) iterations when \(\veh = \veg^{1/2}\) in the Newton-CG step. Similarly to Lemma~\ref{lem:dmeo}, the following result holds when Algorithm~\ref{alg:meo} is invoked by Algorithm~\ref{alg:pncg}.

\begin{lemma}\label{lem:dmeo2}
Let Assumption~\ref{assum:problem} hold and suppose that at iteration \(k\), Algorithm~\ref{alg:meo} is invoked by Algorithm~\ref{alg:pncg} and identifies a direction with curvature less than or equal to \(-\frac{1}{2}\varepsilon_h\). Then we have 
\[
\varphi(x^k) - \varphi(x^{k+1}) > \frac{1}{8}c_{nc}\min\{\veg\veh, \veh^3\},
\]
where \(c_{nc} = \eta\theta^2\min\{1, \frac{9(1 - 2\eta)^2}{L_H^2}\}\). 
\end{lemma}

We are now ready to state our main result for the second-order iteration complexity. 

\begin{theorem}\label{th:soic}
Suppose that Assumption~\ref{assum:problem} holds for Problem~\eqref{eq:l1normcom}. Consider Algorithm~\ref{alg:pncg} with \(0 < \veg < 1\) and \(\veg^{1/2} \leq \veh < 1\). Define 
\begin{equation*}
K^{so}_{pgn2cm} := \left\lceil \frac{5(\varphi(x^0) - \varphi_{\rm low})}{\min\{c_{nc}/8, c^1_{pg}, c_{sol}\}}\max\{\veg^{-1}\veh^{-1}, \veh^{-3}, \veg^{-2}\veh, \veg^{-3/2}\}\right\rceil + 2,
\end{equation*}
where \(c^1_{pg}\), \(c_{nc}\), and \(c_{sol}\) are the same as in Lemmas~\ref{lem:dphiik0},~\ref{lem:dphiik+-nc}, and~\ref{lem:dphiik+-sol}, respectively. Then with probability at least \((1 - \sigma)^{K^{so}_{pgn2cm}}\), Algorithm~\ref{alg:pncg} terminates at a weak \((\varepsilon_g, \varepsilon_h)\)-2o point within at most \(K^{so}_{pgn2cm}\) iterations, where \(\sigma \in [0, 1)\) is the probability of failure of Algorithm~\ref{alg:meo}. In particular, if \(\veg = \varepsilon\) and \(\veh = \varepsilon^{1/2}\), then Algorithm~\ref{alg:pncg} outputs a weak \((\varepsilon, \varepsilon^{\frac{1}{2}})\)-2o point with probability at least \((1 - \sigma)^{\overline{K}^{so}_{pgn2cm}}\) within 
 \[
\overline{K}^{so}_{pgn2cm} := \left\lceil \frac{5(\varphi(x^0) - \varphi_{\rm low})}{\min\{c_{nc}/8, c^1_{pg}, c_{sol}\}}\varepsilon^{-3/2}\right\rceil + 2
\]
iterations. 
\end{theorem}
\begin{proof}
Suppose, for contradiction, that Algorithm~\ref{alg:pncg} runs for \(K^{so}_{pgn2cm}\) iterations without terminating. Consider the following sets of iteration indices:
\[
\mathcal{K}_1 \cup \mathcal{K}_2 \cup \mathcal{K}_3 \cup \mathcal{K}_4 = \{0, 1, \cdots, K^{so}_{pgn2cm}- 1\},
\]
where \(\mathcal{K}_1\), \(\mathcal{K}_2\), \(\mathcal{K}_3\), and \(\mathcal{K}_4\) are defined as follows. 

\smallskip 

\noindent\textbf{Case 1.} \(\mathcal{K}_1 := \left\{l = 0, 1, \cdots, K^{so}_{pgn2cm} - 1: {\rm the}~\textbf{Proximal~gradient~step}~{\rm is~invoked}\right\}\). Similar to \textbf{Case 1} in the proof of Theorem~\ref{th:foic},~\eqref{eq:dvarphif2} holds. 

\smallskip 

\noindent\textbf{Case 2.} \(\mathcal{K}_2 := \left\{l = 0, 1, \cdots, K^{so}_{pgn2cm} - 1: {\rm the}~\textbf{Newton-CG step}~{\rm is~invoked~and}\right.\)
\(\left.\|g^{l+1}_{\neq 0\varepsilon}\| >\veg\right\}\). Similar to \textbf{Case 2} in the proof of Theorem~\ref{th:foic},~\eqref{eq:dvarphif3} holds. 

\smallskip 

\noindent\textbf{Case 3.} \(\mathcal{K}_3 := \left\{l = 0, 1, \cdots, K^{so}_{pgn2cm} - 1: {\rm the}~\textbf{Newton-CG step}~{\rm is~invoked~but}\right.\)
\(\left.\|g^{l+1}_{\neq 0\varepsilon}\| \leq \veg\right\}\). For any \(l \in\mathcal{K}_3\), we have \(\varphi(x^{l+1}) \leq \varphi(x^l)\) and \(l + 1\in \mathcal{K}_1\) or \(l + 1\in \mathcal{K}_4\) whenever \(l + 1\leq K^{so}_{pgn2cm} - 1\). Therefore, \(\vert \mathcal{K}_3\vert \leq \vert\mathcal{K}_1\vert + \vert\mathcal{K}_4\vert + 1\). 

\smallskip 

\noindent\textbf{Case 4.} \(\mathcal{K}_4 := \left\{l = 0, 1, \cdots, K^{so}_{pgn2cm} - 1: {\rm Algorithm~\ref{alg:meo}~is~ invoked~at~the}~l\text{-}th~{\rm step} \right\}\). At each iteration \(l \in \mathcal{K}_4\), Algorithm~\ref{alg:pncg} calls Algorithm~\ref{alg:meo}, which returns a direction of negative curvature for \(S^l_{\neq 0}(\nabla^2 f(x^l))_{\neq 0}S^l_{\neq 0}\) since \(l < K^{so}_{pgn2cm}\) and the algorithm continues to iterate. By Lemma~\ref{lem:dmeo2}, we have 
\begin{equation}\label{eq:dvarphi1}
\varphi(x^l) - \varphi(x^{l+1}) > \frac{1}{8}c_{nc}\min\{\veg\veh, \veh^3\}, \quad \forall l\in\mathcal{K}_4.
\end{equation}

Using the monotonicity of \(\{\varphi(x^l)\}_{l\in\mathbb{N}}\),~\eqref{eq:dvarphif2},~\eqref{eq:dvarphif3}, and \eqref{eq:dvarphi1}, we obtain 
\begin{equation}
\begin{split}\label{eq:dvarphik}
&\varphi(x^0) \!-\! \varphi_{\rm low} \!\geq\!\!\! \sum_{l \in \mathcal{K}_1}(\varphi(x^l) \!-\! \varphi(x^{l+1})) \!\geq\! c^1_{pg}\varepsilon_g^{3/2}\vert\mathcal{K}_1\vert; \\ 
&\varphi(x^0) \!-\! \varphi_{\rm low} \!\geq\!\!\! \sum_{l \in \mathcal{K}_2}(\varphi(x^l) \!-\! \varphi(x^{l+1}))  \!\geq\! \min\{c_{nc}, c_{sol}\}\!\min\{\veg\veh, \veh^3, \veg^2\veh^{-1}\}\!\vert\mathcal{K}_2\vert;\\ 
&\varphi(x^0) \!-\! \varphi_{\rm low} \!\geq\!\!\!  \sum_{l \in \mathcal{K}_4}(\varphi(x^l) - \varphi(x^{l+1})) \!\geq\! \frac{1}{8}c_{nc}\min\{\veg\veh, \veh^3\}\vert\mathcal{K}_4\vert.  
\end{split}
\end{equation}
From~\eqref{eq:dvarphik}, we have 
\begin{align*}
&\vert\mathcal{K}_1\vert \leq \frac{\varphi(x^0) - \varphi_{\rm low}}{c^1_{pg}}\veg^{-3/2},\quad \vert\mathcal{K}_2\vert  \leq \frac{\varphi(x^0) - \varphi_{\rm low}}{\min\{c_{nc}, c_{sol}\}}\max\{\veg^{-1}\veh^{-1}, \veh^{-3}, \veg^{-2}\veh\},\\
&\vert\mathcal{K}_4\vert  \leq \frac{8(\varphi(x^0) - \varphi_{\rm low})}{c_{nc}}\max\{\veg^{-1}\veh^{-1}, \veh^{-3}\}.
\end{align*}
Hence,  
\begin{align*}
K^{so}_{pgn2cm} =& \sum_{i=1}^4\vert \mathcal{K}_i\vert \leq 2\vert \mathcal{K}_1\vert + \vert \mathcal{K}_2\vert + 2\vert \mathcal{K}_4\vert + 1\leq 5\max\{\vert \mathcal{K}_1\vert, \vert \mathcal{K}_2\vert, \vert \mathcal{K}_4\vert\} + 1\\
\leq& \frac{5(\varphi(x^0) - \varphi_{\rm low})}{\min\{c_{nc}/8, c^1_{pg}, c_{sol}\}}\max\{\veg^{-1}\veh^{-1}, \veh^{-3}, \veg^{-2}\veh, \veg^{-3/2}\} + 1\leq K^{so}_{pgn2cm} - 1,
\end{align*}
which yields the required contradiction. 
 
The probability that Algorithm~\ref{alg:pncg} fails to identify a point satisfying \(\lambda_{\min}(S^k_{\neq 0}(\nabla^2 f(x^k))_{\neq 0}S^k_{\neq 0}) \geq -\veh\) within \(K^{so}_{pgn2cm}\) iterations is at most \( 1- (1 - \sigma)^{K^{so}_{pgn2cm}}\), by the same argument as in Theorem 3 in~\cite{XW23}. This completes the proof.
\end{proof}

We now give the operation complexity of Algorithm~\ref{alg:pncg} in terms of gradient evaluations and Hessian-vector products. In addition to the complexity analysis of MEO, the proof also relies on the complexity analysis about the Capped CG method in~\cite{RNW20}. 

\begin{corollary}\label{coro:operpgn2cm}
Suppose that the assumptions of Theorem~\ref{th:soic} and Assumption~\ref{assume:meo} hold. Let \(K^{so}_{pgn2cm}\) be defined as in Theorem~\ref{th:soic}. 
Suppose that the value of \(M\) used in the Capped CG method satisfies \(M \leq L_g\). 
Then, with probability at least \((1 - \sigma)^{K^{so}_{pgn2cm}}\), Algorithm~\ref{alg:pncg} terminates at a weak \((\veg, \veh)\)-2o point after at most
\[
\mathcal{O}(K^{so}_{pgn2cm}\min\{n, \veh^{-1/2}\ln(n\sigma^{-1}\veh^{-1})\})
\]
Hessian-vector products and/or gradient evaluations. The probability that the algorithm terminates incorrectly, outputting a weak $\varepsilon_g$-1o point that is not a weak $(\varepsilon_g, \varepsilon_h)$-2o point, is at most $1-(1-\sigma)^{K^{so}_{pgn2cm}}$.

For \(n\) sufficiently large, the following statements hold: 
\begin{itemize}
\item[(i)] With probability at least \((1 - \sigma)^{K^{so}_{pgn2cm}}\), the bound is 
\[
\widetilde{\mathcal{O}}(\max\{\veg^{-1}\veh^{-3/2}, \veh^{-7/2}, \veg^{-2}\veh^{1/2}, \veg^{-3/2}\veh^{-1/2}\}).
\]
\item[(ii)] If \(\veg \!=\! \varepsilon\) and \(\veh \!=\! \varepsilon^{1/2}\), then, with probability at least \((1 \!-\! \sigma)^{K^{so}_{pgn2cm}}\), the bound is \(\widetilde{\mathcal{O}}(\varepsilon^{-7/4})\).
\end{itemize}
\end{corollary}
\begin{proof}
According to \cite[Lemma 1, (8)]{RNW20}, the number of Hessian-vector products required by the Capped CG method at iteration \(k\) is bounded by \(2\min\{n, J_k(L_g, \veh, \zeta)\} + 1\), where 
\[
J_k(L_g, \veh, \zeta) \leq \min\{n, \mathcal{O}(\veh^{-1/2}\ln(\veh^{-1}))\}= \min\{n, \widetilde{\mathcal{O}}(\veh^{-1/2})\}. 
\]
Recalling the proof of Corollary~\ref{cor:hpgnc}, the number of Hessian-vector products required by MEO at iteration \(k\) is \(N_k^{\rm meo} \leq \min\{n, \mathcal{O}(\veh^{-1/2}\ln(n/\sigma))\}\). 

Therefore, using Theorem~\ref{th:soic}, the total number of Hessian-vector products and gradient evaluations required by Algorithm~\ref{alg:pncg} is bounded by
\begin{align*}
&\sum_{k=0}^{K^{so}_{pgn2cm} - 1}\max\{2\min\{n, J_k(L_g, \veh, \zeta)\} + 1, N_k^{\rm meo}\} \\
= &\mathcal{O}(K^{so}_{pgn2cm}\min\{n, \veh^{-1/2}\max\{\ln(\veh^{-1}), \ln(n/\sigma)\}\})\\
=&\mathcal{O}(K^{so}_{pgn2cm}\min\{n, \veh^{-1/2}\ln(n\sigma^{-1}\veh^{-1})\})\\
=&\mathcal{O}(\max\{\veg^{-1}\veh^{-1}, \veh^{-3}, \veg^{-2}\veh, \veg^{-3/2}\}\min\{n, \veh^{-1/2}\ln(n\sigma^{-1}\veh^{-1})\}).
\end{align*}

(i) For \(n\) sufficiently large, the above bound becomes 
\[
\mathcal{O}(\max\{\veg^{-1}\veh^{-3/2}, \veh^{-7/2}, \veg^{-2}\veh^{1/2}, \veg^{-3/2}\veh^{-1/2}\}\ln(n\sigma^{-1}\veh^{-1})),
\]
which yields the desired result. 

(ii) When \(\veg = \varepsilon\) and \(\veh = \varepsilon^{1/2}\), the bound in statement (i) becomes  \(\mathcal{O}(\varepsilon^{-7/4}\ln(n\sigma^{-1}\varepsilon^{-1/2}))\), which yields the desired result. 
\end{proof}

\section{Numerical experiments}\label{sec:num}

In this section, we test the performance of Algorithms~\ref{alg:hpgnc} and~\ref{alg:pncg} with \(\veh = \veg^{1/2}\). All numerical experiments were implemented in MATLAB R2024b and conducted on a computer equipped with an Intel(R) Core(TM) i9-10885U CPU @ 2.40GHz \(\times\) 2.4 and 32GB of RAM. %
The following settings are shared by all experiments involving the corresponding algorithm. For Algorithm~\ref{alg:hpgnc}, we set $\eta=10^{-4}$. For Algorithm~\ref{alg:pncg}, we set \(\eta = \min\{\frac{1}{4}(1 - \zeta), 10^{-4}\}\), where \(\zeta = 0.999\) in the Capped CG method (Algorithm~\ref{alg:ccg}), \(\tau_k = \frac{2\veh}{\min\{1, \|g^{k}_{\neq 0 \varepsilon}\|^{\delta}\}}\) with \(\delta = 1\) in the regularized Newton equation~\eqref{eq:rnecopycopy}. All other parameters are specified separately in the descriptions of the individual experiments.
When the MEO procedure is invoked, we implement it using the randomized Lanczos approach described in \cite[Algorithm 5]{RNW20}, with the oracle failure tolerance set to $\sigma=0.01$. Experiments involving only the first phase do not invoke the MEO procedure.

\subsection{A toy example}\label{subsec:toyexample}

The goal of this toy experiment is to illustrate the negative-curvature escape mechanism analyzed in Theorems~\ref{th:soicpgnc} and~\ref{th:soic}. Specifically, we show that Algorithm~\ref{alg:hpgnc} can escape from a strong \(\veg\)-1o point that is not a strong* \((\veg, \veh)\)-2o point, and that Algorithm~\ref{alg:pncg} can escape from a weak \(\veg\)-1o point that is not a weak \((\veg, \veh)\)-2o point.

Let \(f(x) = \sum_{i=1}^3f_i(x)\), where \(f_1(x) = (x_1 - 1)^2(x_1 - 3)^2 - (10^{-4} + \frac{1}{\sqrt{3}}\times 10^{-6})x_1\), \(f_2(x) = (x_2 + 1)^2(x_2 + 3)^2 + (10^{-4} + \frac{1}{\sqrt{3}}\times 10^{-6})x_2\), and \(f_3(x) = x_3^2(x_3 - 1)^2 - (10^{-4} - \frac{1}{\sqrt{3}}\times 10^{-6})x_3\). Consider the problem
\begin{equation}\label{eq:toyn3}
\min_{x\in\mathbb{R}^3} f(x) + 10^{-4} \|x\|_1.
\end{equation}
We note that $\varphi(x) := f(x) + 10^{-4}\|x\|_1$ is a coercive function, which ensures that for any $x^0 \in \mathbb{R}^3$, the level set $\mathcal{L}_{\varphi}(x^0) = \{x \in \mathbb{R}^3 \mid \varphi(x) \leq \varphi(x^0)\}$ is bounded. Since $f(x)$ is a quartic polynomial, it is twice uniformly Lipschitz continuously differentiable on any bounded set, thereby satisfying  Assumption~\ref{assum:problem}. Specifically, let \(R := \max_{x\in\mathcal{N}}\{\|x\|_2\}\), where \(\mathcal{N}\) is a compact neighborhood of $\mathcal{L}_{\varphi}(x^0)$. By calculating the second derivatives $f_1''(x_1) = 12x_1^2 - 48x_1 + 44$, $f_2''(x_2) = 12x_2^2 + 48x_2 + 44$, $f_3''(x_3) = 12x_3^2 - 12x_3 + 2$, and the third derivatives $f_1'''(x_1) = 24x_1 - 48$, $f_2'''(x_2) = 24x_2 + 48$, $f_3'''(x_3) = 24x_3 - 12$, the Lipschitz constants can be bounded by $L_g = 12R^2 + 48R + 44$ and $L_H = 24R+48$.  
We set \(\beta = 2\), \(\bar{\eta} = 0.5\), and \(\theta = 0.25\) in Algorithm~\ref{alg:hpgnc}. Since the definition of \(g(x)\) is numerically ill-conditioned, especially when \(x_i\) is close to zero for some \(i\), incorrect signs can propagate substantial errors in \((g(x))_i\). To address this instability, we compute \(g\) in Algorithm~\ref{alg:hpgnc} as follows: \(g_i = (\nabla f(x))_i + \lambda\) if \(x_i > 10^{-16}\); \(g_i = (\nabla f(x))_i - \lambda\) if \(x_i < -10^{-16}\); otherwise, \(g_i = (\nabla f(x))_i - \min\{\max\{-\lambda, (\nabla f(x))_i\}, \lambda\}\). We set \(\beta = 2\), \(\theta = 0.7\) in~\eqref{eq:lsreq}, \(\theta = 0.25\) in~\eqref{eq:lsreq2}. 

\begin{figure}[h!]
\begin{minipage}[t]{1\linewidth}
\centering
\includegraphics[width = 0.325\textwidth]{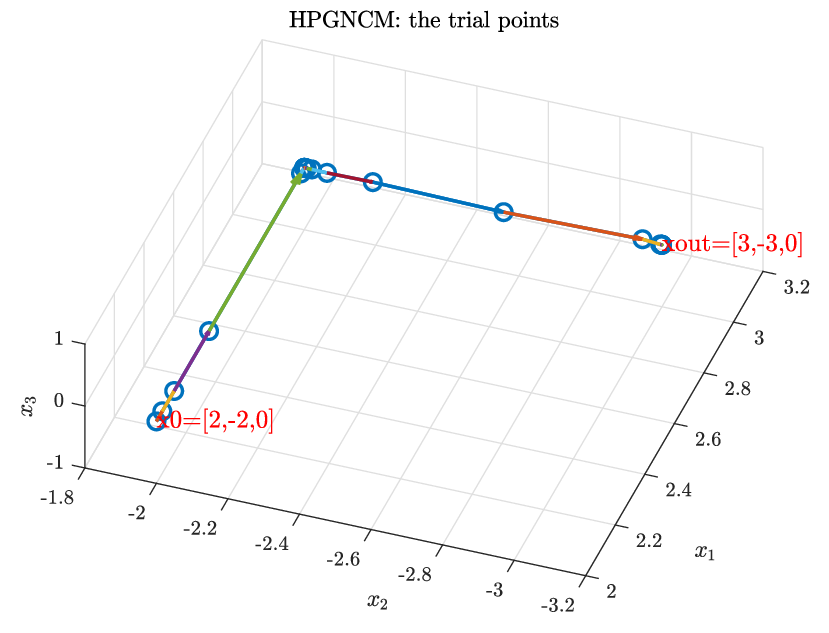}
\includegraphics[width = 0.325\textwidth]{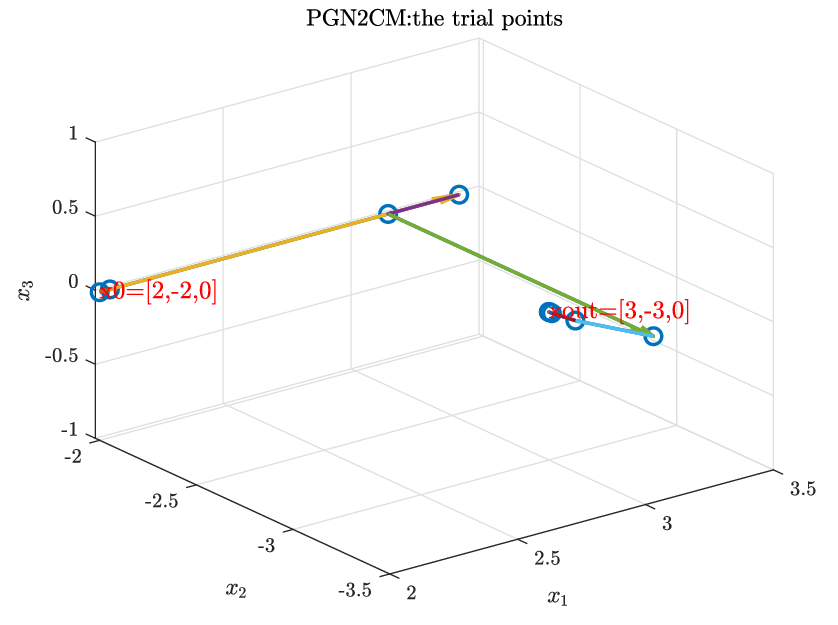}
\includegraphics[width = 0.325\textwidth]{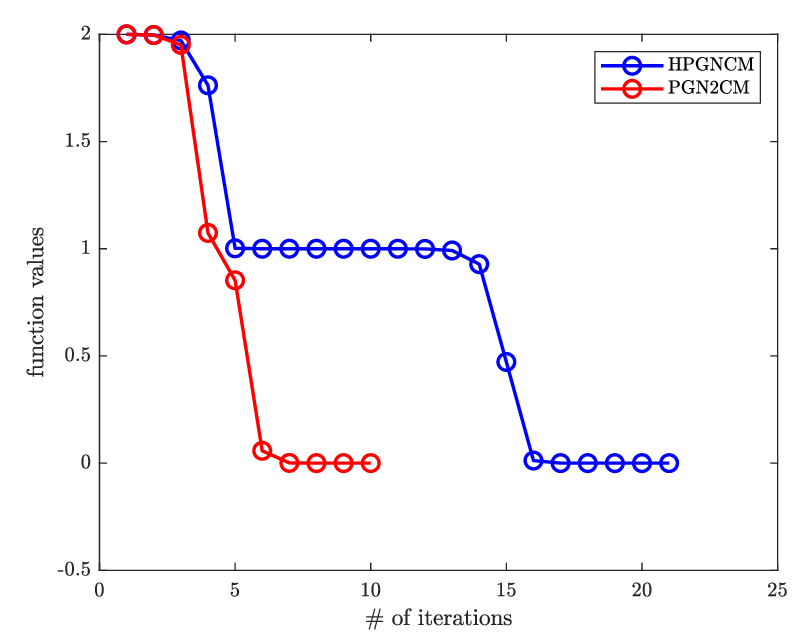}
\end{minipage}
\caption{Left: the trial points generated by Algorithm~\ref{alg:hpgnc}; middle: the trial points generated by Algorithm~\ref{alg:pncg}; right: the objective function value plotted against iteration.}\label{fig:toyexamplen3}
\end{figure}

Figure~\ref{fig:toyexamplen3} displays the trial points generated by Algorithms~\ref{alg:hpgnc} and~\ref{alg:pncg} from \(x^0 = (2, -2, 0)^\top\), along with the corresponding objective function \(\varphi(x^k)\) across iterations. It can be verified that \(g(x^0) = g^{\varepsilon}(x^0) = (-\frac{1}{\sqrt{3}}\times 10^{-6}, \frac{1}{\sqrt{3}}\times 10^{-6}, 0)^\top\) and \(\lambda_{\min}(S^0_{I_{\neq 0}}(\nabla^2 f(x^0))_{I_{\neq 0}}S^0_{I_{\neq 0}}) = -4\). These results imply that \(x^0\) is a strong* (weak) \(\veg\)-1o point with \(\veg \leq 10^{-5}\) but is not a strong* (weak) \((\veg, \veh)\)-2o point for any \(\veh\in(0, 1)\).
It can be observed from Figure~\ref{fig:toyexamplen3} that the trial points generated by Algorithms~\ref{alg:hpgnc} and~\ref{alg:pncg} differ from one another. The variations in the function values indicate that the Newton-CG step accelerates the decrease in the objective function. For Problem~\eqref{eq:toyn3}, the output \((3, -3, 0)^\top\) is an optimal solution.

\begin{figure}[h!]
\begin{minipage}[t]{1\linewidth}
\centering
\includegraphics[width = 0.325\textwidth]{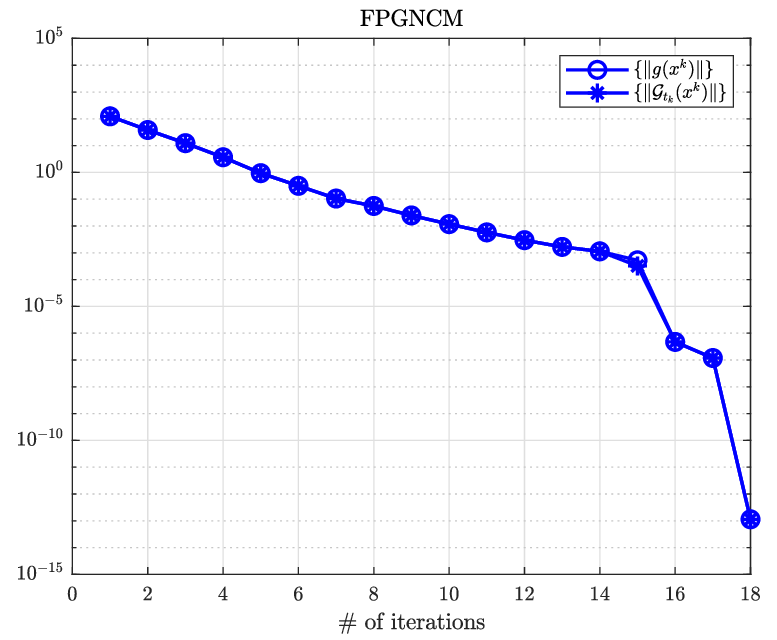}
\includegraphics[width = 0.325\textwidth]{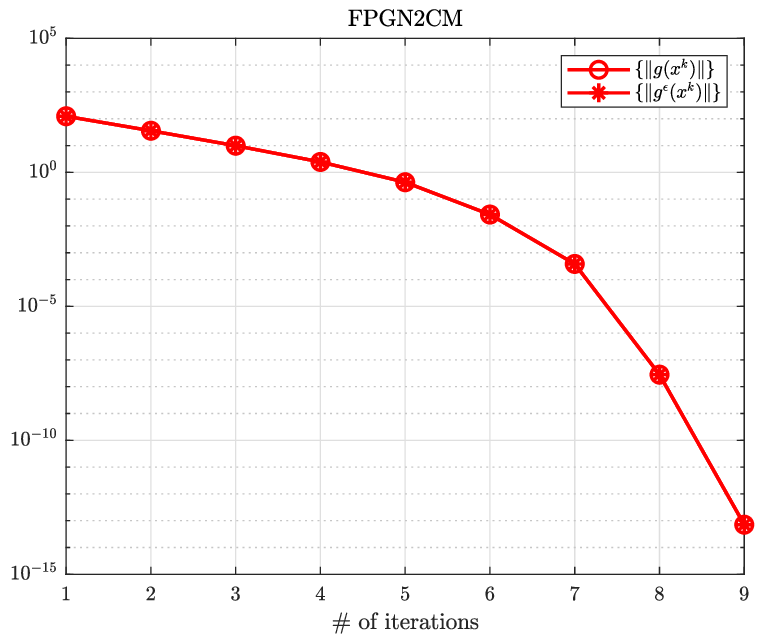}
\includegraphics[width = 0.335\textwidth]{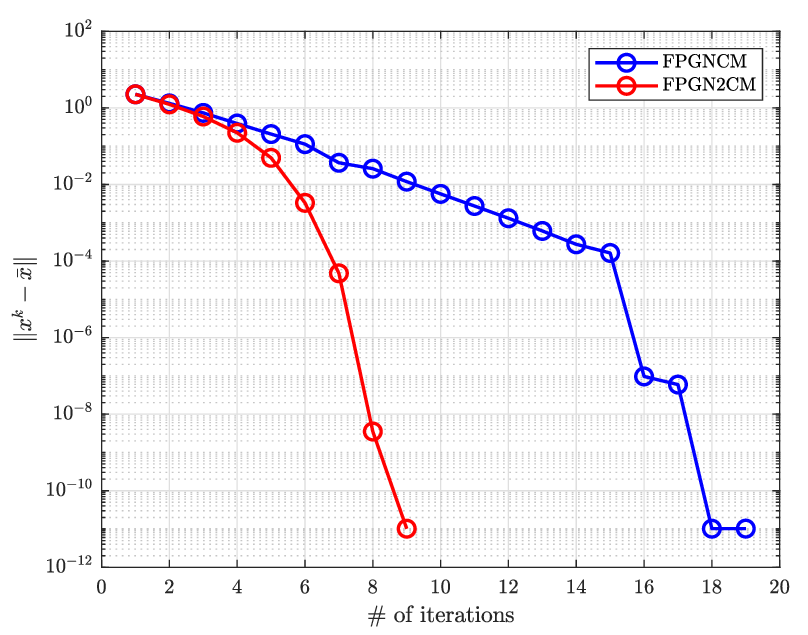}
\end{minipage}
\caption{Left: the performance of \(\|g(x^k)\|\) and \(\|\mathcal{G}_{t_{k+1}}(x^k)\|\) across iterations generated by Algorithm~\ref{alg:hpgnc}; middle: the performance of \(\|g(x^k)\|\) and \(\|g^{\varepsilon}(x^k)\|\) across iterations generated by Algorithm~\ref{alg:pncg}; right: the performance of \(\|x^k - \bar{x}\|\) plotted against iteration.}\label{fig:toyngdx}
\end{figure}

Figure~\ref{fig:toyngdx} displays the performance of the first phase of Algorithms~\ref{alg:hpgnc} (FPGNCM) and~\ref{alg:pncg} (FPGN2CM) in terms of \(\|\mathcal{G}_{t_{k+1}}(x^k)\|\), \(\|g^{\varepsilon}(x^k)\|\), and \(\|x^k - \bar{x}\|\) across iterations, where \(x^0 = (-1.245334, -1.054100,-0.318778)^\top\) is generated randomly and \(\bar{x} = (1, -1, 0)^\top\) denotes the output of the first phase of Algorithm~\ref{alg:pncg} after \(1000\) iterations. It can be verified that \(\bar{x}\) is a strong* (weak) \((\veg, \veh)\)-2o point for any \(\veh\in(0, 1)\) but not an optimal solution of Problem~\eqref{eq:toyn3}. 
The left panel of Figure~\ref{fig:toyngdx} displays the decay of $\|\mathcal{G}_{t_{k+1}}(x^k)\|$ and $\|g(x^k)\|$ generated by FPGNCM, and the middle panel displays the decay of $\|g(x^k)\|$ and $\|g^\varepsilon(x^k)\|$ generated by FPGN2CM.  In this numerical test, $\|\mathcal{G}_{t_{k+1}}(x^k)\|$, $\|g^\epsilon(x^k)\|$, and $\|x^k-\bar{x}\|$ exhibit empirically accelerated decreases.

\subsection{\(\ell_1\)-regularized Student's \(t\)-regression}\label{subsec:st}

The \(\ell_1\)-regularized Student's \(t\)-regression~\cite{AFHL10} problem takes the form of 
 \begin{equation}\label{eq:st}
 \min_x\sum_{i=1}^m\log(1 + (Ax-b)_i^2/\nu) + \lambda \|x\|_1,
 \end{equation}
 where \(\lambda > 0\) is the regularization parameter and \(\nu\) can be interpreted as a tuning parameter: for low values of \(\nu\), one expects a high degree of non-normality, but as \(\nu\) increases, the distribution behaves more like a Gaussian distribution. This interpretation is highlighted in~\cite{LLT89,AFHL10}. Problem~\eqref{eq:st} was used to test the numerical performance of proximal Newton-type method~\cite{LPWY24,Z25}.  Defining \(\phi(u) = \sum_{i=1}^m\log(1 + u_i^2/\nu)\), we have \(f(x) = \phi(Ax - b)\). To verify Assumption~\ref{assum:problem}, consider the scalar function \(\psi(t) = \log(1 + t^2/\nu)\). Its derivatives are \(\psi''(t) = \frac{2(\nu - t^2)}{(\nu + t^2)^2}\) and \(\psi'''(t) = \frac{4t(t^2 - 3\nu)}{(\nu + t^2)^3}\). It can be shown that \(\sup_{t\in\mathbb{R}}\vert \psi''(t)\vert = \frac{2}{\nu}\) and \(\sup_{t\in\mathbb{R}}\vert \psi'''(t)\vert = \frac{1}{(6 - 4\sqrt{2})\nu^{3/2}}\). Consequently, the Lipschitz constants of \(\nabla f(x)\) and \(\nabla^2 f(x)\) are bounded by \(L_g = \frac{2}{\nu}\|A\|^2\) and \(L_H = \frac{1}{(6 - 4\sqrt{2})\nu^{3/2}}\|A\|^3\), respectively, ensuring that Assumption~\ref{assum:problem} is satisfied for Problem \eqref{eq:st}. Furthermore, the Hessian is given by \(\nabla^2 f(x) = A^\top \nabla^2\phi(Ax - b)A\), where \(\nabla^2\phi\) is a diagonal matrix with entries \((\nabla^2\phi(u))_{ii} = \frac{2(\nu - u_i^2)}{(\nu + u_i^2)^2}\). Notably, a smaller value of \(\nu\) increases the likelihood that \((\nabla^2\phi(u))_{ii} < 0\), introducing more pronounced nonconvexity into the problem. We set \(\nu = 0.001\) and \(\lambda = 0.1\|\nabla f(0)\|_{\infty}\) in the following tests. 

 The test examples are randomly generated as follows. The true sparse signal \(x^{\rm true}\in\mathbb{R}^n\), with \(n \in \{2^8, 2^9, 2^{10}\}\), is \(k\)-sparse, where \(k \in \{\lfloor n/20 \rfloor, \lfloor n/40 \rfloor, \lfloor n/60 \rfloor\}\). The support of \(x^{\rm true}\) is chosen uniformly at random from \(\{1, \cdots, n\}\), and the magnitude of each nonzero entry is set to \(x^{\rm true}_i = \eta_1(i)10^{d\eta_2(i)/20}\). Here, \(\eta_1(i) \in\{-1, +1\}\) is a symmetric random sign, and \(\eta_2(i)\) is uniformly distributed in \([0, 1]\). The signal has a dynamic range of \(d\) dB. The difficulty of solving Problem~\eqref{eq:st} increases with \(d\). The matrix \(A\in\mathbb{R}^{m\times n}\) is a partial discrete cosine transform (DCT) matrix with 
\(m = n/8\), formed by randomly selecting a row index set \(J \subset \{1, \cdots, n\}\) with \(\vert J\vert = m\) and retaining the corresponding rows of the \(n\times n\) DCT matrix. The measurement \(b\) is obtained by adding Student's t-noise with degrees of freedom \(5\), rescaled by \(0.1\), to \(Ax^{\rm true}\). 

\subsubsection{Comparing HPGNCM and PGN2CM}

We compare the empirical iteration counts and operation counts of Algorithms~\ref{alg:hpgnc} and~\ref{alg:pncg} in light of the global complexity bounds established in Theorems~\ref{th:soicpgnc} and~\ref{th:soic} and Corollaries~\ref{cor:hpgnc} and~\ref{coro:operpgn2cm}.

We evaluate the performance of Algorithms~\ref{alg:hpgnc} and~\ref{alg:pncg} on instances with \(d \in \{20, 40\}\), using \(\veg = 10^{-6}\) and \(\veh = \veg^{1/2}\). In this test, we set \(\beta = 3\), \(\bar{\eta} = 0.475\), and \(\theta = 0.3\) in Algorithm~\ref{alg:hpgnc}. We terminate Algorithm~\ref{alg:hpgnc} if \(\|\mathcal{G}_{t_{k+1}}(x^k)\| \leq \veg\) and \(\lambda_{\min}(S^k_{\neq 0}(\nabla^2 f(x^k))_{\neq 0}S^k_{\neq 0}) \geq -\veh\), or if the number of iterations reaches to \(2{,}000{,}000\). 
For Algorithm~\ref{alg:pncg}, we set \(\beta = 2.75\), \(\theta = 0.75\) in~\eqref{eq:lsreq}, and \(\theta = 0.3\) in~\eqref{eq:lsreq2}. We terminate Algorithm~\ref{alg:pncg} if \(\|g^{\varepsilon}(x^k)\| \leq \veg\) and \(\lambda_{\min}(S^k_{\neq 0}(\nabla^2 f(x^k))_{\neq 0}S^k_{\neq 0}) \geq -\veh\), or if the number of iterations reaches \(2,000,000\). Table~\ref{table:proxgncvsproxncg}, which is presented in two panels, reports the results averaged over \(20\) independent trials, including the number of iterations (Iter), the objective function value (Fval), and running time in seconds (time(s)), stationarity residuals, minimum eigenvalues, and the total number of operations (\#Ops). Specifically, for the final output, we report the residual norms (\(\|\mathcal{G}_t\|\), \(\|g(x)\|\), \(\|g^{\varepsilon}(x)\|\)), as well as the minimum eigenvalues of the sub-Hessian (\(\lambda_{\min}^{H_{\neq 0}}\)) and the scaled sub-Hessian (\(\lambda_{\min}^{S_{\neq 0}H_{\neq 0}S_{\neq 0}}\)). The total operation count (\#Ops) summarizes all gradient evaluations, matrix-vector products, and MEO calls during the entire process. From Table~\ref{table:proxgncvsproxncg}, we can see that (i) Across all test instances, PGN2CM consistently requires significantly fewer iterations than HPGNCM. Furthermore, in the vast majority of cases, PGN2CM achieves a substantially lower total running time. This computational advantage in terms of runtime becomes increasingly prominent as the problem scale $n$, the dynamic range $d$, and the signal density $k$ increase. (ii) In most cases, HPGNCM failed to terminate within the maximum limit of $2{,}000{,}000$ iterations in several trials, as evidenced by $\|\mathcal{G}_t\| > 10^{-6}$. In the reported averages, PGN2CM reaches residual levels close to or below the prescribed tolerance in all test groups, whereas HPGNCM exceeds the iteration limit in several groups. (iii) We observe that PGN2CM requires more \#Ops than HPGNCM for most instances when $d = 20$, but fewer \#Ops when $d = 40$. This trend is consistent across all three sparsity levels and aligns with our theoretical results in Corollaries~\ref{cor:hpgnc} and~\ref{coro:operpgn2cm}, which establish the operation complexities of $\mathcal{O}(\varepsilon^{-9/4})$ for HPGNCM and $\tilde{\mathcal{O}}(\varepsilon^{-7/4})$ for PGN2CM. Although  PGN2CM involves more intensive computations per Newton-CG step, such as Hessian-vector products, its significantly lower iteration complexity dominates the total cost as the problem difficulty increases, eventually leading to a lower overall operation count in more complex scenarios. iv) The overall relative performance of the two algorithms is qualitatively robust to the choice of sparsity level. An exception is observed when $d = 20$ and $k = \lfloor n/60 \rfloor$, where HPGNCM achieves a shorter running time than PGN2CM despite requiring more iterations. This suggests that, in this simpler and sparser setting, the frequent invocation of the Newton-CG step in PGN2CM may actually hinder computational efficiency. Although Newton-CG steps target a theoretical $\mathcal{O}(\varepsilon^{3/2})$ objective decrease, they often require small step sizes $t \in (0, \min\{t_k^+, t_k^-\}]$, as specified in Lemma~\ref{lem:dphi2}, to ensure stability. In simpler landscapes, this restricted progress, combined with the high computational overhead, makes them less cost-effective than standard proximal gradient steps. Moreover, we also observe that HPGNCM can obtain strong* \((\veg, \veh)\)-2o points in some instances with \(d = 20\), whereas all outputs of  PGN2CM are weak \((\veg, \veh)\)-2o points. 

\begin{sidewaystable*}[p] 
{\small
\caption{Numerical comparisons of Algorithms~\ref{alg:hpgnc} (HPGNCM) and~\ref{alg:pncg} (PGN2CM) on regularized Student's \(t\)-regression.}\label{table:proxgncvsproxncg} {Panel A: Iteration counts, running times, and residual norms}\\
\setlength{\tabcolsep}{1pt}
\begin{tabular*}{\textwidth}{@{\extracolsep{\fill}}ccc|ccc|ccc|ccc|ccc|ccc}\hline\hline
\multirow{2}{*}{\(d\)} &   \multirow{2}{*}{\(n\)}     &\multirow{2}{*}{Algs.} & \multicolumn{3}{c|}{Iter} & \multicolumn{3}{c|}{time(s)}  & \multicolumn{3}{c|}{\(\|\mathcal{G}_t\|\)}          &\multicolumn{3}{c|}{\(\|g\|\)} &\multicolumn{3}{c}{\(\|g^{\varepsilon}\|\)}                                                                                                                                               \\ \cline{4-18}
& & & \(\lfloor n/20 \rfloor\) & \(\lfloor n/40 \rfloor\) & \(\lfloor n/60 \rfloor\) & \(\lfloor n/20 \rfloor\) & \(\lfloor n/40 \rfloor\) & \(\lfloor n/60 \rfloor\) &\(\lfloor n/20 \rfloor\) & \(\lfloor n/40 \rfloor\) & \(\lfloor n/60 \rfloor\) &\(\lfloor n/20 \rfloor\) & \(\lfloor n/40 \rfloor\) & \(\lfloor n/60 \rfloor\) &\(\lfloor n/20 \rfloor\) & \(\lfloor n/40 \rfloor\) & \(\lfloor n/60 \rfloor\) \\\hline\hline  
\multirow{6}{*}{\(20\)} & \multirow{2}{*}{\(2^8\)} & HPGNCM & 485201 & 160136 & 199755 & 7.40 & {\color{blue}2.59} & {\color{blue}3.04}& {\color{red}2.06E-6} & 9.89E-07 & 9.93E-7& 1.76E-6 & 9.90E-4 & 1.04E-6& $-$ & $-$ & $-$ \\
& & PGN2CM & {\color{blue}53872} & {\color{blue}68069} & {\color{blue}87056} & {\color{blue}4.79} &6.13 & 7.19& $-$ & $-$ &$-$ & 4.28E-3& 4.10E-3 & 6.29E-3& 7.98E-7 & 7.54E-7 & 7.04E-7 \\ \cline{2-18}
& \multirow{2}{*}{\(2^9\)} & HPGNCM & 712823 & 388263 & 190832 & 30.27 & 16.70 & {\color{blue}6.58} & {\color{red}1.31E-6} & {\color{red}1.35E-6} & 9.94E-7& 1.41E-6 & 1.60E-6 & 9.96E-7& $-$ & $-$ & $-$ \\
& & PGN2CM & {\color{blue}31993} & {\color{blue}38320} & {\color{blue}57398} & {\color{blue}7.72} & {\color{blue}9.09} & 12.31& $-$ & $-$ & $-$ & 1.04E-2& 1.64E-3 & 2.69E-3& 7.32E-7 & 7.50E-7 & 7.30E-7 \\ 
\cline{2-18}
& \multirow{2}{*}{\(2^{10}\)}& HPGNCM & 973504 & 463595 & 311028 & 62.17 & 32.39 & {\color{blue}19.95} & {\color{red}4.85E-6} &  {\color{red}4.24E-6} & {\color{red}1.18E-6} & 4.67E-6 & 2.10E-6 & 1.19E-6& $-$ & $-$ & $-$ \\
& & PGN2CM & {\color{blue}20248} & {\color{blue}42581} & {\color{blue}46285} & {\color{blue}15.78} & {\color{blue}29.92} & 31.53& $-$ & $-$ & $-$& 1.03E-3& 1.83E-3 & 1.77E-3& 7.73E-7 & 7.52E-7 & 7.20E-7 \\ \hline
\multirow{6}{*}{\(40\)} & \multirow{2}{*}{\(2^8\)} & HPGNCM & 1777548 & 1020257 & 1094973 & 26.12 & 15.34 & 16.13& {\color{red}1.02E-3} & {\color{red}1.71E-4} & {\color{red}1.18E-4} & 8.94E-4 & 2.37E-4 & 2.10E-4& $-$ & $-$ & $-$ \\
& & PGN2CM & {\color{blue}31710} &{\color{blue} 44526} & {\color{blue}30848} & {\color{blue}2.75} & {\color{blue}3.75} & {\color{blue}2.58}& $-$ & $-$ & $-$& 2.06E-3& 2.41E-3 & 1.35E-3& 7.44E-7 & 7.09E-7 & 7.21E-7 \\
\cline{2-18}
&  \multirow{2}{*}{\(2^9\)}& HPGNCM & 1970732 & 1313245 & 851807 & 75.17 & 50.04 & 28.88& {\color{red}1.51E-3} & {\color{red}2.00E-5} & {\color{red}1.88E-6}& 1.22E-3 & 2.10E-5 & 1.92E-6& $-$ & $-$ & $-$ \\
& & PGN2CM & {\color{blue}17370} & {\color{blue}27653} & {\color{blue}21814} & {\color{blue}3.07} & {\color{blue}6.62} & {\color{blue}4.61} & $-$ & $-$ & $-$& 8.79E-4& 1.28E-3 & 6.51E-4& 7.01E-7 & 7.53E-7 & 7.54E-7 \\
\cline{2-18}
&  \multirow{2}{*}{\(2^{10}\)}
& HPGNCM & 1973182 & 1467053 & 888424 & 117.45 & 90.43 & 56.33& {\color{red}1.81E-3} & {\color{red}1.95E-5} & {\color{red}1.51E-6} & 1.54E-3 & 1.82E-5 & 1.72E-6& $-$ & $-$ & $-$ \\
& & PGN2CM & {\color{blue}24919} & {\color{blue}27203} & {\color{blue}26263} & {\color{blue}15.04} & {\color{blue}16.38} & {\color{blue}18.84} & $-$ & $-$ & $-$& 5.60E-4& 7.07E-4 & 8.00E-4& 6.97E-7 & 7.21E-7 & 7.28E-7 \\\hline\hline 
\end{tabular*}

\vspace{4mm}
Panel B: Objective values, minimum eigenvalues, and operation counts \\
\setlength{\tabcolsep}{1pt}
\begin{tabular*}{\textwidth}{@{\extracolsep{\fill}}ccc|ccc|ccc|ccc|ccc}\hline\hline
\multirow{2}{*}{\(d\)} &   \multirow{2}{*}{\(n\)}     &\multirow{2}{*}{Algs.} & \multicolumn{3}{c|}{Fval} & \multicolumn{3}{c|}{\(\lambda_{\min}^{H_{\neq 0}}\)}  & \multicolumn{3}{c|}{\(\lambda_{\min}^{S_{\neq 0}H_{\neq 0}S_{\neq 0}}\)}          &\multicolumn{3}{c}{\#Ops}                                                                                                                                               \\ \cline{4-15}
& & & \(\lfloor n/20 \rfloor\) & \(\lfloor n/40 \rfloor\) & \(\lfloor n/60 \rfloor\) &\(\lfloor n/20 \rfloor\) & \(\lfloor n/40 \rfloor\) & \(\lfloor n/60 \rfloor\) &\(\lfloor n/20 \rfloor\) & \(\lfloor n/40 \rfloor\) & \(\lfloor n/60 \rfloor\) &\(\lfloor n/20 \rfloor\) & \(\lfloor n/40 \rfloor\) & \(\lfloor n/60 \rfloor\)  \\\hline\hline 
 \multirow{6}{*}{\(20\)} & \multirow{2}{*}{\(2^8\)} & HPGNCM & 3.21 & 2.39 & 2.09& 4.71E{-1} & 5.88E{-1} & 4.89E{-1}& 4.71E{-1} & 4.89E{-1} & 4.09E{-1}& {\color{blue}4.85E+{5}} & {\color{blue}1.60E+{5}} & {\color{blue}2.00E+{5}} \\
& & PGN2CM & 3.21 & 2.39 & 2.09& 4.75E{-2} & -6.79E{-15} & 4.92E{-3}& 4.75E{-2} & 5.60E{-7} & 6.69E{-7}& 1.04E+{6} & 1.21E+{6} & 1.48E+{6} \\ \cline{2-15}
& \multirow{2}{*}{\(2^9\)}& HPGNCM & 7.57 & 5.81 & 4.59& 1.10E{-1} & 1.26E{-1} & 1.33E{-1}& 1.06E{-1} & 1.13E{-1} & 1.23E{-1}& {\color{blue}7.13E+{5}} & {\color{blue}3.88E+{5}} & {\color{blue}1.91E+{5}} \\
& & PGN2CM & 7.57 & 5.81 & 4.59& -2.48E{-15} & -2.51E{-15} & -4.53E{-16}& 5.36E{-20} & -4.51E{-21} & 1.33E{-19}& 9.05E+{6} & 1.10E+{6} & 1.33E+{6} \\
\cline{2-15}
& \multirow{2}{*}{\(2^{10}\)}& HPGNCM & 17.66 & 12.87 & 10.01& 2.66E{-2} & 2.71E{-2} & 3.34E{-2}& 2.63E{-2} & 2.33E{-2} & 2.94E{-2}& 9.74E+{5} & {\color{blue}4.64E+{5}} & {\color{blue}3.11E+{5}} \\
& & PGN2CM & 17.66 & 12.87 & 10.01& 1.14E{-15} & -2.20E{-15} & -1.36E{-15}& -2.35E{-19}& 4.87E{-20} & 2.61E{-20}& {\color{blue}7.21E+{5}} & 1.19E+{6} & 1.38E+{6} \\ \hline
\multirow{6}{*}{40}
& \multirow{2}{*}{\(2^8\)}& HPGNCM & 5.67 & 3.16 & 2.34& 1.37E{-1} & 2.18E{-1} & 3.33E{-1}& 1.37E{-1} & 3.67E{-1} & 3.21E{-1}& 1.78E+{6} & 1.02E+{6} & 1.09E+{6} \\
& & PGN2CM & {\color{blue}5.66} & {\color{blue}3.15} & 2.34& -3.08E{-15} & 6.05E{-16} & -2.22E{-15}& -1.25E{-20} & 3.44E{-20} & -1.18E{-20}& {\color{blue}6.10E+{5}} & {\color{blue}7.58E+{5}} & {\color{blue}5.62E+{5}} \\
\cline{2-15}
& \multirow{2}{*}{\(2^9\)}& HPGNCM & 15.45 & 8.80 & 5.77& 1.37E{-2} & 7.71E{-2} & 8.98E{-2}& 1.37E{-2}& 2.53E{-2} & 8.06E{-2}& 1.97E+6 & 1.31E+{6} & 8.52E+{5} \\
& & PGN2CM & {\color{blue}15.43} & 8.80 & 5.77& -4.56E{-15} & 4.21E{-4} & -2.02E{-16}& -1.08E{-18} & 4.21E{-4} & -2.97E{-20}& {\color{blue}3.43E+{5}} & {\color{blue}7.49E+{5}} & {\color{blue}5.53E+{5}} \\
\cline{2-15}
& \multirow{2}{*}{\(2^{10}\)}& HPGNCM & 33.79 & 20.21 & 17.67& 2.23E{-3} & 1.58E{-2} & 2.72E{-2}& 2.23E{-3} & 1.58E{-2}& 2.12E{-2}& 1.97E+{6} & 1.47E+{6} & 8.88E+{5} \\
& & PGN2CM & {\color{blue}33.77} & 20.21 & 17.67& -4.22E{-15} & -1.06E{-15} & -1.72E{-15}& -8.78E{-20} & 1.05E{-19} & -2.01E{-19}& {\color{blue}6.44E+{5}} & {\color{blue}7.07E+{5}} & {\color{blue}8.57E+{5}} \\ \hline \hline
\end{tabular*}}
\end{sidewaystable*}

\begin{figure}[h!]
\begin{minipage}[t]{1\linewidth}
\centering
\includegraphics[width = 0.325\textwidth]{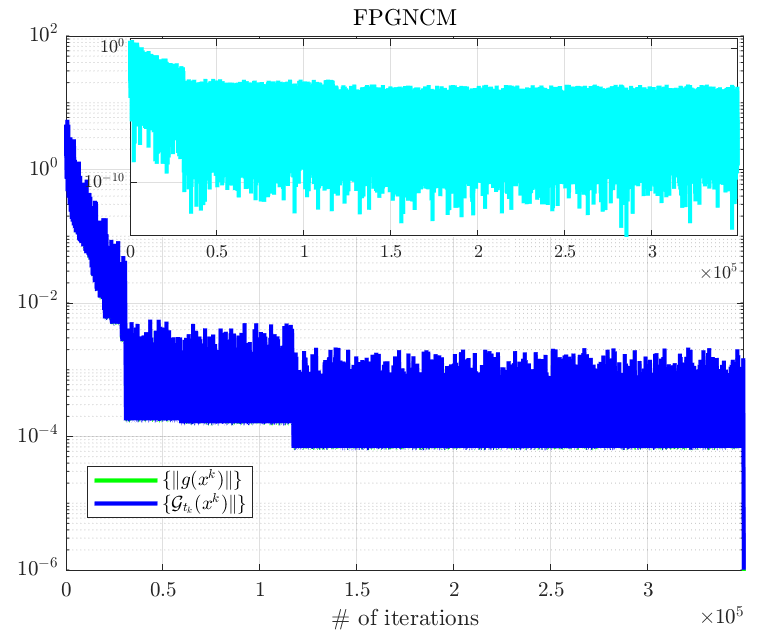}
\includegraphics[width = 0.325\textwidth]{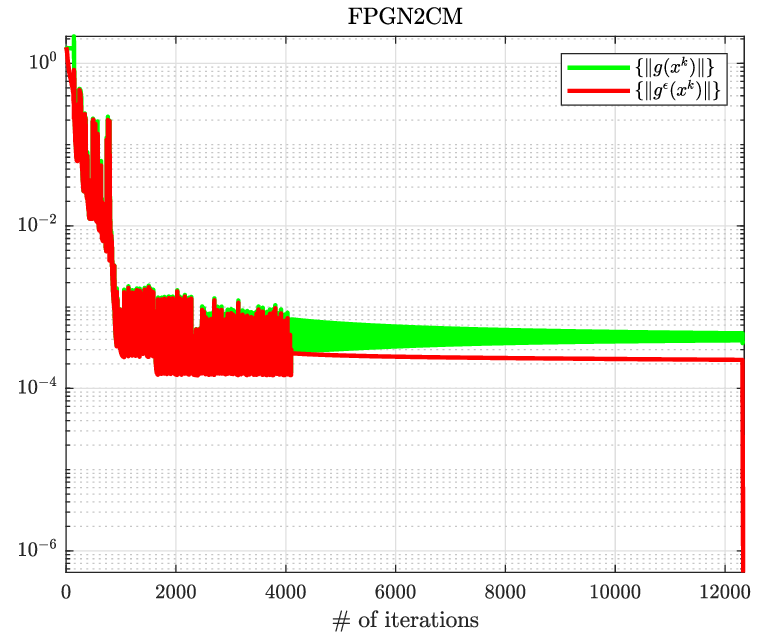}
\includegraphics[width = 0.335\textwidth]{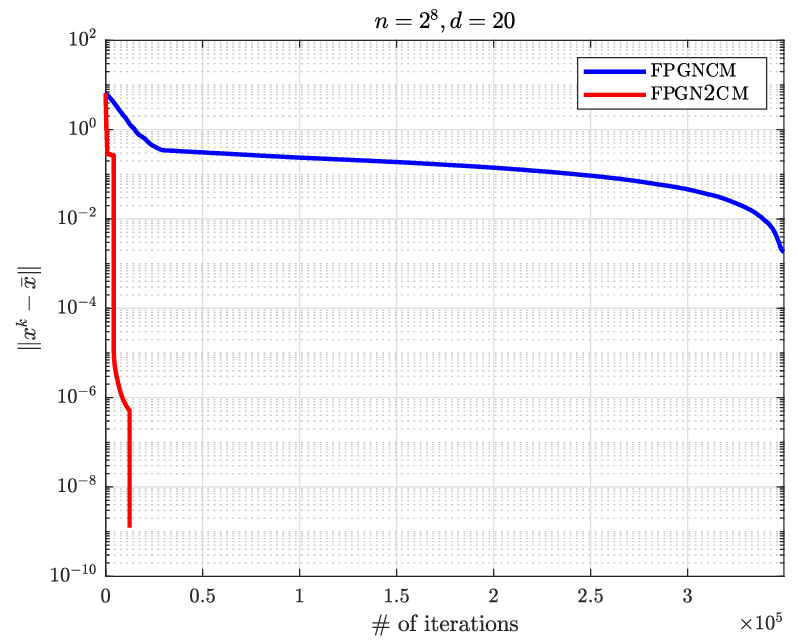}
\end{minipage}
\caption{Left: the performance of \(\|g(x)\|\) and \(\|\mathcal{G}_{t_{k+1}}(x^k)\|\) across iterations generated by Algorithm~\ref{alg:hpgnc}; middle: the performance of \(\|g(x)\|\) and \(\|g^{\varepsilon}(x)\|\) across iterations generated by Algorithm~\ref{alg:pncg}; right: the performance of \(\|x^k - \bar{x}\|\) plotted against iteration.}\label{fig:stngdx}
\end{figure}

Figure~\ref{fig:stngdx} displays the performance of the first phase of Algorithms~\ref{alg:hpgnc} (FPGNCM) and~\ref{alg:pncg} (FPGN2CM) in terms of \(\|\mathcal{G}_{t_{k+1}}(x^k)\|\), \(\|g^{\varepsilon}(x^k)\|\), and \(\|x^k - \bar{x}\|\) against iteration count for an instance with \(n = 2^8\) and \(d = 20\), where \(\bar{x}\) denotes the output of FPGN2CM after \(500\, 000\) iterations. The left and middle subplots show that the measured residuals are consistent with the valid inequalities \(\|\mathcal{G}_{t_{k+1}}(x^k)\| \leq \|g(x^k)\|\) and \(\|g^{\ve}(x^k)\| \leq \|g(x^k)\|\). The right subplot illustrates empirically accelerated decay of \(\{\|x^k - \bar{x}\|\}\).

\subsubsection{Comparing PGN2CM with baseline methods}
We compare Algorithm~\ref{alg:pncg} with  the forward-backward splitting (FBS) method~\cite{LM79,CR97}, ZeroFPR(L-BFGS)~\cite{TSP18} downloaded from the website~\url{https://github.com/kul-optec/ForBES}, and the Proximal Newton method~\cite{Z25} to assess whether the second-order information exploited by PGN2CM translates into a practical computational advantage over these methods. For the baseline methods, FBS and ZeroFPR(L-BFGS) are run with the default parameter settings of their respective implementations, whereas the parameters of the Proximal Newton method are tuned to obtain satisfactory numerical results. To ensure a fair comparison, all baseline methods use the stopping criterion \(\|g^{\varepsilon}(x^k)\| \leq \veg\) with \(\veg = 10^{-6}\) and a maximum of $1{,}000{,}000$ iterations. Apart from the stated stopping criterion, maximum iteration limit, and parameter tuning of the Proximal Newton method, no additional modifications were made to the baseline implementations. PGN2CM is terminated according to its own weak second‑order stopping criterion with \(\veg = 10^{-6}\) and the same maximum iteration limit.
For Algorithm~\ref{alg:pncg}, we set \(\beta = 2\), \(\theta = 0.25\) in~\eqref{eq:lsreq}, and \(\theta = 0.25\) in~\eqref{eq:lsreq2}. 
Table~\ref{table:student}, which is presented in two panels, reports the average number of iterations (Iter), the objective function value (Fval), the running time in second (time(s)), the norm of \(g^{\varepsilon}(x)\) at the output (\(\|g^{\varepsilon}\|\)), and the minimum eigenvalues of \((\nabla^2 f(x))_{I_{\neq 0}}\) and \(S_{I_{\neq 0}}(\nabla^2 f(x))_{I_{\neq 0}}S_{I_{\neq 0}}\) at the output (\(\lambda_{\min}^{H_{\neq 0}}\) and \(\lambda_{\min}^{S_{\neq 0}H_{\neq 0}S_{\neq 0}}\)) over \(20\) independent trials. We exclude HPGNCM from the comparison. As already observed in Table~\ref{table:proxgncvsproxncg}, HPGNCM struggles to satisfy the stopping criterion within the maximum iteration limit even for simpler instances (e.g., $d=40$). For more challenging cases with $d=60$ and $d=80$, HPGNCM consistently reaches the iteration limit without satisfying the stopping criterion and thus provides no meaningful numerical benchmark. From Table~\ref{table:student}, we make the following observations: (i) In terms of computational time, the overall relative performance of PGN2CM is more robust to the choice of sparsity level than the comparison methods; (ii) For all test instances, the iterate returned by PGN2CM satisfies the weak-2o optimality condition; (iii) For the challenging test problems with \(d=80\), PGN2CM yields objective function values no higher than those produced by the competing methods, and outperforms them strictly on several instances; (iv) Across all numerical tests, although PGN2CM does not always yield the smallest iteration count, it consumes less computational time than all comparison methods. Moreover, its runtime advantage becomes increasingly pronounced as the test problems grow more difficult; (v) For all test instances, the iterate generated by PGN2CM satisfies the weak \(\varepsilon_g\)-1o optimality condition (i.e., \(\|g^{\varepsilon}\| \leq \varepsilon_g\)). For \(d = 80\), L-BFGS and FBS fail to satisfy the weak \(\varepsilon_g\)-1o condition in most trials. 
The degraded performance of the baseline methods on the \(d=80\) instances is consistent with the increased numerical difficulty induced by the larger dynamic range. 

\begingroup  
\setlength{\tabcolsep}{1pt}  
\begin{table}[!h]
\centering
\caption{Numerical comparisons on regularized Student's \(t\)-regression. \\
Panel A: Iteration counts, running times, and objective values}\label{table:student}
\vspace{-3mm}
\begin{minipage}{\linewidth} 
\renewcommand{\arraystretch}{1.05}
\begin{tabular*}{\textwidth}{@{\extracolsep{\fill}}ccc|ccc|ccc|ccc}\hline\hline
\multirow{2}{*}{\(d\)} &   \multirow{2}{*}{\(n\)}     &\multirow{2}{*}{Algs.} & \multicolumn{3}{c|}{Iter} & \multicolumn{3}{c|}{time(s)}  & \multicolumn{3}{c}{Fval} \\ \cline{4-12}
& & & \(\lfloor n/20 \rfloor\) & \(\lfloor n/40 \rfloor\) & \(\lfloor n/60 \rfloor\) & \(\lfloor n/20 \rfloor\) & \(\lfloor n/40 \rfloor\) & \(\lfloor n/60 \rfloor\) &\(\lfloor n/20 \rfloor\) & \(\lfloor n/40 \rfloor\) & \(\lfloor n/60 \rfloor\) \\ \hline\hline
\multirow{8}{*}{\(60\)}& \multirow{4}{*}{\(2^8\)}& L-BFGS & 212534 & 56124 & 47419 & 30.96 & 8.31 & 7.20 & 6.00 & 3.63 & 2.42 \\
& & FBS & 505073 & 201772 & 243783 & 43.96 & 17.68 & 21.80 & 6.00 & 3.63 & 2.42 \\
& & ProxNewton & {\color{blue}21981} & {\color{blue}25407} & {\color{blue}42902} & 3.67 & 3.52 & 5.75 & 6.00 & 3.63 & 2.42 \\
& & PGN2CM & 64612 & 45340 & 56254 & {\color{blue}2.78} & {\color{blue}2.23} & {\color{blue}2.57} & 6.00 & 3.63 & 2.42 \\\cline{2-12}
& \multirow{4}{*}{\(2^9\)}& L-BFGS & 195527 & {\color{blue}63593} & {\color{blue}24867} & 40.80 & 13.88 & 5.26 & 19.41 & 10.62 & 8.49 \\
& & FBS & 483427 & 303426 & 187710 & 62.82 & 41.26 & 25.16 & 19.41 & 10.62 & 8.49 \\
& & ProxNewton & 682472 & 427095 & 400659 & 206.03 & 91.49 & 88.33 & 19.41 & 10.62 & 8.49 \\
& & PGN2CM & {\color{blue}118146} & 89366 & 48304 & {\color{blue}12.97} & {\color{blue}7.78} & {\color{blue}5.07} & 19.41 & 10.62 & 8.49 \\ \cline{1-12}
\multirow{8}{*}{\(80\)}& \multirow{4}{*}{\(2^8\)}& L-BFGS & 838896 & 333506 & 344848 & 123.11 & 48.22 & 51.67 & 10.77 & 6.08 & 4.43 \\
& & FBS & 964529 & 593587 & 748829 & 84.72 & 50.08 & 66.33 & 10.73 & 6.08 & 4.47 \\
& & ProxNewton & 460348 & {\color{blue}111463} & {\color{blue}98109} & 88.21 & 19.19 & 15.26 & {\color{blue}10.66} & {\color{blue}6.04} & {\color{blue}4.42} \\
& & PGN2CM & {\color{blue}334820} & 150711 & 172604 & {\color{blue}28.78} & {\color{blue}7.74} & {\color{blue}9.14} & {\color{blue}10.66} & {\color{blue}6.04} & {\color{blue}4.42} \\ \cline{2-12}
& \multirow{4}{*}{\(2^9\)}& L-BFGS & 971536 & 415753 & {\color{blue}158586} & 200.41 & 91.59 & 33.43 & 12.01 & 13.93 & {\color{blue}7.93} \\
& & FBS & 1000000 & 910086 & 741323 & 129.73 & 123.89 & 96.87 & 12.07 & 13.98 & 7.98 \\
& & ProxNewton & 1000000 & 965545 & 1000000 & 537.05 & 396.93 & 361.05 & 11.95 & 13.91 & {\color{blue}7.93} \\
& & PGN2CM & {\color{blue}481862} & {\color{blue}307613} & 218395 & {\color{blue}96.89} & {\color{blue}31.03} & {\color{blue}20.10} & {\color{blue}11.82} & {\color{blue}13.90} & {\color{blue}7.93} \\ \hline \hline
\end{tabular*}

\vspace{2mm}

Panel B: stationarity residual and minimum eigenvalues
\setlength{\tabcolsep}{0.5pt}
\begin{tabular*}{\textwidth}{@{\extracolsep{\fill}}ccc|ccc|ccc|ccc}\hline\hline
\multirow{2}{*}{\(d\)} &   \multirow{2}{*}{\(n\)}     &\multirow{2}{*}{Algs.} &  \multicolumn{3}{c|}{\(\|g^{\varepsilon}\|\)}    &  \multicolumn{3}{c|}{\(\lambda_{\min}^{H_{\neq 0}}\)}  & \multicolumn{3}{c}{\(\lambda_{\min}^{S_{\neq 0}H_{\neq 0}S_{\neq 0}}\)} \\ \cline{4-12}
& & & \(\lfloor n/20 \rfloor\) & \(\lfloor n/40 \rfloor\) & \(\lfloor n/60 \rfloor\) &\(\lfloor n/20 \rfloor\) & \(\lfloor n/40 \rfloor\) & \(\lfloor n/60 \rfloor\) & \(\lfloor n/20 \rfloor\) & \(\lfloor n/40 \rfloor\) & \(\lfloor n/60 \rfloor\) \\ \hline\hline
\multirow{8}{*}{\(60\)}
& \multirow{4}{*}{\(2^8\)} & L-BFGS& 8.39E{-7} & 8.45E{-7} & 8.36E{-7}& 3.79E{-1} & 4.05E{-1} & 4.72E{-1}& 3.79E{-1} & 4.05E{-1} & 4.61E{-1} \\
& & FBS& {\color{red}2.77E{-5}} & 9.11E{-7} & 9.35E{-7}& 3.41E{-1} & 3.89E{-1} & 4.60E{-1}& 3.41E{-1} & 3.89E{-1} & 4.49E{-1} \\
& & ProxNewton& 5.15E{-7} & 6.70E{-7} & 5.76E{-7}& 4.04E{-1} & 4.22E{-1} & 4.66E{-1}& 4.04E{-1} & 3.93E{-1} & 4.55E{-1} \\
& & PGN2CM& 7.54E{-7} & 7.16E{-7} & 7.29E{-7}& -2.76E{-15} & -4.83E{-15} & -4.93E{-15}& 1.55E{-20} & -3.67E{-20} & -8.08E{-21} \\ \cline{2-12}
& \multirow{4}{*}{\(2^9\)} & L-BFGS& 8.79E{-7}& 8.71E{-7} & 9.01E{-7}& 8.67E{-2} & 8.94E{-2} & 1.08E{-1}& 8.67E{-2} & 8.94E{-2} & 9.85E{-2} \\
& & FBS& {\color{red}4.32E{-6}} & 9.80E{-7} & 9.79E{-7}& 7.64E{-2} & 8.26E{-2} & 1.06E{-1}& 7.64E{-2} & 8.26E{-2} & 9.71E{-2} \\
& & ProxNewton& {\color{red}2.14E{-5}} & {\color{red}3.83E{-6}} & {\color{red}3.32E{-6}}& 8.67E{-2} & 9.09E{-2} & 1.08E{-1}& 8.67E{-2} & 9.09E{-2} & 9.94E{-2} \\
& & PGN2CM& 7.19E{-7} & 7.85E{-7} & 7.22E{-7}& -4.85E{-15} & -2.35E{-15} & 1.11E{-15}& -2.17E{-19} & 5.63E{-20} & 1.30E{-20} \\\cline{1-12}
\multirow{8}{*}{\(80\)}& \multirow{4}{*}{\(2^8\)}& L-BFGS& {\color{red}6.69E{-4}} & {\color{red}3.28E{-4}} & {\color{red}5.50E{-5}} & 2.09E{-1} & 2.96E{-1} & 3.20E{-1}& 2.09E{-1} & 2.96E{-1} & 3.20E{-1} \\
& & FBS& 2.49E{-3} & 2.17E{-3} & 2.55E{-3}& 8.65E{-2} & 1.98E{-1} & 8.10E{-2}& 8.65E{-2}& 1.98E{-1} & 8.10E{-2} \\
& & ProxNewton& {\color{red}5.98E{-5}} & 5.96E{-7} & 6.24E{-7}& 3.75E{-1} & 4.23E{-1} & 4.11E{-1}& 3.75E{-1} & 4.23E{-1} & 4.00E{-1} \\
& & PGN2CM& 7.28E{-7} & 7.35E{-7} & 7.09E{-7}& -5.79E{-15} & -8.59E{-15} & -9.14E{-15}& 4.29E{-20} & 1.64E{-20} & 9.50E{-21} \\ \cline{2-12}
& \multirow{4}{*}{\(2^9\)}& L-BFGS& {\color{red}2.44E{-3}} & {\color{red}9.57E{-5}} & 8.81E{-7}& 4.14E{-3} & 5.56E{-2} & 7.20E{-2}& 4.14E{-3} & 5.56E{-2} & 6.51E{-2} \\
& & FBS& {\color{red}8.33E{-3}} & {\color{red}6.66E{-3}} & {\color{red}4.72E{-3}} & -3.15E{-14} & 1.82E{-2} & 2.94E{-2}& -3.57E{-14} & 1.82E{-2} & 2.25E{-2} \\
& & ProxNewton& {\color{red}4.52E{-4}} & {\color{red}9.91E{-5}} & {\color{red}1.66E{-5}} & 2.10E{-3} & 7.48E{-2} & 8.81E{-2}& 2.10E{-3} & 7.48E{-2} & 8.12E{-2} \\
& & PGN2CM& 7.14E{-7} & 7.12E{-7} & 7.06E{-7}& -1.15E{-14} & -7.72E{-15} & -6.41E{-15}& -1.28E{-20} & -2.57E{-20} & 1.59E{-21} \\ \hline\hline
 \end{tabular*}
 \end{minipage}
\end{table}
\endgroup 

In summary, the numerical experiments yield three main observations. First, the toy example illustrates that both HPGNCM and PGN2CM are capable of escaping first-order stationary points that violate the second-order condition, providing numerical evidence supporting the motivation to incorporate negative-curvature directions. Second, the comparison between the two proposed methods on Student's $t$-regression problem shows that PGN2CM consistently requires fewer iterations than HPGNCM, and the computational advantage becoming more pronounced as the problem difficulty increases, in agreement with the improved iteration complexity established in Theorem~\ref{th:soic}. Third, PGN2CM exhibits competitive performance against FBS, L-BFGS and the Proximal Newton method, especially on challenging instances where baseline methods fail to converge within the prescribed iteration limit.

\section{Conclusions}\label{sec:con}

In this paper, we present two second-order methods for \(\ell_1\)-norm regularized nonconvex composite problems. We introduce the concepts of strong* and weak \((\veg, \veh)\)-second-order stationary points. We show that the iteration and computational complexity of the hybrid proximal gradient and negative curvature method for identifying a strong* \((\varepsilon, \varepsilon^{\frac{1}{2}})\)-2o point are \(\mathcal{O}(\varepsilon^{-2})\) and \({\mathcal{O}}(\varepsilon^{-\frac{9}{4}})\) with high probability, respectively. In contrast, the iteration and computational complexities of the proximal gradient-Newton-CG method for identifying a weak \((\varepsilon, \varepsilon^{\frac{1}{2}})\)-2o point are \(\mathcal{O}(\varepsilon^{-3/2})\) and \(\widetilde{\mathcal{O}}(\varepsilon^{-\frac{7}{4}})\) with high probability, respectively. The latter result is consistent with the analysis of the Newton-CG method for unconstrained optimization~\cite{RNW20} and the projected Newton-CG method for bound-constrained optimization~\cite{XW23}.

\appendix

\section{Proof of Lemma~\ref{lem:proofggt}}\label{app:proofggt}
\begin{proof}
(a) The Inequality \(\|\mathcal{G}_t(x)\| \leq \|g(x)\|\) follows from~\cite[Theorem 3.5]{DL18} by noting that \(\|g(x)\| = {\rm dist}(0, \nabla f(x) + \lambda\partial \|x\|_1)\).

(b) Notice that, for any \(t > 0\),  
\begin{equation}\label{eq:dxx}
\mathcal{G}_t(x)_i = \left\{
\begin{array}{ll}
(\nabla f(x))_i + \lambda, & {\rm if}~x_i > \frac{1}{t}((\nabla f(x))_i + \lambda);\\
tx_i, & {\rm if}~\vert x_i - (\frac{1}{t}\nabla f(x))_i\vert \leq \frac{\lambda}{t};\\
(\nabla f(x))_i - \lambda, & {\rm if}~x_i < \frac{1}{t}((\nabla f(x))_i - \lambda). 
\end{array}
\right. \quad i = 1, \ldots, n.
\end{equation}
We denote \(\nabla_i f := (\nabla f(x))_i\) to simplify the notation. 

\noindent Case I. \(i \in I_0\). In this case, 
\[
g_i = \nabla_if - \min\{\max\{-\lambda, \nabla_if\}, \lambda\} = \left\{
\begin{array}{ll}
\nabla_if - \lambda, & {\rm if}~\nabla_if > \lambda;\\
0, & {\rm if}~\nabla_if \in[-\lambda, \lambda];\\
\nabla_if + \lambda, & {\rm if}~\nabla_if < -\lambda. 
\end{array}
\right.
\]
\begin{itemize}
\item[a)] If \(\nabla_i f > \lambda\), then \(x_i = 0 < \frac{1}{t}(\nabla_i f - \lambda)\), and hence \(\mathcal{G}_t(x)_i = g_i\);
\item[b)] If \(\nabla_i f \in [-\lambda, \lambda]\), then \(x_i = 0 \in[\frac{1}{t}(\lambda - \nabla_i f), \frac{1}{t}(\nabla_i f + \lambda)]\) and \(\mathcal{G}_t(x)_i = tx_i = 0 = g_i\);  
\item[c)] If \(\nabla_i f < -\lambda\), then \(x_i = 0 > \frac{1}{t}(\nabla_i f + \lambda)\) and hence \(\mathcal{G}_t(x)_i = g_i\). 
\end{itemize}

\noindent Case II. \(i \in I_+\). In this case, we have \(g_i = \nabla_i f + \lambda\). 
\begin{itemize}
\item[a)] If \(\nabla_i f > \lambda\), then, by~\eqref{eq:dxx}, we have 
\[
\mathcal{G}_t(x)_i = \left\{
\begin{array}{ll}
\nabla_if + \lambda = g_i, & {\rm if}~x_i > \frac{1}{t}(\nabla_if + \lambda);\\
tx_i \leq g_i, & {\rm if}~\vert x_i - \frac{1}{t}\nabla_if\vert \leq \frac{\lambda}{t};\\
\nabla_if - \lambda \leq g_i, & {\rm if}~0 < x_i < \frac{1}{t}(\nabla_if - \lambda);
\end{array}
\right. 
\]
\item[b)] If \(\nabla_i f \in [-\lambda, \lambda]\), then, by~\eqref{eq:dxx}, we have 
\[
\mathcal{G}_t(x)_i = \left\{
\begin{array}{ll}
\nabla_if + \lambda = g_i, & {\rm if}~x_i > \frac{1}{t}(\nabla_if + \lambda);\\
tx_i \leq g_i, & {\rm if}~0 < x_i \leq \frac{1}{t}(\nabla_i f + \lambda);
\end{array}
\right. 
\]
\item[c)] If \(\nabla_i f < -\lambda\), then \(\nabla_if + \lambda < 0 < x_i\). By~\eqref{eq:dxx}, we have \(\mathcal{G}_t(x)_i = g_i\). 
\end{itemize}
From the definition of \(\hat{t}\), we have \(x_i \geq \frac{1}{t}g_i = \frac{1}{t}(\nabla_i f + \lambda)\), \(\forall t \geq \hat{t}\). Recalling~\eqref{eq:dxx}, we have \(\mathcal{G}_t(x)_i = g_i\).

\noindent Case III. \(i\in I_-\). In this case, we have \(g_i = \nabla_i f - \lambda\). 
\begin{itemize}
\item[a)] If \(\nabla_i f > \lambda\), then \(\nabla_i f - \lambda > 0 > x_i\). By~\eqref{eq:dxx}, we have \(\mathcal{G}_t(x)_i = g_i\).
\item[b)] If \(\nabla_i f \in [-\lambda, \lambda]\), then, by~\eqref{eq:dxx},  
\[
\vert \mathcal{G}_t(x)_i\vert = \left\{
\begin{array}{ll}
t\vert x_i\vert  \leq  \vert g_i\vert, & {\rm if}~\frac{1}{t}g_i =\frac{1}{t}(\nabla_if - \lambda) \leq x_i < 0;\\
\vert g_i\vert, & {\rm if}~x_i < \frac{1}{t}(\nabla_i f - \lambda);
\end{array}
\right. 
\]
\item[c)] If \(\nabla_i f < -\lambda\), then, by~\eqref{eq:dxx}, 
\[
\mathcal{G}_t(x)_i = \left\{
\begin{array}{ll}
\vert \nabla_i f + \lambda\vert \leq \vert g_i\vert , & {\rm if}~0 > x_i > \frac{1}{t}(\nabla_if + \lambda);\\
t\vert x_i\vert \leq \vert g_i\vert, & {\rm if}~\vert x_i - \frac{1}{t}\nabla_if\vert \leq \frac{\lambda}{t};\\
\vert \nabla_i f - \lambda\vert = \vert g_i\vert, & {\rm if}~x_i < \frac{1}{t}(\nabla_if - \lambda) = \frac{1}{t}g_i; 
\end{array}
\right. 
\]
\end{itemize}
From the definition of \(\hat{t}\), we have \(x_i \leq \frac{1}{t}g_i = \frac{1}{t}(\nabla_if - \lambda)\), \(\forall t \geq \hat{t}\). Recalling~\eqref{eq:dxx}, we have \(\mathcal{G}_t(x)_i = g_i\).

(c) Let $I_{\neq 0} := \text{supp}(\bar{x})$ and $I_0 := \{i : \bar{x}_i = 0\}$. Define $x^+ = \text{prox}_{\frac{\lambda}{t}|\cdot|_1}(x - \frac{1}{t}\nabla f(x))$. We first show that there exists $r_1 > 0$ such that for all $x \in \mathbb{B}(\bar{x}, r_1)$, the following sign and support conditions hold:
\begin{equation}
\label{eq:sign_ident}
\text{sgn}(x_i) = \text{sgn}(x_i^+) = \text{sgn}(\bar{x}_i), \ \forall i \in I_{\neq 0}; \quad \text{and} \quad x_i^+ = 0, \ \forall i \in I_0.
\end{equation}
Specifically, under the SC condition, the proximal gradient mapping possesses the finite identification property. As established in \cite{L02, HL07}, there exists $\delta_1 > 0$ such that for all $x \in \mathbb{B}(\bar{x}, \delta_1)$, $x^+$ correctly identifies the active manifold, including the support and sign pattern, of $\bar{x}$. This implies $\text{sgn}(x_i^+) = \text{sgn}(\bar{x}_i)$ for $i \in I_{\neq 0}$ and $x_i^+ = 0$ for $i \in I_0$. Additionally, by the continuity of the coordinate projections, there exists $\delta_2 > 0$ such that, $\text{sgn}(x_i) = \text{sgn}(\bar{x}_i)$ for all $i \in I_{\neq 0}$ and $x \in \mathbb{B}(\bar{x}, \delta_2)$. Thus, \eqref{eq:sign_ident} holds by setting $r_1 = \min\{\delta_1, \delta_2\}$.

Now, consider $x \in \mathbb{B}(\bar{x}, r_1) \cap \mathcal{M}$. For any $i \in I_{\neq 0}$, the definition of the proximal operator yields
\begin{equation*}
    x_i^+ = x_i - \frac{1}{t}\left( (\nabla f(x))_i + \lambda \text{sgn}(x_i) \right).
\end{equation*}
Consequently, the $i$-th component of the proximal gradient is
\begin{equation*}
    (\mathcal{G}_{t}(x))_i = t(x_i - x_i^+) = (\nabla f(x))_i + \lambda \text{sgn}(x_i) = g_i(x).
\end{equation*}

For any $i \in I_0$, we have $x_i^+ = 0$, which implies $(\mathcal{G}_{t}(x))_i = t x_i$. Since $x \in \mathcal{M}$, it follows from the definition of \({\rm supp}(x)\) that $x_i = 0$ for all $i \in I_0$, leading to $(\mathcal{G}_{t}(x))_i = 0$. Furthermore, the SC property ensures that $|(\nabla f(\bar{x}))_i| < \lambda$, \(\forall i\in I_0\). By continuity, after shrinking \(r_1\) if necessary, we have $|(\nabla f(x))_i| < \lambda$ for all $x \in \mathbb{B}(\bar{x}, r_1)$. In this case, $g_i(x) = 0$. 

Thus, $(\mathcal{G}_{t}(x))_i = g_i(x)$ for all $i \in \{1, \ldots, n\}$, which completes the proof.
\end{proof}

\section{Minimum Eigenvalue Oracle (MEO)}\label{appendix:meo}

Algorithm~\ref{alg:meo}, originally proposed in~\cite{RNW20}, is a minimum-eigenvalue oracle for a symmetric matrix \(H\). Given a tolerance \(\epsilon > 0\) and a failure probability \(\sigma \in [0, 1)\), it either returns a unit vector \(v\) such that \(v^\top Hv \leq -\epsilon/2\) or certifies that \(\lambda_{\min}(H) \geq -\epsilon\). The latter certificate is correct with probability at least \(1 - \sigma\). The MEO can be implemented using randomized Lanczos or conjugate gradient (CG) iterations with a random starting vector~\cite{RNW20}.

\begin{algorithm}[h!]
\caption{Minimum Eigenvalue Oracle (MEO).} \label{alg:meo}
\begin{algorithmic}[1]
\State{Inputs: Symmetric matrix \(H\), tolerance \(\epsilon > 0\), error probability \(\sigma\in[0, 1)\);}
\State{Optional input: upper bound on Hessian norm \(M\);}
\State{Outputs: An estimate \(\hat{\lambda}\) of \(\lambda_{\min}(H)\) such that \(\hat{\lambda} \leq -\epsilon/2\), and vector v with \(\|v\| = 1\) such that \(v^\top Hv = \hat{\lambda}\) or a certificate that \(\lambda_{\min}(H) \geq - \epsilon\). The probability that the certificate is issued but \(\lambda_{\min}(H) < -\epsilon\) is at most \(\sigma\).}
\end{algorithmic}
\end{algorithm}

\section{Capped CG method}\label{appendix:ccg}

The main loop of the Capped CG method~\cite{RNW20} consists of classical CG iterations while detecting indefiniteness in the Hessian. 
Algorithm~\ref{alg:ccg} presents the Capped CG method for solving the regularized Newton equation~\eqref{eq:rnecopycopy}. It is a direct adaptation of the method in~\cite{RNW20};  the only differences are  the values of \(\bar{H}\) and \(\kappa^{cg}\) at line~\ref{line1} in Algorithm~\ref{alg:ccg}, which stems from the regularization term of Equation~\eqref{eq:rnecopycopy}. All other steps of the algorithm are identical to those in~\cite{RNW20}.

\begin{algorithm}[th!]
\caption{Capped CG method for the regularized Newton equation~\eqref{eq:rnecopycopy}.} \label{alg:ccg}
\begin{algorithmic}[1]
\Require \(x\gets x_k\), operator \(H\), \(g\), \(\epsilon\), \(\zeta\in(0, 1)\), \(\delta\), \(\bar{\tau}\) and \(M \geq 0\).
\Ensure \([d, d_{\rm type}]\);
\State \textbf{Set} \(\bar{H} := H + \bar{\tau}\|g\|^{\delta}I\),~\(\kappa^{cg} := \frac{M + \bar{\tau}\|g\|^{\delta}}{\epsilon}\),~\(\hat{\zeta} := \frac{\zeta}{3\kappa^{cg}}\),~\(\tau := \frac{\sqrt{\kappa^{cg}}}{\sqrt{\kappa^{cg}} + 1}\), \(T:= \frac{4(\kappa^{cg})^4}{(1-\sqrt{\tau})^2}\); \(y_0 \gets 0\), \(r_0 \gets g\), \(p_0 \gets -g\), \(j \gets 0\); \label{line1}
\State \textbf{if} \(\langle p_0, \bar{H}[p_0]\rangle < \epsilon\|p_0\|^2\)~\textbf{then}\\
\quad Set \(d = p_0\) and terminate with \(d_{\rm type} = {\rm NC}\);\\
\textbf{else if} \(\|H[p_0]\| > M\|p_0\|\)~\textbf{then}\\
\quad Set \(M \gets \|H[p_0]\|/\|p_0\|\) and update \(\kappa^{cg}, \hat{\zeta}, \tau, T\) accordingly;\\
\textbf{end if} 
\While{TRUE}
\State \(\alpha_j \gets \langle r_j, r_j\rangle/\langle p_j, \bar{H}[p_j]\rangle\); 
\State \(y_{j+1} \gets y_j + \alpha_jp_j\);
\State \(r_{j+1} \gets r_j + \alpha_j\bar{H}[p_j]\);
\State \(\beta_{j+1}\gets \langle r_{j+1}, r_{j+1}\rangle/\langle r_j, r_j\rangle\);
\State \(p_{j+1} \gets -r_{j+1} + \beta_{j+1}p_j\);
\State \(j\gets j+1\); 
\State \textbf{if} \(\|H[p_j]\| > M\|p_j\|\)~\textbf{then}\\
\quad\quad\quad Set \(M \gets \|H[p_j]\|/\|p_j\|\) and update \(\kappa^{cg}, \hat{\zeta}, \tau, T\) accordingly;\\
\quad\quad\textbf{end if}
\State \textbf{if} \(\|H[y_j]\| > M\|y_j\|\)~\textbf{then}\\
\quad\quad\quad Set \(M \gets \|H[y_j]\|/\|y_j\|\) and update \(\kappa^{cg}, \hat{\zeta}, \tau, T\) accordingly;\\
\quad\quad\textbf{end if}
\State \textbf{if} \(\|H[r_j]\| > M\|r_j\|\)~\textbf{then}\\
\quad\quad\quad Set \(M \gets \|H[r_j]\|/\|r_j\|\) and update \(\kappa^{cg}, \hat{\zeta}, \tau, T\) accordingly;\\
\quad\quad\textbf{end if}
\State \textbf{if} \(\langle y_j, \bar{H}[y_j]\rangle < \epsilon\|y_j\|^2\)~\textbf{then}\\
\quad\quad\quad Set \(d \gets y_j\) and terminate with \(d_{\rm type} = NC\);\\
\quad\quad\textbf{else if} \(\|r_j\| \leq \hat{\zeta}\|r_0\|\)~\textbf{then}\\
\quad\quad\quad Set \(d \gets y_j\) and terminate with \(d_{\rm type} = SOL\);\\
\quad\quad\textbf{else if} \(\langle p_j, \bar{H}[p_j]\rangle < \epsilon\|p_j\|^2\)~\textbf{then}\\
\quad\quad\quad Set \(d \gets p_j\) and terminate with \(d_{\rm type} = NC\);\\
\quad\quad\textbf{else if} \(\|r_j\| > \sqrt{T}\tau^{j/2}\|r_0\|\)~\textbf{then}\\
\quad\quad\quad Compute \(a_j, y_{j+1}\) as in the main loop above;\\
\quad\quad\quad Find \(i\in\{0, \cdots, j-1\}\) such that 
\(\frac{\langle y_{j+1} - y_i, \bar{H}[y_{j+1} - y_i]\rangle}{\|y_{j+1} - y_i\|^2} < \epsilon\);\\
\quad\quad\quad Set \(d \gets y_{j+1} - y_i\) and terminate with \(d_{\rm type} = NC\); \\
\quad\quad\textbf{end if}
\EndWhile
\end{algorithmic}
\end{algorithm}

\bibliographystyle{unsrtnat}
\bibliography{manuscript}

\end{document}